\documentclass[11pt]{article}
\usepackage{geometry}
\usepackage{xcolor}
\usepackage{xspace}

\usepackage{enumitem}
\setlist[itemize]{label=--}

\RequirePackage{amsthm,amsmath,amsfonts,amssymb}
\RequirePackage[numbers]{natbib}
\definecolor{LinkBlue}{HTML}{234B6C}
\definecolor{CiteTeal}{HTML}{17635A}
\definecolor{URLSlate}{HTML}{455A78}
\RequirePackage{hyperref}
\hypersetup{
  colorlinks=true,
  linkcolor=LinkBlue,
  citecolor=CiteTeal,
  urlcolor=URLSlate,
  filecolor=URLSlate
}
\RequirePackage{graphicx}
\usepackage{subcaption}
\usepackage{bbm}
\usepackage{algorithm2e}
\RestyleAlgo{ruled}
\usepackage{multicol}
\usepackage{mathtools}

\usepackage{accents}

\usepackage{textgreek}

\theoremstyle{plain}

\newtheorem{theorem}{Theorem}[section]
\newtheorem{lemma}[theorem]{Lemma}

\newtheorem{corollary}[theorem]{Corollary}

\theoremstyle{definition}
\newtheorem{definition}[theorem]{Definition}
\newtheorem{rem}[theorem]{Remark}

\def\beq{\begin{equation}} 
\def\eeq{\end{equation}}
\def\beqn{\begin{eqnarray*}}
\def\eeqn{\end{eqnarray*}}
\def\Bal{\begin{align}}
\def\Eal{\end{align}}
\def\Bitem{\begin{itemize}\setlength{\itemsep}{.2in}}
\def\bitem{\begin{itemize}\setlength{\itemsep}{.05in}}
\def\eitem{\end{itemize}}
\def\blatin{\begin{enumerate}\setlength{\itemsep}{.05in}\renewcommand{\labelenumi}{\roman{enumi}.}}
\def\elatin{\end{enumerate}}
\def\Benum{\begin{enumerate}\setlength{\itemsep}{.2in}}
\def\benum{\begin{enumerate}\setlength{\itemsep}{.05in}}
\def\eenum{\end{enumerate}}
\def\bmult{\begin{multline*}}
\def\emult{\end{multline*}}
\def\bcenter{\begin{center}}
\def\ecenter{\end{center}}
\def\bframe{\begin{frame}}
\def\eframe{\end{frame}}

\newcommand{\thmref}[1]{Theorem~\ref{thm:#1}}

\newcommand{\lemref}[1]{Lemma~\ref{lem:#1}}

\def\cB{\mathcal{B}}
\def\cC{\mathcal{C}}

\def\cE{\mathcal{E}}

\def\cG{\mathcal{G}}
\def\cH{\mathcal{H}}

\def\cM{\mathcal{M}}

\def\cP{\mathcal{P}}
\def\cQ{\mathcal{Q}}
\def\cR{\mathcal{R}}
\def\cS{\mathcal{S}}
\def\cT{\mathcal{T}}
\def\cU{\mathcal{U}}
\def\cV{\mathcal{V}}

\def\cX{\mathcal{X}}

\def\cZ{\mathcal{Z}}

\def\bx{\mathbf{x}}
\def\by{\mathbf{y}}

\def\bbE{\mathbb{E}}

\def\bbI{\mathbb{I}}

\def\bbM{\mathbb{M}}
\def\bbN{\mathbb{N}}

\def\bbP{\mathbb{P}}

\def\bbR{\mathbb{R}}

\def\bbT{\mathbb{T}}

\def\bbV{\mathbb{V}}

\def\ind{\mathbbm{1}}

\DeclareMathOperator{\diam}{diam}
\DeclareMathOperator{\dist}{dist}

\DeclareMathOperator{\op}{op} 
\DeclareMathOperator{\Card}{Card}

\DeclareMathOperator{\Id}{Id}
\DeclareMathOperator{\Vect}{Vect}

\DeclareMathOperator{\imag}{Im}
\renewcommand{\Im}{\imag}
\DeclareMathOperator{\II}{I\!I}

\DeclareMathOperator{\pr}{\pi}

\DeclareMathOperator{\supp}{support}
\DeclareMathOperator{\vol}{vol}

\DeclareMathOperator{\rad}{rad}

\DeclareMathOperator{\spec}{spec}

\DeclareMathOperator{\conv}{conv}
\DeclareMathOperator{\Conv}{Conv}

\let\lac\{
\let\rac\}
\renewcommand{\{}{\left\lac}
\renewcommand{\}}{\right\rac}
\newcommand{\inner}[2]{\langle #1, #2 \rangle}

\def\({\left(}
\def\){\right)}

\newcommand{\ve}{\varepsilon}

\newcommand*\diff{\mathop{}\!\mathrm{d}}
\newcommand\wt{\widetilde}
\newcommand\wh{\widehat}
\renewcommand{\leq}{\leqslant}
\renewcommand{\geq}{\geqslant}
\newcommand*\ball{\mathop{}\mathsf{B}}

\newcommand{\interior}[1]{\mathring{#1}} 
\newcommand*\oball{\interior{\ball}}

\DeclareMathOperator{\dt}{\mathrm{d}}

\newcommand{\modelM}{\cM^{(d)}(\kappa_{\max},V_{\max})}

\makeatletter
\newcommand{\para}{\mathrel{\mathpalette\new@parallel\relax}}
\newcommand{\new@parallel}[2]{%
  \begingroup
  \sbox\z@{$#1T$}
  \resizebox{!}{\ht\z@}{\raisebox{\depth}{$\m@th#1/\mkern-5mu/$}}%
  \endgroup
}
\makeatother

\newcommand{\dH}{\mathrm{d}_{\mathrm{H}}}
\DeclareMathOperator{\Span}{span}
\renewcommand{\parallel}{\mathrel{/\mkern-5mu/}}

\DeclareMathOperator{\slab}{\mathrm{S}}
\DeclareMathOperator{\VC}{VC}

\newcommand{\codetection}{\hyperref[algo:strat]{\texttt{Slabeling}}\xspace}

\makeatletter
\newcommand{\bauthor}[1]{{%
  \let\binits\@firstofone
  \let\bsnm\@firstofone
  \let\bfnm\@gobble
  \let\bparticle\@firstofone
  \let\bsuffix\@firstofone
  \scshape #1}}
\newcommand{\beditor}[1]{{%
  \let\binits\@firstofone
  \let\bsnm\@firstofone
  \let\bfnm\@gobble
  \let\bparticle\@firstofone
  \let\bsuffix\@firstofone
  \biblioeditorstyle #1}}
\let\biblioeditorstyle\scshape

\let\byear\@firstofone
\let\bpages\@firstofone
\let\bpublisher\@firstofone
\let\baddress\@firstofone

\def\biblio@withoutMR#1{\biblio@stripMR#1MR#1MR\relax}
\def\biblio@stripMR#1MR#2MR#3\relax{#2}
\def\biblio@MR#1 #2\relax#3{%
  \href{https://mathscinet.ams.org/mathscinet-getitem?mr=#1}{MR\biblio@withoutMR{#3}}}
\newcommand{\bmrnumber}[1]{\ \biblio@MR#1 \relax{#1}\ignorespaces}
\def\biblio@checkdoi#1#2#3#4#5\end{%
  \def\biblio@doifull{0}%
  \def\biblio@http{http}%
  \def\biblio@start{#1#2#3#4}%
  \ifx\biblio@start\biblio@http
    \def\biblio@doifull{1}%
  \fi}
\newcommand{\bdoi}[1]{\ 
  \biblio@checkdoi#1\end
  \ifnum\biblio@doifull>0\relax
    \href{#1}{\nolinkurl{#1}}%
  \else
    \href{https://doi.org/#1}{\nolinkurl{https://doi.org/#1}}%
  \fi
  \ignorespaces}
\AtBeginDocument{\urlstyle{rm}\setlength{\bibsep}{0pt}}
\makeatother

\begin{document}

\thispagestyle{empty}

\title{A theory of stratification learning\\
}
\author{
Eddie Aamari\thanks{D\'epartement de Math\'ematiques et Applications, \'Ecole Normale Sup\'erieure, Univ. PSL, CNRS, Paris, France. \texttt{eddie.aamari@ens.fr}}
\and
Clément Berenfeld\thanks{PreMeDICaL, Inria, Inserm, Univ. Montpellier, Montpellier, France. \texttt{clement.berenfeld@inria.fr}}
}
\date{}
\maketitle

\begin{abstract}
Given i.i.d. sample from a stratified mixture of immersed manifolds of different dimensions, we study the minimax estimation of the underlying stratified structure.
We provide a constructive algorithm allowing to estimate each mixture component at its optimal dimension-specific rate adaptively.
The method is based on an ascending hierarchical co-detection of points belonging to different layers, which also identifies the number of layers and their dimensions, assigns each data point to a layer accurately, and estimates tangent spaces optimally.
These results hold regardless of any ambient assumption on the manifolds or on their intersection configurations.
They open the way to a broad clustering framework, where each mixture component models a cluster emanating from a specific nonlinear correlation phenomenon.
\end{abstract}

\section{Introduction}

\subsection{Context}

\paragraph*{Union of manifold hypothesis}

With the ubiquitous availability of high-dimensional and unstructured data, the so-called \textit{manifold hypothesis} is often raised as an explanation for the surprisingly good performances of prediction~\cite{bengio2013representation} and generative methods~\cite{croitoru2023diffusion}.
It posits that high-dimensional data typically reside close to some~(unknown) submanifold $M$ of low intrinsic dimension~$d$, significantly smaller than the ambient dimension~$D$.
The support $M$ models all the possibly non-linear constraints or correlations that data may present, and $d$ represents its \textit{true} local degrees of freedom.
Statistically speaking, the success of modern machine learning methods is then attributed to estimation rates depending mostly (if not only) on this lower intrinsic dimension $d$, hence enabling effective learning on reasonably-sized datasets $X_1,\ldots,X_n$ even when $D \gg n$.
This approach has proven fruitful in a quite wide range of applications, including cancer classification~\cite{aukerman2022persistent}, RNA sequencing~\cite{moon2018manifold}, image analysis~\cite{carlsson2008local}, location-based services~\cite{jia2018dimension} and thermodynamics~\cite{ruppeiner1995riemannian}.

The geometric nature of the manifold hypothesis brings an entire field of possible data modelings based on differential geometry~\cite{doCarmo92} and geometric measure theory~\cite{Federer59}.
Manifold estimation being fundamentally non-parametric, one can see each geometric modeling (or statistical model) as a specific regularity constraint on the problem.
The most crude (and studied) one consists in assuming that $M$ is a connected $\mathcal{C}^2$-submanifold without boundary~\cite{Genovese12b}.
Most existing manifold estimation techniques assume that $M$ is a single well-behaved such manifold. 

However, recent empirical studies suggest that real data actually exhibit mixed dimensionalities~\cite{brown2023verifying,goldt2020modeling,abdolali2021beyond,li2022neural,pope2021intrinsic,robinson2024structure, ding2023unsupervised,lim2023hades}. 
The present work addresses this latter general scenario, where data actually arise from a stratification of manifolds of possibly different dimensions.
In probabilistic terms, we consider the estimation and clustering problem of arbitrary mixtures of distributions with individual components having $\cC^2$ submanifolds as a support.

\paragraph*{A manifold approach to generalized clustering}

Clustering (or unsupervised learning) is a statistical task that aims at grouping data points into \textit{simple} groups (called clusters) that are as \textit{separated} as possible from each other.
Cluster separation is most commonly quantified via the distance between clusters, while simplicity shall be measured through some notion of smallness or connectedness. 
However, achieving separation based on location only may not always be relevant, especially in high dimension.
Following~\cite{Jegelka09}, the \textit{generalized clustering} problem consists in decomposing the distribution $P$ of data into ``simple" components $(P_k)_{1 \leq k \leq K}$ while maximizing a specified distance function between them.

\begin{figure}[!b]
\centering
\includegraphics[page=5,width=\textwidth]{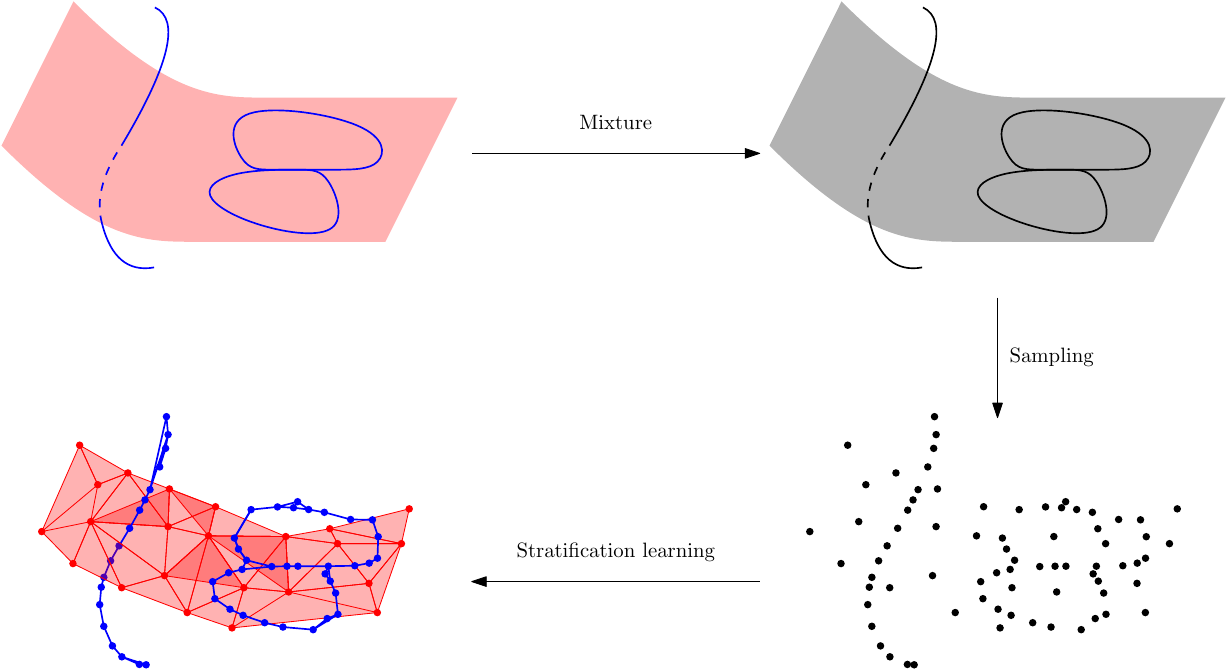}
\caption{
A union of curves $M_1$ and a surface $M_2$ in dimension $D=3$ intersecting arbitrarily (top right). None of them have boundary, and drawings are to be continued outside the picture if necessary.
The two completely visible closed curves intersect tangentially, coincide over a non-trivial interval, and are both completely included in the surface.
Despite the topology of a pre-image of $M_1 \subset \bbR^D$ by an immersion being not uniquely defined when $M_1$ self-intersects tangentially, the dimension-wise \textit{clustering} and \textit{reconstruction} of the overall support still makes sense (top left).
After sampling (bottom right), learning the stratification consists in labeling points dimension-wise, as well as approximating the underlying manifold structures (bottom left).
}
\label{fig:cluster-from-manifolds}

\end{figure}

As mentioned above, local manifold parametrizations model local data correlations arising from numerous non-linear constraints of the generating distribution.
Hence, membership in different parametrizations yields a natural geometric criterion for cluster separation.
The sample counterpart is a fine non-linear clustering of a point cloud that both includes positional and directional information given by tangent spaces locally~\cite{Wu_2022_WACV}.
See Figure~\ref{fig:cluster-from-manifolds} for an illustration.

\subsection{Contribution}
\label{sec:contribution}

\subsubsection*{Informal results}

This work presents a first attempt to build an estimation theory of general unions of compact $\cC^2$-submanifolds of $\bbR^D$.
Its framework steps away from ambient convexity-type assumptions on the target set such as the reach, $\mu$-reach or rolling ball conditions~\cite{Chazal06,rodriguez2007set}.
It shows that
bounded volume and directional curvature are sufficient to obtain fast convergence rates as well as dimension-wise clustering.

Let
$X_1,\ldots,X_n \in \bbR^D$ be drawn i.i.d. from a union $\bigcup_{k=1}^K M_k$ of immersed compact $\cC^2$-manifolds, with unknown number of mixture components $K$ and unknown dimensions $d_1 < \ldots < d_K$ (see Figure~\ref{fig:cluster-from-manifolds}).
With $\cX_n := \{X_1,\ldots,X_n\}$ and scale parameters as input, we propose a constructive algorithm called \codetection, which outputs an approximation of the stratified structure of $M$.
This output has a layer-by-layer structure, each of which corresponds to a fixed dimension.
It allows us to recover the structure $(M_k)_{1 \leq k \leq K}$ in the following sense.

\smallskip
\noindent
\textit{(Discrete stratification learning)}
With high probability, the discrete structure of $M = \bigcup_{k=1}^K M_k$ can be identified by \codetection through
\begin{align}
\hat{K} = K \text{~and~}\{\hat{d}_1,\ldots,\hat{d}_K\} = \{d_1,\ldots,d_K\}.
\tag{Theorem~\ref{thm:dimension-and-component-estimation}}
\end{align}
Furthermore, data points $\cX_{k,n} := \ M_k \cap \cX_n$ from each layer $M_k$ can be clustered by dimension through a constructive decomposition $\cX_n = \bigsqcup_{1 \leq k \leq K} \wh\cX_{k,n}$ of the data set such that
\begin{align*}
\frac{\#\bigl(\wh\cX_{k,n}  \triangle \cX_{k,n}\bigr)}{\max\{ \# \cX_{k,n},1\}} 
\lesssim 
\(\frac{\log n}{n}\)^{2/d_{k}} 
,
\tag{Theorem~\ref{thm:main-clustering-dimensions}}
\end{align*}
with accurate dimension labeling (Theorem~\ref{thm:main-dimension-labeling}).

\smallskip
\noindent
\textit{(Approximation)}
On the same event of high probability,
all the manifolds $M_k$ can be estimated in Hausdorff distance at the layer-adaptive rate
\begin{align}
\dH(M_k,\wh{M}_k)
\lesssim
\left(\frac{\log n}{n}\right)^{2/d_k}
,
\tag{Theorem~\ref{thm:main-manifold-estimation}}
\end{align}
and tangent spaces at $X_i$'s can be estimated,
in the sense that points $x \in M_k$ and $y \in \wh M_k$ at distance $\| x -y \| \lesssim (\log n / n)^{2/d_k}$ can be put in correspondence such that
\begin{align}
\angle\bigl(
T_x M_k
,
T_y \wh M_k
\bigr)
\lesssim
\left(\frac{\log n}{n}\right)^{1/d_k}
\tag{Theorem~\ref{thm:main-tangent-estimation}}
\end{align}
{\color{black}The estimate uses directional information already computed by \codetection. These directions help reject misaligned tuples and control the geometric error of accepted cross-branch simplices.}

The algorithm does not require any prior information about the number $K$ of mixture components, but only a maximal intrinsic dimension $d_{\max}$ past which points shall be considered coming from a too-high dimension and hence possibly considered clutter noise (see Remark~\ref{rem:clutter-as-ambient-manifold}). Naturally, the above results hold if scale parameters are properly chosen and $d_{\max} \geq d_K$.
Otherwise, the results extend to the learning of the thresholded structure 
$\bigcup_{\substack{1 \leq k \leq K \\ d_k \leq d_{\max} }} M_k$.

\paragraph*{Consequences}

The main high-level consequence of our results is that for manifold estimation, intrinsic geometric regularity constraints (namely curvature and volume bounds) are sufficient to obtain fast estimation rates. That is, no assumption on the reach or any other ambient regularity parameters is required.
Possible (self-)intersections and interplay between layers of different dimensions are also harmless for the final rate.
This phenomenon also applies to tangent space estimation, but is in sharp contrast to other geometric and topological features such as homology or geodesic distances (see Section~\ref{sec:estimation-impossibilities}), which become impossible to estimate without reach condition.

Going from single-source to mixtures, manifolds turn out to be much easier to handle than classical density-based (Gaussian) mixture models where source separation requirements are needed to maintain distinguishability~\cite{laurent2018multidimensional}.
Our algorithm resolves and separates all the layers up to dimension $d$ in time $O(n (\log n)^d)$, that is the same order of magnitude as dealing with the $d$-dimensional layer taken separately.

\subsection{Related works}

\paragraph*{Single-manifold estimation and reach condition}

The minimax theory of support estimation from data drawn on \textit{smooth} sets is now well understood.
In full dimension, all the studied regularity constraints build upon some notion of convexity. 
This includes standard convexity~\cite{dumbgen1996rates}, $r$-convexity~\cite{aaron2016local}, rolling ball conditions~\cite{cuevas2012statistical} and reach positiveness~\cite{Federer59}, all yielding rates of order $O\bigl((\log n / n)^{2/D}\bigr)$ or $O\bigl((\log n / n)^{2/(D+1)}\bigr)$ in Hausdorff distance.

For embedded structures of dimension $d<D$, regularity constraints always borrow from \textit{both} ambient convexity-type constraints \textit{and} differential geometry simultaneously.
Usually, the ambient constraint is formalized through a positive lower bound on the reach.
Intrinsic differential smoothness takes the form of $\cC^k$ parametrizations for $k\geq 2$~\cite{Genovese12b,Aamari18,Aamari19,Divol21} --- yielding rate $O\bigl((\log n / n)^{k/d}\bigr)$ --- with possibly non-empty boundary for $k=2$~\cite{Aamari21} --- yielding rate $O\bigl((\log n / n)^{2/(d+1)}\bigr)$.
These methods rely on local polynomials or local convex hulls depending on the context.

Seen as a regularity parameter, the reach prevents quasi-self intersections to be arbitrarily narrow~\cite{Aamari19b}.
It has been used as a key scaling factor in all the works cited above, besides (and sometimes instead of) curvature bounds.
Beyond the technical easiness induced by reach constraints~\cite{Federer59}, the systematic assumption of reach positiveness probably is to be attributed to the fact that some geometric and topological quantities become impossible to estimate under intrinsic constraints only (see Section~\ref{sec:estimation-impossibilities}). 
Yet, as proved in the present article, several key geometric quantities remain estimable without convexity-type constraints and even for general unions of immersed manifolds, with no deterioration in the estimation rates.
This includes, among others, the support and tangent spaces.

\paragraph*{A concise review of the multi-manifold literature}

To our knowledge, the only works coming with provable geometric or statistical guarantees require either linearity, or mixture component to intersect transversally.
The others are of heuristic or numerical nature.
Theses pieces of work can be grouped as follows.
\begin{itemize}

\item
\textit{(Multiple PCA)}
The cluster recovery of mixture components formed by low-dimensional objects in high dimension has first been carried out in the case of linear subspaces \cite{lerman2010probabilistic,thomas2014learning}.
The main method relies on generalized Principal Component Analysis (PCA)~\cite{vidal2005generalized}, which consists in the minimization of a tweaked least square functional.

\item
\textit{(Local PCA)}
When coping with non-linearity, the most studied methods are based on local PCA~\cite{arias2011spectral,tinarrage2023recovering,lim2023hades} and kernelized version of it~\cite{chen2009spectral}.
The rationale is that properly-localized covariance matrices estimate local tangent directions through their largest singular values and associated singular vectors.
Conversely, if the dimension is known to be $d$, a covariance with at least $d+1$ singular values significantly larger than the others indicates the presence of a local $(d+1)$-dimensional space of variations, and hence a branching (or intersection) point.
Somehow related is the kernel-based clustering method of~\cite{Jegelka09}, which builds upon a preliminary featurization step onto a Hilbert space, intended to linearize the geometric setting globally.

\item
\textit{(Homology-based clustering)}  
By definition, the geometry of neighborhoods at an intersection point is not that of a single $d$-dimensional ball, but that of a union of several $d$-dimensional balls.
This can be characterized topologically by the $(d-1)$-dimensional homology of annuli being that of a sphere or not.
\newline
This idea has been exploited in dimension $d=1$, where the union of closed curves is often called \textit{metric graph}~\cite{aanjaneya2011metric,bendich2012local, lecci2014statistical}. The method uses persistent homology, which compares topology at different scales.
Although the method manages to handle noise, the curved graph edges are assumed to intersect transversally at vertices.
We will show that in fact, no transversality assumption is required to estimate a union of curves with fast rate, as long as each of them is~$\cC^2$ smooth.

\item
\textit{(Dynamic methods)}
Instead of finely clustering raw data points via local geometric or topological descriptors, the \textit{mean-shift} algorithm~\cite{fukunaga1975estimation} moves them incrementally towards high density regions and then uses a simple $k$-mean clustering method.
In a mixture of manifolds setting,~
\cite{Wu_2022_WACV} proposed to perform an optimization scheme that alternates between local structure preservation and cluster distinguishability to account for tangential information and not only positions. 

\end{itemize}
Stepping slightly aside from the manifold-based hypotheses, low-dimensional structures can also be modeled in a relaxed way.
Data can be assumed to be generated \textit{near} a submanifold $M$, i.e. drawn according to a distribution having a very peaky ambient density in a neighborhood of $M$. 
In such approximate models, the existence of a likelihood has been leveraged in at least two orthogonal directions.
\begin{itemize}

\item 
\textit{(Density ridges)}
The existence of a density allows to define the so-called \textit{density ridges}.
They form a stratification of the ambient space through generalizations of density modes, called ridges, and indexed by dimension. 
Their estimation is now well understood, both theoretically and algorithmically~\cite{qiao2023confidence,MR4255222}. 
Although this line of work parallels quite closely the notion of mixture of manifolds, note that density ridges --- as name suggests --- only make sense if data is generated by a distribution having a density with respect to the ambient Lebesgue measure.
In the same vein, let us also mention poissonized non-linear likelihood models~\cite{haro2006stratification}, allowing for likelihood maximization.

\item
\textit{(Bayesian approaches)}
Putting a mixture prior on the density leads to estimators that are structurally stable to mixtures.
As a result, the super-level sets of the posterior density behaves well when distributions are mixed.
In the linear case, (nearly-)singular Gaussian mixtures have been studied in~\cite{thomas2014learning} on union of (thickened) subspaces.
In the non-parametric curved case, local Gaussian mixtures have been put to use in~\cite{berenfeld2022estimating}.

\end{itemize}
Let us insist on the fact that all the above cited articles either require transversality at intersection points, or do not come with a statistical analysis.
Furthermore, we are not aware of any work that goes beyond the clustering-by-dimensionality problem (or dimension labeling), and that would allow for a geometric reconstruction of the mixture components.
\section{Geometric and statistical model}

\subsection{Euclidean immersions, curvature and volume}
\label{sec-euclidean-immersions}

We recall that a Euclidean $\cC^2$-immersion of a compact differentiable manifold $\mathbb{M}_0$ is a twice continuously differentiable map $\Psi_0 : \mathbb{M}_0 \to \bbR^D$ whose differential is everywhere injective.
The image set $M := \Psi_0(\mathbb{M}_0) \subset \bbR^D$ is called an \textit{immersed manifold}.

The differential structure of $\bbM_0$ at any point $x_0 \in \bbM_0$ is canonically pushed forward via the differential of $\Psi_0$.
We write $T_{x_0}^{(\Psi_0)} M := \diff \Psi_0(x_0)[T_{x_0} \bbM_0]$, where $T_{x_0} \bbM_0$ stands for the tangent space of $\bbM_0$ at $x_0$.
Accordingly, the \emph{tangent cone} at a point $x \in M = \Psi_0(\bbM_0)$ is the set
\begin{align*}
T_x^{(\Psi_0)} M
:=
\bigcup_{x_0 \in \Psi_0^{-1}(x)}
\diff \Psi_0(x_0)[T_{x_0} \bbM_0] 
\end{align*}
of pushed forward tangents of the pre-image $\Psi_0^{-1}(x)$.
Now looking at structure of order two, let us recall that the Euclidean metric on $\bbR^D$ pulls back canonically onto $\bbM_0$ via $\Psi_0$, so that $\bbM_0$ inherits a Riemannian structure $(\bbM_0,g_0)$ making $\Psi_0 : (\bbM_0,g_0) \to (\bbR^D,\langle \cdot, \cdot\rangle)$ an isometric immersion in the Riemannian sense.
This allows to define the \textit{second fundamental} form~\cite[Chapter~6]{doCarmo92}
$
\II_{x_0}^{(\Psi_0)} : T_{x_0} \bbM_0 \times T_{x_0} \bbM_0 \mapsto T_{x_0}^{(\Psi_0)} M^{\perp}
.
$
We write
$$
\kappa^{(\Psi_0)}(M) := \sup_{x_0 \in \bbM_0} \|\II_{x_0}^{(\Psi_0)}\|_{\op}
$$
for the maximal curvature of the immersion of $\bbM_0$, and $\vol(\bbM_0)$ its Riemannian volume.

\begin{figure}[!b]
        \centering
        \includegraphics[page=2,width=0.9\linewidth]{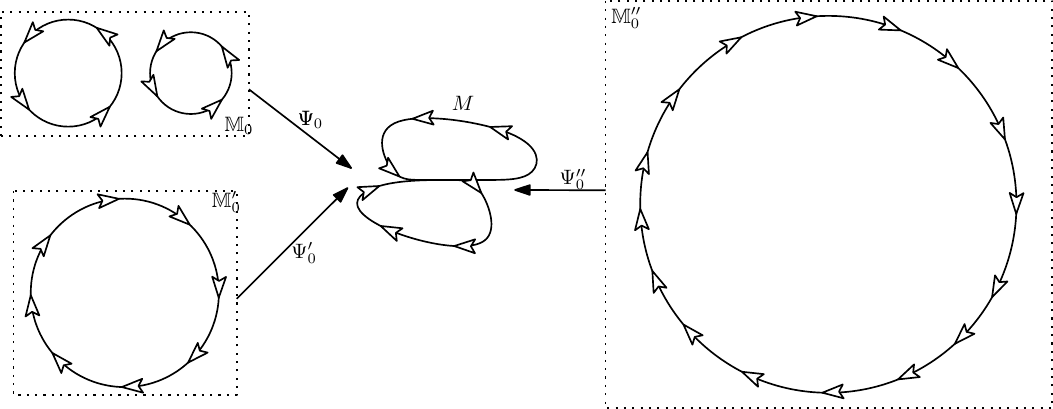}
\caption{
The immersed $8$-shape curve $M \subset \bbR^2$ (middle) can originate from either one or two topological circles, shown as $\mathbb{M}_0$ and $\mathbb{M}'_0$ (left).
Furthermore, running through the curve $M$ twice also yields an immersion with the same image, but from a manifold $\mathbb{M}''_0$ having twice the common volume of $\mathbb{M}_0$ and $\mathbb{M}'_0$ (right).
}
\label{fig:topology-not-identifiable}
\end{figure}

In what follows, we will infer properties of some (a priori unstructured) geometric set $M \subset \bbR^D$ under regularity assumptions induced by the \emph{existence} of some $\cC^2$ immersion $\Psi_0 : \bbM_0 \to M$. 
Hence, $\bbM_0$, its Riemannian and differential structures shall be seen as a latent unobserved structure.
However, there may not be a unique manifold $\bbM_0$ and immersion $\Psi_0 : \bbM_0 \to \bbR^D$ with geometric image $M$.
Actually, neither the geometry nor the topology of $(\bbM_0,g_0)$ is uniquely defined (see Figure~\ref{fig:topology-not-identifiable}).

If $M \subset \bbR^D$,
we let $\bbI(M)$ denote the set of couples $(\bbM_0,\Psi_0)$ such that $\bbM_0$ is an abstract $\cC^2$ manifold, and $\Psi_0 : \bbM_0 \to \bbR^D$ is a $\cC^2$ immersion with image $\Psi_0(\bbM_0) = M$.
Despite $\bbI(M)$ possibly containing several very different geometric structures, the following result asserts that their pushforward tangent cones and maximal curvature actually are all equal.

\begin{lemma}\label{lem:maximal-curvature-identification}
Given $M \subset \bbR^D$ with $\bbI(M) \neq \emptyset$ and $x \in M$, then for all immersions $(\bbM_0,\Psi_0)$ and $(\bbM_0',\Psi_0')$ in $\bbI(M)$, there holds 
\begin{center}
$T_x^{(\Psi_0)} M = T_x^{(\Psi_0')} M$
\hspace{1em} and \hspace{1em}
$\kappa^{(\Psi_0)}(M) = \kappa^{(\Psi_0')}(M)$.
\end{center}
\end{lemma}
The proof can be found in Section~\ref{sec:appendix-curvature}.
Lemma~\ref{lem:maximal-curvature-identification} allows to define unambiguously the \emph{tangent cones of $M$ at $x \in M$ as}
\begin{align}
\label{eq:tangents}
T_x M := T^{(\Psi_0)}_x M
,
\end{align}
and the \emph{maximal curvature} of $M$ as
\begin{align}
\label{eq:max-curvature}
\kappa(M) := \kappa^{(\Psi_0)}(M)
,
\end{align}
which do not depend on the choice of an immersion of $(\bbM_0,\Psi_0) \in \bbI(M)$.
On the other hand, the volume $\vol(\bbM_0)$ obviously depends on the immersion (see Figure~\ref{fig:topology-not-identifiable}). 
We define the \emph{minimal immersion-realizable volume} of $M \subset \bbR^D$ as
\begin{align}
\label{eq:min-volume}
\bbV(M) := \inf_{(\bbM_0,\Psi_0) \in \bbI(M)} \vol(\bbM_0)
.
\end{align}
By construction, notice in particular that $\vol(M) \leq \bbV(M)$, where $\vol(M)$ stands for the $d$-dimensional Hausdorff measure of $M$ seen as a rectifiable subset of $\bbR^D$.
We are now in position to define the volume- and curvature-constrained model of $\cC^2$ immersed manifolds at a given dimension.

\begin{figure}[!b]
\centering
\includegraphics[width=0.8\linewidth]{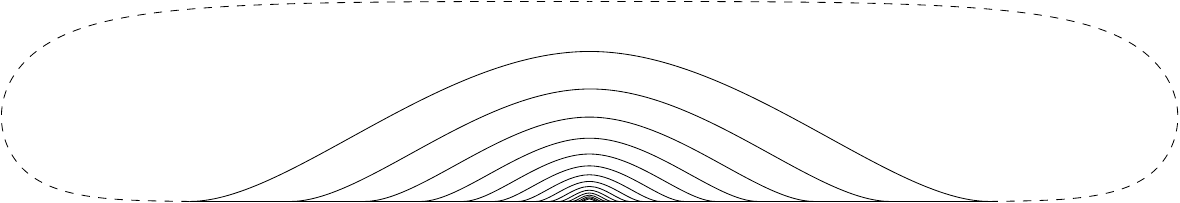}
\caption{
An increasing family of immersed curves $(M_\ell)_{\ell\in \bbN}$, each with $\ell$ nested $\cC^2$ bumps of exponentially decaying height and width.
The curves have uniformly bounded lengths $\vol(M_\ell) \leq \vol(\bigcup_{\ell' \in \bbN} M_{\ell'}) < \infty$.
However, the nestedness of the bumps forces immersions $\Psi_\ell : \bbM_\ell \to M_\ell$ describing $\ell$ bumps to go around the upper unbumped part $M^{(\mathrm{u})}$ (dashed) at least $\ell$ times, so that $\vol(\bbM_\ell) \geq \ell \vol(M^{(\mathrm{u})}) \xrightarrow[\ell \to \infty]{} \infty$.
}
\label{fig:fractal-mayhem}
\end{figure}

\begin{definition}[Geometric model for immersed $\cC^2$ manifolds]
\label{def:model-immersion}
For $d \in \{1,\ldots,D-1\}$, $\kappa_{\max} > 0$ and $V_{\max} > 0$, we let $\cM^{(d)}(\kappa_{\max},V_{\max})$ denote the set of compact immersed $\cC^2$ manifolds $M$ of $\bbR^D$ of dimension $d$ and such that 
\bitem
\item[i)] $\kappa(M) \leq \kappa_{\max}$;
\item[ii)] $\bbV(M) \leq V_{\max}$.
\eitem

\end{definition}
Constraining $\vol(M)$ rather that $\bbV(M)$ may first seem more natural.
In fact, bounding $\vol(M)$ only would still allow for pathological subsets $M \subset \bbR^D$ that exhibit points $x \in M$ with arbitrarily wild local neighborhood fibers $\Psi_0^{-1}(\ball(x,r))$ (see Figure~\ref{fig:fractal-mayhem}).
In contrast, a (stronger) bound on $\bbV(M)$ conveys enough rigidity to assert that neighborhoods $M \cap \ball(x,r)$ all resemble a quantitatively bounded union of Euclidean $d$-balls (see Lemma~\ref{lem:sheaf_number}).

\subsection{Manifold mixtures}
Let $M \in \cM^{(d)}(\kappa_{\max},V_{\max})$ have an associated $\cC^2$ immersion $(\bbM_0,\Psi_0)$. 
As a rectifiable subset of $\bbR^D$ parametrized by a compact $\cC^2$ manifold, $M$ has finite non-zero $d$-dimensional Hausdorff measure $\cH^d(M) = \vol(M)$~\cite[Section~3.2]{Federer69}.
We let $\vol_M$ denote the \emph{volume measure} on $M$, defined as 
\begin{align*}
\vol_M(\diff x) := \ind_{M} \cH^d(\diff x).
\end{align*}
In general, the volume measure $\vol_M$ on $M$ and the pushforward $\Psi_0^{\#} \vol_{\bbM_0}$ of the Riemannian volume measure on $\bbM_0$ for some $(\bbM_0,\Psi_0) \in \bbI(M)$ do \emph{not} coincide.
For example in Figure~\ref{fig:topology-not-identifiable}, one can take $\Psi_0$ to describe the upper circle twice and the other only once, leaving the latter with half the mass relatively to the upper circle.
However, these two measures do coincide locally, in the sense that for subsets $A \subset \bbM_0$ of small enough diameter, $\vol_{\bbM_0}(A) = \vol_M(\Psi_0(A))$.

In what follows, we will consider distributions over $M$ that are essentially uniform, in the sense that they have densities with respect to $\vol_M$ bounded away from $0$ and $\infty$.

\begin{definition}[Single-manifold model]
\label{def:model-single-component}
For $d \in \{1,\ldots,D-1\}$, $0 < \nu_{\max} < \infty$ and $0 < a_{\min} \leq a_{\max} < \infty$, we let $\cP^{(d)}_{\kappa_{\max},\nu_{\max}}(a_{\min},a_{\max})$ denote the set of Borel probability distributions over $\bbR^D$ such that
\begin{itemize}
\item[i)]
$M := \supp (P)$ belongs to $\cM^{(d)}(\kappa_{\max},\nu_{\max}(\kappa_{\max})^{-d} )$;
\item[ii)]
$P$ admits a density $f : M \to \bbR_+$ with respect to $\vol_M$ such that for all $x\in M$,
$$
a_{\min} \kappa_{\max}^{d} \leq f(x) \leq a_{\max} \kappa_{\max}^{d}
.
$$
\end{itemize}
\end{definition}

\begin{rem}[On model parameters]
\label{rem:model-parameters}
At this stage, let us emphasize two points.
\begin{itemize}
\item
The above model could as well be parametrized by $f_{\min} := a_{\min} \kappa_{\max}^{d}$, $f_{\max} = a_{\max} \kappa_{\max}^{d}$ and $V_{\max} = \nu_{\max}(\kappa_{\max})^{-d}$.
We choose instead to use a single scale parameter $\kappa_{\max}$ and nondimensionalized constants $a_{\min}, a_{\max},\nu_{\max}$ so that they can be shared across several models $\cP^{(d)}_{\kappa_{\max},\nu_{\max}}(a_{\min},a_{\max})$ with different values of $d$.
\item
Closed manifolds with curvature bounded by $\kappa_{\max}$ have volume bounded below by that of a sphere of radius $1/\kappa_{\max}$~\cite{almgren1986optimal}. 
Hence, one needs to request that $\nu_{\max} \geq \sigma_d$ to have a non-empty model in dimension $d$, where $\sigma_d$ is the volume of the $d$-dimensional sphere. 
In order to make sure that all our models are non-empty for $d \in \{1,\ldots,D-1\}$, we will hence assume without loss of generality that
\beq
\nu_{\max} \geq \max_{d \geq 1} \sigma_d = \sigma_6 = 16\pi^{3}/15, 
\eeq 
for the rest of the paper.

\end{itemize}
\end{rem}
This model generalizes that of \cite{Genovese12b,Aamari18,Divol21}, in which the support $M$ is additionally required to be an \textit{embedded} submanifold with positive reach, instead of an \textit{immersed} one with no reach constraint.
In addition, the present article considers general mixtures of such distributions with varying intrinsic dimensions~$d$.

\begin{definition}[Manifold mixture model]
\label{def:model-mixture}
Given $\alpha_{\min} > 0$, let $\bar{\cP}(\alpha_{\min},\kappa_{\max},\nu_{\max},\allowbreak a_{\min}, a_{\max})$ denote the set of mixture probability distributions $P$ over $\bbR^D$ of the form
\begin{align*}
P
=
\sum_{k=1}^K \alpha_k P_k
,
\end{align*}
where $K \in \{1,\ldots,D-1\}$ varies, and
\begin{itemize}
\item[i)] $1 \leq d_1 < d_2 < \ldots < d_K \leq D-1$,
\item[ii)] $\min_{1 \leq k \leq K} \alpha_k \geq \alpha_{\min}>0$,
\item[iii)]
$P_k \in  \cP_{\kappa_{\max},\nu_{\max}}^{(d_k)}(a_{\min}, a_{\max})$ for all $k \in \{1,\ldots,K\}$.
\end{itemize}
Model $\bar{\cP}$ is parametrized by $\alpha_{\min}, \kappa_{\max}, a_{\min}, a_{\max}$ and $\nu_{\max}$, all left implicit for sake of brevity.
\end{definition}
Model $\bar{\cP}$ includes mixture structures of different intrinsic dimensions with varying sequences of dimensions $(d_k)_{k \leq K}$ in a single unified model.
If present, a mixture component $P_k$ is assigned significant weight $\alpha_k \geq \alpha_{\min} > 0$.
A random variable $X \sim 
\sum_{k=1}^K \alpha_k P_k$ can be represented in the following way:
\begin{itemize}
\item
Draw layer index $L \in \{1,\ldots,K\}$ according to
$\bbP(L = k) = \alpha_k$;
\item
Conditionally on $L$, draw $X$ according to $P_L$, that is
$X \mid \{L = k\} \sim P_k$.
\end{itemize}

As any reasonable non-parametric model built upon a $\mathcal{C}^2$-type geometric model, the statistical model $\bar{\cP}$ can be shown to be stable (up to constants) under ambient $\mathcal{C}^2$-diffeomorphisms (see~\cite[Proposition~1]{Aamari18}).
Additionally, model $\bar{\cP}$ is structurally stable (up to constants) under mixture.
Indeed, the union of stratified manifolds in the model still belongs to the model.
At the level of data points, this means that merging several datasets generated from different distributions yields a dataset that can itself be seen as being generated from the model.

\begin{rem}[Reinterpreting clutter noise]
\label{rem:clutter-as-ambient-manifold}
A random variable $X$ drawn from model $\bar{\cP}$ originates from a mixture distribution $\sum_{k=1}^K \alpha_k P_{k}$ over manifolds of different dimensions.
It generalizes the \textit{clutter noise model} studied in \cite{Aamari18,Genovese12b}, where data could be drawn either from (signal) a smooth $d$-dimensional submanifold $M \subset \bbR^D$ with positive reach, or (clutter) from a uniform distribution over an ambient $D$-dimensional ball $\ball(0,R)$ such that $\ball(0,R/2) \supset M$.
To see this, one simply has to consider $\ball(0,R)$ as a $D$-dimensional (flat) submanifold of $\bbR^D$ and take $\{d_1,d_2\} = \{d,D\}$.
Note the slight abuse here: because we only deal with manifolds \emph{without} boundary, our model only formally allows for $d < D$.
In the context of this remark, this is of no importance, since we are not interested in the estimation of the $D$-dimensional stratum.
\end{rem}

\subsection{Statistical setting}
\label{sec:stat-setting}
Let $(X_1,L_1),\ldots,(X_n,L_n)$ be an i.i.d. sequence with $X_i \sim \sum_{k=1}^K \alpha_k P_k \in \bar{\cP}$ and $X_i \mid L_i \sim P_{L_i}$. 
We observe the \emph{unlabeled} i.i.d. sample $X_1,\ldots,X_n$ distributed according to the mixture $P = \sum_{k=1}^K \alpha_k P_k$, and we write
\begin{align*}
N_{k}
:=
\sum_{i = 1}^n \ind_{L_i = k}
\end{align*}
for the (unobserved) number of data points generated by the $d_k$-dimensional measure $P_k$.
In what follows, we denote by
\begin{align*}
P_{k,n}
:=
\frac{1}{{\color{black}1\vee N_{k}}} \sum_{i = 1}^n \ind_{L_i = k} \delta_{X_i}
\end{align*}
for the (unobserved) empirical distribution of the population distribution $P_k = P_{k,\infty}$ obtained with the true layer labels $L_i$'s. {\color{black}This sets $P_{k,n}=0$ when $N_k=0$.}

\paragraph*{What can and cannot be estimated on immersed manifolds}
\label{sec:estimation-impossibilities}

As mentioned above, reach positiveness \cite{Federer59} is violated when immersions are not embeddings.
As it has long been seen as a key regularity parameter and even considered to yield a \textit{minimal} smoothness notion (see references in ~\cite{Aamari19b,aamari2023optimal}), one may wonder whether various geometric quantities of a support $M = \bigcup_{k=1}^K M_k$ can be estimated from sample without reach condition, at least in the large sample limit.
We will say that a given distribution and data dependent quantity of interest $\theta_{\cX_n}(P)$ \textit{cannot} be estimated uniformly over the model $\bar{\cP}$ for the loss $\ell(\cdot,\cdot)$ if
\begin{align*}
\liminf_{n \to \infty}
\inf_{\hat{\theta}_n}
\sup_{P \in \bar{\cP}}
\mathbb{E}_{P^n}
\left[ \ell\bigl( \theta_{\cX_n}(P), \hat{\theta}(\cX_n) \bigr) \right]
>
0
,
\end{align*}
where the infimum ranges among all the estimators $\hat{\theta}(\cX_n)$ based on $n$ sample $\cX_n := \{X_1,\ldots,X_n\}$.
Let us list briefly some important geometric quantities that become impossible to estimate on mixtures of immersions.

\begin{itemize}
\item
The dimension from which a given point (say $X_1$) cannot be estimated. For instance, a small-dimensional component $M_1$ may be completely included in another higher-dimensional one $M_2$ (see orange points in Figure~\ref{fig:slabs-in-action}).
In this case, points from $M_2$ drawn very close to $M_1$ are statistically indistinguishable from that drawn exactly on $M_1$, and the dimension at $X_1$ becomes impossible to estimate.
Yet, we shall show that the dimensions involved in the mixture can actually be estimated (Theorem~\ref{thm:dimension-and-component-estimation}), as well as the individual dimensions of the $X_i$'s far enough from intersections (Theorem~\ref{thm:main-dimension-labeling}), yielding an overall precise dimension-wise clustering of the points (Theorem~\ref{thm:main-clustering-dimensions}).

\item Tangent spaces \textit{at} $X_i$'s cannot be estimated consistently, because of points near transversal self-intersections, see \cite[Theorem~1 and Figure~2]{Aamari19b}.
This fact is also clear given the fact that estimating a tangent space requires to identify the manifold dimension at a given point.
Yet, as we shall show that accurate tangent directions of the manifold mixture can be estimated in small neighborhoods (Theorem~\ref{thm:main-tangent-estimation}).

\item
For the exact same reasons, the locations of self-intersections cannot be estimated, or even detected. 
Therefore, any method that would rely on a preliminary preprocessing of possible singular neighborhoods (see \cite{lim2023hades}) is bound to fail.

\item
Homology also becomes inaccessible even asymptotically without reach condition, because topology is unstable without global constraints \cite[Theorem~1]{balakrishnan2012minimax}.
While the homology of general immersed manifolds can be estimated when self-intersections are transversal~\cite{bendich2007inferring,tinarrage2023recovering}, this feature becomes not identifiable for general immersions (see Figure~\ref{fig:topology-not-identifiable}).

\item
Geodesic distances are inestimable for similar reasons as homology. 
Since connectedness (i.e. zeroth order Betti number) cannot be detected without reach condition, one cannot even test for $\mathrm{d}_M(X_1,X_2) < \infty$ consistently.
Under an additional connectedness assumption, one may also construct branching structures with very different behaviors at arbitrarily small Euclidean scales with bottleneck structures~\cite{Aamari19b}.

\item
The volume behaves the same, as the volume of a single sphere can hardly be distinguished from the union of two spheres with same center, with one of them having an arbitrarily close (but different) radius. 

\item
If non-empty, the boundary $\partial M$ cannot be estimated without reach condition~\cite{Aamari21}, even if its own curvature is forced to remain bounded.
To see this, consider a fixed circle $M$ and remove an infinitesimal segment to it.
Changing the location of this deletion yields significant displacement of $\partial M$, while its smallness makes it undetectable.
Cartesian products with spheres yields curvature-bounded examples in higher dimensions.
\end{itemize}
We insist on the fact that all the above (counter-)examples yielding inestimability of some functionals can also be built under $\mathcal{C}^\infty$ constraints instead of $\mathcal{C}^2$ curvature bounds.
\section{The \codetection algorithm}

\subsection{Subspace detection with slabs}
\label{sec:heuristic-slab-infinite-sample}

The overall strategy of \codetection is to identify intrinsic neighborhoods of data points, that is to say aligned subsets of data likely to emanate from the same low-dimensional parametrization.
Let us develop on how to make this idea precise.

Let $P_k \in \cP^{(d_k)}_{{\kappa_{\max},\nu_{\max}}}(a_{\min},a_{\max})$ be the $d_k$-dimensional mixture component of $P = \sum_{k=1}^K \alpha_k P_k \in \bar{\cP}$.
In order to learn the stratification of $P$ from sample $\{x_1,\ldots,x_n\} = \cX \subset M = \bigcup_{k=1}^K M_k$, we will be interested in the empirical mass of specific neighborhoods $S \in \cS$ of points in space.
To simplify the exposition, a three-step reduction is in order.
\begin{itemize}
\item
\textit{(From sample to distributions)}
Given a class $\cS$ of measurable subsets of $\bbR^D$, a Vapnik–Chervonenkis argument (see Section~\ref{sec:appendix-concentration}) yields that with high probability, for all $S \in \cS$ and all $k\in \{1,\ldots,K\}$ simultaneously,
$$ -\sqrt{P_k(S) \VC(\cS) \log n/n} \lesssim P_{k,n}(S)-P_k(S) \lesssim 
\sqrt{P_{k,n}(S) \VC(\cS) \log n/n}
.
$$
Hence, if relying only on empirical counts over elements of a class $\cS$ with bounded Vapnik-Chernovenkis dimension $\VC(\cS)$, one can reason at the population level $P_k = P_{k,\infty}$ up to fluctuations.
\item
\textit{(From distributions to manifolds)}
By definition of the model, distributions $P_k(\diff x)$ behave roughly like the volume measure $\vol_{M_k}(\diff x) = \ind_{M_k} \cH^{d_k}(\diff x)$, up to constants depending on $a_{\min}, a_{\max}$ and $\kappa_{\max}$ (Definition~\ref{def:model-mixture}).
Distributional considerations hence boil down to volumic ones.
\item
\textit{(From manifolds to linear spaces)}
{\color{black}As shown in Section~\ref{sec:sheaf-number},} because each $M_k$ is an immersed manifold with curvature bounded by $\kappa_{\max}$, the neighborhood $\ball(x,h_{\para})\cap M_k$ of each point $x \in M_k$ looks like a finite union 
$\bigcup_{\ell \leq N_0} \ball_{T_\ell}(y_\ell,h_\ell)$ of $d_k$-dimensional tangent balls {\color{black}(Definition~\ref{dfn:slab}).} 
This approximation holds for $h_{\para} \lesssim 1/\kappa_{\max}$ up to an error of order $\kappa_{\max} h_{\para}^2$ in the direction normal to the tangent spaces nearby $x$ (Section~\ref{sec:geodesic-estimates}).
\end{itemize}
To best identify the local structure $M_k \cap \ball(x,h_{\para}) \simeq \bigcup_{\ell \leq N_0} \ball_{T_\ell}(y_\ell,h_\ell)$ up to quadratic terms, we make use of thickened low-dimensional balls, called \textit{slabs}.

\begin{figure}[!ht]
\centering
\includegraphics[width=1\linewidth,page=5]{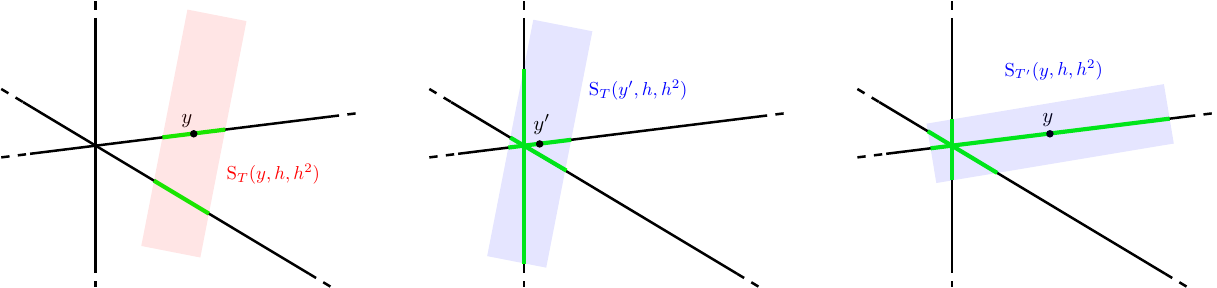}
\caption{Stratification detection for a union of lines ($d=1$) in the plane ($D=2$) via slabs.
The green lines correspond to the intersection with the hidden structure. 
To yield large enough intersected volume, a slab needs to be both well-aligned \emph{and} nearly centered at the structure (blue). Otherwise, it yields small intersected volume (red).
}
\label{fig:linear-intersections-explained}
\end{figure}

\begin{figure}
\centering
\includegraphics[width=0.9\linewidth]{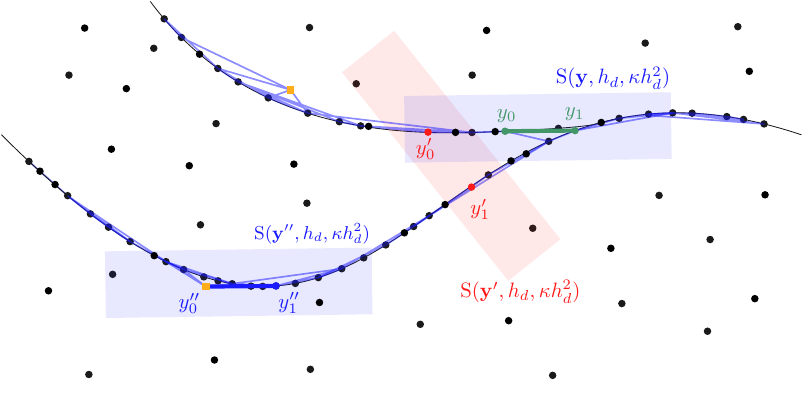}
\caption{
A typical configuration described by Lemma~\ref{lem:recalgo}.
At step $d=1$, all the close-enough pairs of sample points are jointly considered and added to $\cU^{(1)}$ if the slab they define contains enough empirical mass. 
In the displayed layout, pairs $\mathbf{y}$ and $\mathbf{y}''$ yield co-detection (blue), while $\mathbf{y}'$ does not (red).
When co-detected, the convex hull of points $\mathbf{y} = \{y_0,\ldots,y_d\}$ contribute to $\hat{M}_d$~\eqref{eq:def-clusters-and-manifold-estimator}.
Points from higher-dimensional layers (orange) or pairs from different immersions (green) can actually be co-detected. If so, they provably yield convex hulls that are $O(h_d^2)$-close to $M_1$ (Lemma~\ref{lem:recalgo}), so that they do not impair the estimation rates.
}
\label{fig:slabs-in-action}
\end{figure}

{\color{black}
\begin{definition}[Tangent ball, Slab] \label{dfn:slab}
Given $x \in \bbR^D$ and a linear subspace $T \subset \bbR^D$ and $h > 0$, we write
$$
\ball_T(x,h) := \{y \in (x+T)~|~\|y-x\| \leq h\} = \ball(x,h) \cap (x+T).
$$
For two positive radii $h_{\para},h_\perp$, we define the slab directed by $T$ as 
$$
\slab_T(x,h_{\para},h_\perp)
:=
x + \ball_T(0,h_{\para}) + \ball_{T^\perp}(0,h_\perp)
.
$$
In the rest of the article, we might refer to $T$ as a linear subspace of as a \emph{flat} indistinctively.  
\end{definition}}
For inference, the normal parameter $h_{\perp}$ is to be chosen of order $\kappa_{\max} h_{\para}^2$.
Indeed if $h_{\perp} := a \kappa_{\max}h_{\para}^2$, we show in Section~\ref{sec:volume-estimates} that for all $d$-dimensional linear subspace $T$ (to be detected), and all (candidate) $d'$-dimensional linear subspace $T'$ with $d' \leq d$, then as $h_{\para} \to 0$,
\begin{align*}
\cH^d\bigl( \ball_{T}(y,h_{\para}) \cap \slab_{T'}(x,h_{\para},h_{\perp})\bigr)
\gtrsim h_{\para}^{d'}
\Leftrightarrow
\begin{cases}
d' = d, \\
\angle(T,T') \lesssim h_{\para},
\\
\|y - x\| \lesssim h_{\para}^2,
\end{cases}
\end{align*}
where $\angle(T,T') := \|\pi_T - \pi_{T'}\|_{\mathrm{op}}$ stands for the largest principal angle between $T$ and $T'$.
This volume characterization naturally extends to finite unions of subspaces (see Figure~\ref{fig:linear-intersections-explained}), in the sense that if $d' \leq \min_{\ell \leq N_0} \dim(T_\ell)$, then
\begin{align}
\label{eq:volume-heuristic}
\cH^d\bigl( \bigcup_{\ell \leq N_0} \ball_{T_\ell}(y_\ell,h_\ell) \cap \slab_{T'}(x,h_{\para},h_{\perp})\bigr)
\gtrsim h_{\para}^{d'}
\Leftrightarrow
\exists \ell \leq N_0, 
\begin{cases}
d' = \dim(T_\ell), \\
\angle(T_{\ell},T') \lesssim h_{\para},
\\
\|y_\ell - x\| \lesssim h_{\para}^2.
\end{cases}
\end{align}
To detect the location and orientation of all the structures of dimensions at most $d_{\max}$, this property suggests to evaluate the intersected volumes of slabs of dimension $d' \in \{1,\ldots,d_{\max}\}$ in increasing order, which leads to the \codetection algorithm.

\subsection{Nonlinear subspace detection from sample}

At the sample level, the heuristic on slabs from the previous section naturally casts to sample-aligned slabs, which are defined as follows.

\begin{definition}[Slab of a $(d+1)$-tuple]
\label{def:slab}
Let $d \in \{1,\ldots,D\}$.  
Given $\mathbf{y} := \{y_{0},\ldots,y_{d}\} \subset \bbR^D$, the \emph{slab} defined by $\mathbf{y}$, tangent width $h_{\para}>0$ and normal width $h_\perp >0$ is
\begin{align*}
\slab(\by , h_{\para}, h_{\perp} )
&:=
\slab_{T(\by)}(\bar{\by},h_{\para},h_{\perp}),
\end{align*}
where $\bar{\mathbf{y}} = \frac{1}{d+1} \sum_{j=0}^d y_j$ is the barycenter of $\mathbf{y}$ and $T(\by) := \Span\left(y_1-y_0,\ldots,y_d-y_0 \right)$ its linear span.
\end{definition}

\codetection enumerates dimensions $d$ from $1$ to a predefined one $d_{\max}$ in increasing order. Starting from the full dataset $\widetilde{\cX}^{(1)} := \cX$, it constructs decreasing subsets $\widetilde{\cX}^{(1)} \supset
\ldots
\supset \widetilde{\cX}^{(d_{\max}+1)}$.
When at current step $d$, points from $\widetilde{\cX}^{(d)}$ are added to the output $\widehat{\cX}^{(d)}$ by $(d+1)$-tuples $\by := (y_0,\ldots,y_d)$. 
Tuple $\by$ is added to $\widehat{\cX}^{(d)}$ when
\begin{itemize}
\item
it has radius $\rad(\by) \leq r_d$ (where $\rad(\by)$ is the radius of the smallest enclosing ball of $\by$), and 
\item
its associated slab $\slab(\by,h_d, \kappa_d h_d^2 )$ contains at least $n_d$ points from $\widetilde{\cX}^{(d)}$.
\end{itemize}
In practice, $h_d \asymp (\log n / n)^{1/d}$ corresponds to the typical sampling scale of a $n$-sample drawn uniformly over a $d$-dimensional ball.
Following the intuition from the population level (Section~\ref{sec:heuristic-slab-infinite-sample}), this criterion indicates that $\{y_0,\ldots,y_d\}$ are likely to jointly lie on a single $d$-dimensional stratum of $M = \bigcup_{k=1}^K M_{k}$: we say that they have been \emph{co-detected}.
\codetection also stores this alignment information by adding $\by$ to a list $\cU^{(d)}$ of $(d+1)$-tuples.

At the end of step $d$, \codetection also adds to $\wh\cX^{(d)}$ the set of points $\cR_d \subset \wt\cX_d$ that are at distance at most $\delta_d$ (of order $\kappa_d h_d^2$ in practice) from the co-detected convex hulls of $\cU^{(d)}$.
This extra step is made to ensure that after dimension $d = d_k$ has been parsed, all the points of $M_{k}$ have been included $\wt\cX^{(d)}$.
This avoids any contamination in $\wt\cX^{(d_k+1)} = \wt\cX^{(d_k)} \setminus \wh\cX^{(d_k)}$ from such possibly forgotten points when then parsing higher-dimensional layers.

Some of the lists $(\cU^{(d)})_{1 \leq d \leq d_{\max}}$ may be empty. 
The main output of the \codetection algorithm is the collection of non-empty tuples $(\wh \cU_k)_{1 \leq k \leq \wh K}$. 
From it, the approximate stratified structure of $M$ can be constructed. Namely, the estimated number of layers and intrinsic dimensions are taken to be
\begin{align}
\label{eq:estimators-discrete}
\wh K := \Card\{d \mid \cU^{(d)} \neq \emptyset\}
\text{\quad and \quad}
\{\hat d_1,\dots,\hat d_{\wh K}\} := \{d \mid \cU_d \neq \emptyset\}
,
\end{align}
and the dimension-wise data points clusters and support estimators are, for all $k \in \{1,\ldots,\wh K\}$
\begin{align}
\label{eq:estimators-sets}
\wh{\cX}_k := 
\wh{\cX}^{(\hat{d}_k)}
\text{\quad and \quad}
\wh{M}_k := \bigcup_{\by \in \wh\cU_k} \conv(\by).
\end{align}

\begin{algorithm}
\caption{
\codetection  
}
\label{algo:strat}

\KwData{Point cloud $\cX \subset \bbR^D$
\newline
Maximal dimension $d_{\max}$;
\newline
Slab radii $(h_d)_{1 \leq d \leq d_{\max}}$;
\newline
Slab aspect ratios $(\kappa_d)_{1 \leq d \leq d_{\max}}$;
\newline
Localization radii $(r_d)_{1 \leq d \leq d_{\max}}$;
\newline
Pruning radii $(\delta_d)_{1 \leq d \leq d_{\max}}$;
\newline
Size cutoffs $(n_d)_{1 \leq d \leq d_{\max}}$;
}
\KwResult{Number of components $\hat{K}$;
\newline
Set of dimensions $(\hat d_1,\dots, \hat d_{\hat K})$;
\newline
Labeled datapoints $(\wh\cX_1,\dots,\wh\cX_{\hat K})$;
\newline
Collection of tuples $(\wh \cU_1,\dots,\wh \cU_{\hat K})$;
}
$d \gets 0$\;
$\textrm{dim list} \gets \{\}$\;
$\textrm{points list} \gets \{\}$\;
$\textrm{tuples list} \gets \{\}$\;
$\hat K \gets 0$\;
$\wt\cX \gets \cX$\;

\While{$\wt\cX \neq \emptyset$ \textbf{\emph{and}} $d < d_{\max}$}{
$d \gets d + 1$\;
$\wh \cX \gets\emptyset$\;
$\cU \gets \emptyset$\;
\For{all $(d+1)$-tuple $\by = (y_0,\dots,y_{d})$ of distinct points of $\wt\cX$}
	{
	\If{$\rad(\by) \leq r_{d}$ \textbf{\emph{and}} $\Card\bigl(\wt \cX \cap \slab(\by,h_{d},\kappa_d h_{d}^2)\bigr) \geq n_d$}
		{
		$\cU \gets \cU \bigcup \{\by\}$\;
		$\wh \cX \gets \wh\cX \bigcup \{y_0,\dots,y_d\}$\;
		}
	}
	$\cR \gets \emptyset$\;
	\For{all $z$ in $\wt\cX\setminus\wh\cX$}
		{
		\If{$\dist(z,\bigcup_{\sigma \in \cU} \conv \sigma) \leq \delta_d$}
			{
			$\cR \gets \cR \bigcup \{z\}$\;
			}
		}
	$\wh\cX \gets \wh\cX \bigcup \cR$\;
	\If{$\wh\cX \neq \emptyset$}
		{
		Append $d$ to \textrm{dim list}\;   
		Append $\wh\cX$ to \textrm{points list}\;
		Append $\cU$ to \textrm{tuples list}\;
		$\wt \cX \gets \wt \cX \setminus \wh \cX$\;
		$\hat K \gets \hat K + 1$\;
		}
}
\end{algorithm}

\begin{rem}[On leftover data points] 
The points remaining in $\wh \cX_{d_{\max}+1} = \cX \setminus \bigcup_{1 \leq k \leq \wh K} \wh \cX_{k}$ after sieving all the dimensions $1$ through $d_{\max}$ are likely to have been generated from a distribution that has intrinsic dimension at least $d_{\max}+1$.
As such, they can be considered as originating from too many degrees of freedom, making them hardly analyzable statistically speaking if one has fixed a dimension threshold over which the curse of dimensionality is considered too strong.
Hence, they can naturally be labeled as clutter noise, in the same spirit as Remark~\ref{rem:clutter-as-ambient-manifold}.
\end{rem}
\section{Optimal stratification learning}

In what follows, we will denote by $C,C',C'',c,c',c''$ generic positive numerical constants.
Their value can differ from one line of an equation to another.
Likewise, we will use $C_d,C'_d,C''_d,c_d,c'_d,c''_d$ for generic constants depending on dimension $d \in \{1,\ldots,D\}$ only. 

\subsection{Geometric loss functions}
In the results to come, we are interested in the exact identification of the integers given by the number of clusters $K$, as well as dimensions $\{d_1,\ldots,d_K\}$.
To measure precision on point location and set estimation, the \emph{distance function} $\dist(\cdot \mid A)$ to a subset $A \subset \bbR^D$ is defined for all $x \in \bbR^D$ as
$$
\dist(x \mid A)
:=
\inf_{a \in A} \| x - a \|
.
$$
The \emph{one-sided Hausdorff distance} from $A \subset \bbR^D$ to $B \subset \bbR^D$ is 
$$
\dH(A\mid B)
:=
\sup_{a \in A} \dist(a \mid B)
,
$$
and the \emph{Hausdorff distance} between $A$ and $B$ is the symmetrized quantity
$$
\dH(A,B)
:=
\max\{\dH(A\mid B), \dH(B\mid A)\}
.
$$
When dealing with linear subspaces $T,T' \subset \bbR^D$, we will be interested in the \emph{principal angle}
$$
\angle(T,T') := \| \pi_T - \pi_{T'} \|_{\mathrm{op}},
$$
where $\pi_T,\pi_{T'} : \bbR^D \to \bbR^D$ stand for the orthogonal projections onto $T$ and $T'$ respectively, and $\| \cdot \|_{\mathrm{op}}$ is the (Euclidean) operator norm on endomorphisms of $\bbR^D$. When comparing two union of flats indexed by $\bbT$ and $\bbT'$, we let
$$
\angle(\bigcup_{T \in \bbT} T \mid \bigcup_{T'\in\bbT'} T')  := \sup_{T \in \bbT} \inf_{T' \in \bbT'}\angle(T,T'),
$$
which coincides with the previous definition when $\bbT$ and $\bbT'$ are restricted to single flats.

\subsection{Parameter settings and algorithmic guarantees}
We are now in position to state the main properties of the \codetection algorithm. In all the results below, we shall pick the following collection of parameters to ensure good properties:
\begin{align}
\label{eq:parameter-choice}
\begin{cases*}
d_{\max} = D-1
\\
h_d = \dfrac{48 d}{\kappa_{\max}}\(1+\dfrac1d\) \left(\dfrac{\Upsilon\gamma \log n}{a_{\min} n}\right)^{1/d}
\\
\kappa_d = \kappa_{\max}
\\
r_d = h_d
\\
\delta_d = \zeta_d \kappa_{\max} h_d^2
\\
n_d = \sigma \gamma \log n
\end{cases*}
\end{align}
where
\begin{itemize}
\item
$\Upsilon$ {\color{black}is to satisfy $\Upsilon \geq C_D\nu_{\max}/\alpha_{\min}$ for a sufficiently large constant $C_D$ depending on $D$}, 
\item
$\sigma$ is to be taken of order $D$ up to multiplicative constants,
\item
$\zeta_d$ {\color{black}shall be taken equal to $1/4$},
\item
$\gamma$ is to be taken at least of order $(q \vee D^2 \log D)/\alpha_{\min}$, where

\item
$q > 0$ is a parameter that will tune the decay of probability of success of the algorithm.
\end{itemize}
We refer to \lemref{recalgo} for a more precise statement regarding these constants.
This lemma gathers the key geometric properties kept invariant throughout the execution of \codetection on input $\cX_n = \{X_1,\ldots,X_n\}$ generated independently and identically distributed from some unknown $P \in \bar{\cP}$ (Definition~\ref{def:model-mixture}).
With probability at least $1 - O(n^{-q})$,
at all steps $d \in \{1,\ldots,D-1\}$, the sample counterpart of \eqref{eq:volume-heuristic} holds. That is,
\begin{itemize}
\item 
\textit{(Discrete stratification learning)}
$\cU^{(d)} \neq \emptyset$ if and only if $d \in \{d_1,\ldots,d_K\}$.

\item
\textit{(Approximation)}
If $d = d_k$ for some $k \in \{1,\ldots,K\}$, then
\begin{itemize}
\item
For all $\by = (y_0,\ldots,y_{d_k}) \in \cU_{d_k}$ co-detected by \codetection
and all $y \in \Conv(\by)$, there exists $x \in \ball\bigl(y,O(h_{d_k}^2)\bigr) \cap M_k$ such that
$$
\|y-x\|
\lesssim h_{d_k}^2
\text{\quad and \quad}
\angle\bigl(\Span(y_1-y_0,\dots,y_{d_k}-y_0), T_x M_k\bigr) \lesssim  h_{d_k}
.
$$

\item
For all $x \in M_k$ and all $d_k$-dimensional subspace $T \subset T_x M_k$, there exists $\by = (y_0,\ldots,y_{d_k})\in \cU_{d_k}$ co-detected by \codetection, and $y \in \ball\bigl(x,O(h_{d_k}^2)\bigr) \cap \Conv(\by)$ such that
$$
\|x-y\|
 \lesssim h_{d_k}^2
\text{\quad and \quad}
\angle\bigl(T, \Span(y_1-y_0,\dots,y_{d_k}-y_0)\bigr) \lesssim  h_{d_k}
.
$$
\end{itemize}
\end{itemize}
These properties are sufficient to establish the main results of this article.
\begin{rem}[On average time complexity]
The above choice of parameters happens to turn \codetection into a quasi-linear time algorithm on average, provided that one already has computed local neighborhood information on $\cX_n$~\cite{har2011geometric}. 
Indeed, \lemref{recalgo} asserts with high probability, that no point from any $d'$-dimensional layer is ever forgotten along the way for all $d'<d$.
Therefore, at step $d$, each remaining point  $X_i \in \cX_n$ has at most $O(n r_d^d) = O(\log n)$ candidate neighbors for its co-detection. 
As a result, \codetection enumerates at most $O((\log n)^d)$ slabs containing $X_i$. 
This results in an overall time complexity $O(n (\log n)^D)$ on average. 
{\color{black}
With our choice of parameters and on the high-probability event considered in Lemma~\ref{lem:recalgo}, the local neighborhood can be computed at each step $d$ using standard methods in $O(n\log n)$ time and space, with an implicit constant depending exponentially on the ambient dimension $D$~\cite{bentley1977complexity}.
}

\end{rem}

\subsection{Main results}
Using notation from Section~\ref{sec:stat-setting}, we let $\cX_n := \{X_1,\ldots,X_n\}$ be drawn iid from some unknown $P = \sum_{k=1}^K \alpha_k P_k \in \bar{\cP}$ and associated labels $(L_1,\ldots,L_n)$, and $\cX_{k,n} := \{ X_i \mid L_i =  k\}$.
The proof of the results below can be found in Section~\ref{sec:proof-main-results}.
\paragraph*{Stratification identification}
Under this sampling scheme and parameter choice as above, the estimated discrete structure
\begin{align}
\{\hat d_k\} := \{d \mid \cU^{(d)} \neq \emptyset\}
\text{\quad and \quad}
\wh K := \Card\bigl(\{\hat d_k\}\bigr)
,
\tag{\ref{eq:estimators-discrete}}
\end{align}
identifies the underlying one precisely with high probability. 
In all the subsequent results, $n$ is taken large enough compared to a constant depending on the parameters of $\bar \cP$ {\color{black}and the fixed, $n$-independent choices of $\Upsilon,\gamma,\sigma,\zeta_d,q$}.
\begin{theorem}[Mixture structure identification]
\label{thm:dimension-and-component-estimation}
With parameters chosen as in \eqref{eq:parameter-choice}, for $n$ large enough, with probability at least $1-6Kn^{-q}$,
\begin{center}
$\wh K = K$ \quad and \quad $\hat d_k = d_k$ for all $k \in \{1,\ldots,K\}$.
\end{center}
\end{theorem}
The recovery of the number of clusters $K$ is, to our knowledge, the first mathematically grounded result of this type (see~\cite{Wu_2022_WACV} for empirical evidences on various methods).  
The dimension(s) estimation part extends the results of~\cite{Kim15}, which hold only under reach condition and with a single mixture component ($K=1$).

Beyond dimension(s) identification, \codetection also assigns a dimension to each data point.
Recall that for all $k \in \{1,\ldots,\wh K\}$,
\begin{align}
\label{eq:def-clusters-and-manifold-estimator}
\tag{\ref{eq:estimators-sets}}
\wh{\cX}_{k,n} := 
\wh{\cX}^{(\hat{d}_k)}
\text{\quad and \quad}
\wh{M}_k := \bigcup_{\by \in \cU_k} \conv(\by)
\end{align}
For all $X_i \in \cX_n$, we let $\wh d(X_i)$ denote \emph{the} dimension $d_k$ to which it belongs in the disjoint union $\cX_n = \bigsqcup_{k = 1}^{\wh K} \wh{\cX}_{k,n}$, i.e. $\wh d(X_i) = \wh d_k$ if and only if $X_i \in \wh{\cX}_{k,n}$.
As mentioned in Section~\ref{sec:estimation-impossibilities}, the exact individual dimension $d_{L_i} = \dim(M_k)$ from which $X_i$ has been drawn cannot be identified. 
It may not even be defined uniquely, for instance if $M_k \subset M_{k+1}$.
However, it can be underestimated locally in the following sense.
\begin{theorem}[Local dimension estimation]
\label{thm:main-dimension-labeling}
With parameters chosen as in \eqref{eq:parameter-choice}, for $n$ large enough, with probability at least $1-6Kn^{-q}$,
\begin{align*}
\text{for all~} i\in \{1,\ldots,n\}, \quad
\min\{d_k~|~ \dist(X_i \mid M_k) \leq \tau_k \(\frac{\log n}{n}\)^{2/d_k}\} 
\leq  
\wh d(X_i)
\leq  
d_{L_i},
\end{align*}
where 
$$
\tau_k := C_{d_k} \Upsilon \nu_{\max}\frac{a_{\max}}{a_{\min}\kappa_{\max}} 
\(
\frac{\Upsilon\gamma}{a_{\min}}
\)^{2/d_k}.
$$
\end{theorem}

In accordance with intuition, the above pointwise statement leaves estimation ambiguity only nearby intersections between layers of different dimensions.
In particular, it turns into \emph{exact} dimension labeling for {\color{black}a single embedded submanifold ($K=1$)} with reach bounded below by $1/\kappa_{\max}$.
Also note that the localization radius $O(h_{d_k}^2)$ is dimension-dependent: the lower the intrinsic dimension, the most sharply points can be identified as coming from it.

Coming back to global sample identification properties, the following statement asserts that most of the data points are assigned to the correct layer.
Recall that on average we have that $N_k := \# \cX_{k,n}$ satisfies 
$\bbE_{P}[N_k] = \alpha_k n \asymp n$, because $\alpha_{\min} > 0$ is considered fixed.
{\color{black}For the loss below, we set $\wh\cX_{k,n}:=\emptyset$ when $k>\wh K$.}
\begin{theorem}[Dimension-wise clustering] 
\label{thm:main-clustering-dimensions}
With parameters chosen as in \eqref{eq:parameter-choice}, for $n$ large enough,
for all $k \in \{1,\ldots,K\}$ simultaneously,
\begin{align*} 
&
\bbE_P\left[\frac{\#\(\wh\cX_{k,n} \triangle \cX_{k,n}\)}{N_k \vee1} \right]
\leq 
{\color{black}\eta_k\left\lbrace\left(\frac{\log n}{n}\right)^{2\frac{d_k - d_{k-1}}{d_{k-1}}} \vee \left(\frac{\log n}{n}\right)^{2\frac{d_{k+1} - d_k}{d_k}}\right\rbrace + 6Kn^{-q},}
\\
&\text{with}\quad 
{\color{black}\eta_k = C_D a_{\max} \alpha^{-1}_{\min} \nu_{\max}^2\left\lbrace1\vee(\kappa_{\max}\tau_{k-1})^{d_k-d_{k-1}}\vee(\kappa_{\max}\tau_k)^{d_{k+1}-d_k}\right\rbrace,}
\end{align*} 
where $\cX \triangle \cX' := (\cX \setminus \cX') \sqcup (\cX' \setminus \cX)$ stands for the symmetric difference between $\cX,\cX' \subset \bbR^D$.
Here, we used convention $d_0 = 0$ and $d_{K+1} = \infty$. {\color{black}In $\eta_k$, terms involving absent neighbors are omitted, and $\tau_k$ is as in Theorem~\ref{thm:main-dimension-labeling}}.
\end{theorem}

The above result provides a rare example of a \textit{blessing of dimensionality} phenomenon.
Indeed, varying $(d_{k-1},d_{k+1})$ while leaving $d_k$ fixed, exponent ${\color{black}\beta_k := 2\left\lbrace\frac{d_k-d_{k-1}}{d_{k-1}}\wedge\frac{d_{k+1}-d_k}{d_k}\right\rbrace}$ ranges 
\begin{itemize}
\item
from $\beta_k = 2/d_{k}$ in worst case when $d_{k+1} = d_k+1$, yielding the simpler clustering rates
\[
\bbE_P\left[\frac{\#\(\wh\cX_{k,n} \triangle \cX_{k,n}\)}{N_k \vee1} \right]
\leq 
\eta_k \(\frac{\log n}{n}\)^{2/d_{k}} {\color{black}+6Kn^{-q}}
;
\]

\item
to $\beta_k =  2(d_k-1)$ when $(d_{k-1},d_{k+1}) = (1,\infty)$ (lower layers best case) and $\beta_k =  2\frac{D-1-d_k}{d_k}$ when $(d_{k-1},d_{k+1}) = (0,D-1)$ (higher layers best case).
\end{itemize}
This accords with the intuition that layers of more heterogeneous dimensions shall be more easily separated, and interfere less with each other.

\paragraph*{Stratification estimation}

Now moving to approximation properties at order zero (Hausdorff distance), we show that the support is estimated at dimension-dependent rate by the local convex hulls $\wh M_k := \bigcup_{\by \in \wh\cU_k} \conv(\by)$ of co-detected points.

\begin{theorem}[Stratified manifold estimation]
\label{thm:main-manifold-estimation}
With parameters chosen as in \eqref{eq:parameter-choice}, for $n$ large enough, with probability at least $1-6Kn^{-q}$, for all $k \in \{1,\ldots,K\}$ simultaneously,
\begin{align*}
\dH(M_k,\wh{M}_k)
&\leq 
\rho_k \left(\frac{\log n}{n}\right)^{2/d_k},
\\
\text{where} \quad \rho_k &= C_{d_k} \Upsilon\nu_{\max}\frac{a_{\max}}{a_{\min} \kappa_{\max}}
\(\frac{\Upsilon\gamma}{a_{\min}}\)^{2/d_k}.
\end{align*}
\end{theorem}
This rate is minimax optimal~\cite{Kim15}, even over the smaller model composed of \emph{embedded} submanifolds with reach bounded below and no proper mixture ($K=1$).
Strikingly, the absence of reach condition and the presence of interacting layers of different dimensions do not impact the Hausdorff estimation rate of each layer.
As mentioned in Remark~\ref{rem:clutter-as-ambient-manifold}, this result generalizes~\cite[Theorem~5]{Genovese12b} and~\cite[Theorem~2.9]{Aamari18}, which considered the estimation of a single manifold $M_{d}$ given an i.i.d. sampling of it corrupted with outliers generated by Huber's contamination model.
In addition to having weaker assumptions, note that Theorem~\ref{thm:main-manifold-estimation} also benefits from an explicit algorithmic description (absent of~\cite{Genovese12b}).
It also processes the estimation of the $d$-dimensional layer with a single pass on data, as opposed to the intricate iterative algorithm of \cite{Aamari18}.

At order one, the piecewise linear structure of $\wh M_k = \bigcup_{\by \in \wh\cU_k} \conv(\by)$ also provides consistent estimates of local tangent spaces. 
For $x\in \wh M_k$, we let
$$
T_x \wh M_k := \bigcup_{\substack{\sigma \in \wh \cU_k \\ x \in \conv\sigma}} {\color{black}T(\sigma)}.
$$
{\color{black}This directional information is already computed by \codetection: it helps reject misaligned tuples and controls the error of any accepted cross-branch simplex.}
As noted in Section~\ref{sec:estimation-impossibilities}, the angle criterion for tangent space estimation needs to be weakened for similar reasons to dimension estimation.
We consider a natural localized version of it by introducing a spatial tolerance $\Delta > 0$. Then, we write
\begin{align*}
\angle_\Delta \bigl( T M_k \mid T \wh M_k \bigr)
&:=
\sup_{x \in M_k}
{\color{black}\sup_{\substack{T\subset T_xM_k\text{ linear}\\\dim T=d_k}}}
\inf_{y \in \wh M_k \cap \ball(x,\Delta)}
\angle ({\color{black}T} \mid T_y \wh M_k)
,
\end{align*}
and
\begin{align*}
\angle_\Delta \bigl( T \wh M_k \mid T M_k \bigr)
&:=
\sup_{y \in \wh M_k}
{\color{black}\sup_{\substack{\wh T\subset T_y\wh M_k\text{ linear}\\\dim \wh T=d_k}}}
\inf_{x \in M_k \cap \ball(y,\Delta)}
\angle ({\color{black}\wh T} \mid T_x M_k)
.
\end{align*}
{\color{black}Thus each tangent branch composing the tangent cone is matched to a branch at its own nearby point. A common matching point for all branches is not required.}
We can then define the symmetric loss 
$$
\angle_\Delta \bigl( T \wh M_k, T M_k \bigr) := \angle_\Delta \bigl( T \wh M_k \mid T M_k \bigr)  \vee \angle_\Delta \bigl( T  M_k \mid T \wh M_k \bigr)
,
$$
which we bound in the following result.
\begin{theorem}[Tangent space estimation]
\label{thm:main-tangent-estimation}
With parameters chosen as in \eqref{eq:parameter-choice}, for $n$ large enough and with probability at least $1-6Kn^{-q}$, for all $k \in \{1,\ldots,K\}$ simultaneously,
\begin{align*}
\angle_{\Delta_k} \bigl( T \wh M_k, T M_k \bigr)
&\leq
\iota_k \left(\frac{\log n}{n}\right)^{1/d_k}, \\
\text{where}\quad
\iota_k &:= C_{d_k}
\Upsilon \nu_{\max}\frac{a_{\max}}{a_{\min}}
\(\frac{\Upsilon\gamma}{a_{\min}}\)^{1/d_k}, \\
\text{and} \quad
\Delta_k
&:={\color{black}C_{d_k}}
\Upsilon \nu_{\max}\frac{a_{\max}}{a_{\min}\kappa_{\max}}
\(\frac{\Upsilon\gamma}{a_{\min}}\)^{2/d_k}
 \(\frac{\log n}{n}\)^{2/d_k}.
\end{align*}
\end{theorem}
In words, this result asserts that the $d_k$-dimensional tangent structure of $M = \cup_{k=1}^K M_k$ can be estimated up to angle $O(h_{d_k})$ in neighborhoods of size $O(h_{d_k}^2)$.
Again, this rate is minimax optimal~\cite{Aamari19} for \emph{embedded} submanifolds with reach bounded below and no proper mixture ($K=1$).
{\color{black}For an embedded submanifold with reach bounded below by $1/\kappa_{\max}$, matching estimated points to their nearest-point projections yields the usual tangent-space accuracy, since tangent variation over a distance $O(h_{d_k}^2)$ is $O(h_{d_k}^2)$.}
\section{Discussion and open questions}

\subsection{Some possible practical guidance for parameter selection}
\label{sec:parameter-selection}

\begin{figure}[h!]
  \centering
  \begin{subfigure}[t]{0.48\textwidth}
    \centering
    \includegraphics[height=5cm]{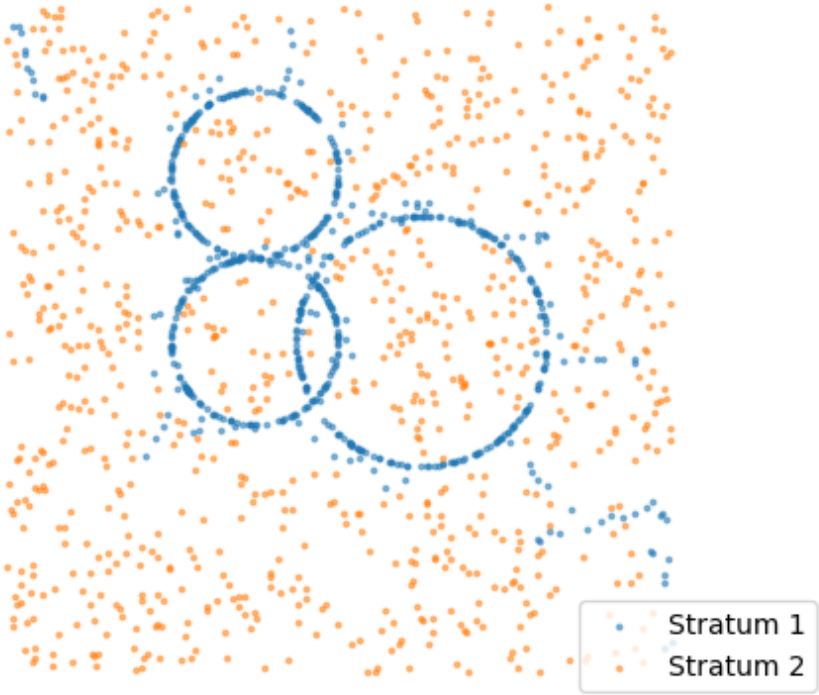}
    \includegraphics[height=5cm]{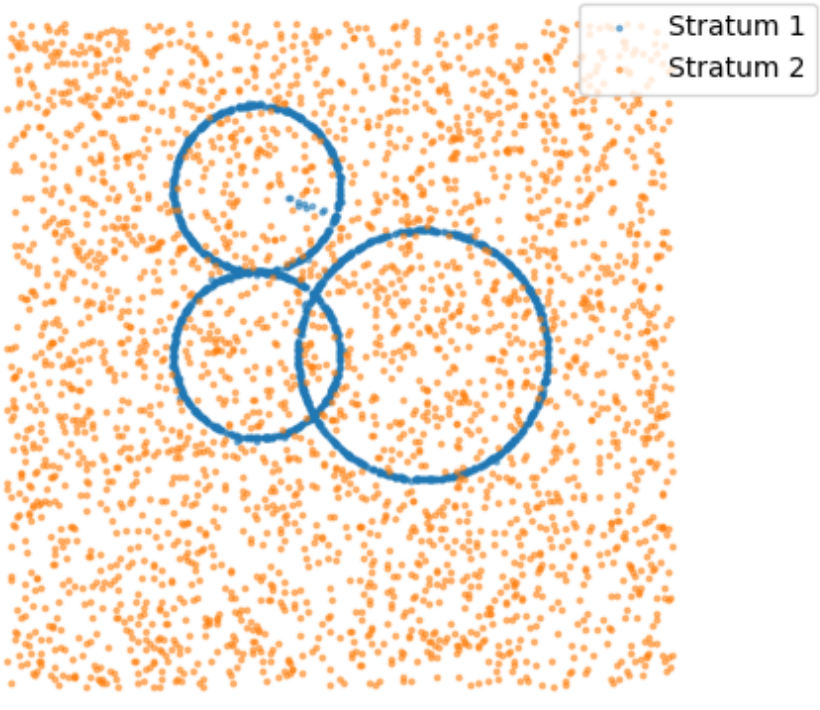}
    \caption{Output clustering.}
    \label{fig:subfig1}
  \end{subfigure}
  \hfill
  \begin{subfigure}[t]{0.48\textwidth}
    \centering
    \includegraphics[height=5cm]{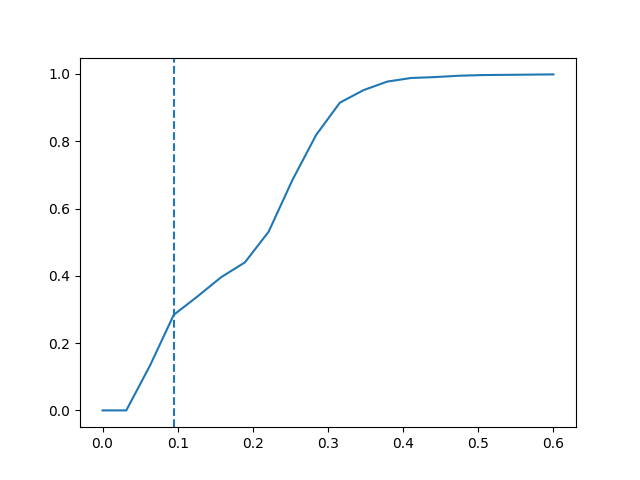}
    \includegraphics[height=5cm]{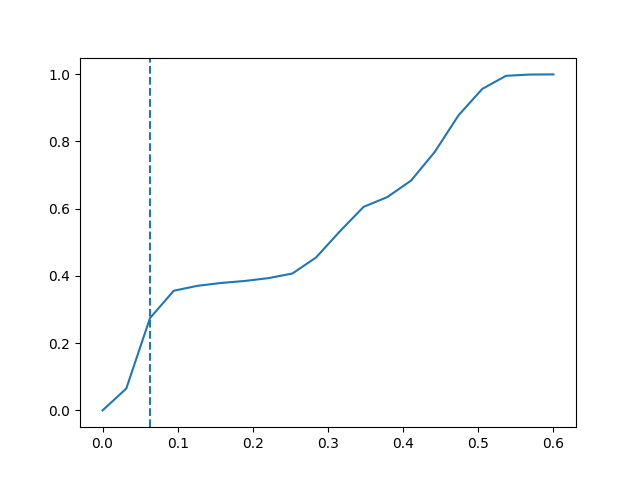}
    \caption{Elbow heuristic for final function $\mathcal{G}_{t}(\mu)$.}
    \label{fig:subfig2}
  \end{subfigure}
  \caption{
  Output clustering $(\widehat{\cX}^{(1)},\widehat{\cX}^{(2)})$ of \codetection, when computed with parameters selected according to the heuristic of Section~\ref{sec:parameter-selection}. (Top) $N_1 = 500$ points lie on the circles, and $N_2 = 1000$ in the ambient square and (Bottom) $N_1 = 1100$ and $N_2 = 3000$.
}
  \label{fig:simu}
\end{figure}

Despite \codetection theoretically requiring the choice of six (dimension-dependent) hyperparameters --- $d_{\max}$, $h_d$, $\kappa_d$, $r_d$, $\delta_d$, and $n_d$ --- let us discuss some simple heuristics that may facilitate its implementation in practice.
These heuristics rely on only two sequences of dimension-dependent hyperparameters, denoted by $\varepsilon_d, \mu_d \in [0,1]$, whose selection will be discussed later.
Namely, a possible parameter selection proceeds as follows.
\begin{enumerate}[leftmargin=*]
\item
First, if no computational constraint has to be accounted for, one may always choose the maximal dimension to be $d_{\max} = D-1$ in order to search for all the possible low-dimensional structures of a dataset.
\item
In practice, to adapt to the local geometry of data, we suggest to use the distance to a fixed number of nearest neighbors $\ell= 5 \lceil \log n \rceil$ to take the place of the localization radius $r = r(X_i)$ and the slab radii $h = h(X_i)$, which may be taken equal to each other. 
\item
The pruning radii $\delta_d$ used to remove any possible false negative can be set to be equal to the orthogonal slab height $\delta(X_i) = \kappa(X_i) h(X_i)^2$.
\item
Finally, the slab aspect ratios $\kappa_d$ and size cutoffs $n_d$ can be set to $\kappa(X_i) = \mu_d / h(X_i)$ and $n_d =  \varepsilon_d \ell$ with $\varepsilon_d,\mu_d \in (0,1]$.
\end{enumerate}
In the end, note that with this parametrization, $\delta(X_i) = h_\perp(X_i) = \mu_d h(X_i)$.
To jointly select hyperparameters $(\varepsilon_d,\mu_d)$, we may use the following scheme.
At step $d$, recall that $\widetilde{\cX}^{(d)}$ denotes the remaining set of points, which we consider fixed.
Write $\widehat{\cX}_{\mu,\varepsilon}^{(d)}$ for the set of points that is output after step $d$,
and $\widehat{\cZ}_{\mu,\varepsilon}^{(d)} := \bigcup_{X_i \in \widehat{\cX}_{\mu,\varepsilon}^{(d)}} \ball(X_i,h(X_i))\cap \widetilde{\cX}^{(d)}$ for a point-dependent offset of it. We proceed as follows:
\begin{enumerate}
\item
First fixing $\varepsilon_{d,0} := 0.5$, plot the curve of $[0,1] \ni \mu \mapsto  \cG_0(\mu) := \Card \bigl( \widehat{\cX}^{(d)}_{\varepsilon_{d,0},\mu}\bigr)$ and set $\mu_{d,1}$ to correspond to its leftmost elbow value.
\item
Then, set $\varepsilon_{d,1}
:=
\Card \bigl( \widehat{\cX}^{(d)}_{\varepsilon_{d,0},\mu_{d,1}}
\bigr)
/
\Card \bigl( \widehat{\cZ}^{(d)}_{\varepsilon_{d,0},\mu_{d,1}}
\bigr)
$
.
\item
Iterate over $\mu_{d,2},\varepsilon_{d,2},\mu_{d,3},\ldots$ until convergence of the couple $(\mu_{d,t},\varepsilon_{d,t})$.
\end{enumerate}
See Figure~\ref{fig:simu} for an illustration of this heuristic.
This procedure is based on the idea that when fixing $\varepsilon_d$, the function $\mu \mapsto \Card \bigl( \widehat{\cX}^{(d)}_{\varepsilon_d,\mu}\bigr)$ should first increase fast, as a lot of low-dimensional points get included early in slabs of reasonably small aspect ratios. 
Higher-dimensional points should both appear for wider slabs (i.e. larger $\mu$), and be more scattered in space (yielding lower slopes).
Taking the leftmost elbow allows to select a value at the frontier of these two regimes.

Then fixing $\mu_d$, we should expect that a good threshold $n_d = \varepsilon_d \ell$ averages somehow uniformly over the whole stratum. That is, if $d=d_k$ and $\widetilde{P} := \alpha_kP_k + \ldots + \alpha_K P_K$, then for $x \in M_k$,
\begin{align*}
\varepsilon_d
=
\frac{n_d}{\ell}
\simeq
\frac{\widetilde{P}\bigl(\slab(x,h(x),\mu_d h(x))\bigr)}{\widetilde{P}\bigl(\ball(x,h(x))\bigr)}
\simeq
\frac{\widetilde{P}\bigl(\bigcup_{x \in M_k} \slab(x,h(x),\mu_d h(x))\bigr)}{\widetilde{P}\bigl(\bigcup_{x \in M_k} \ball(x,h(x))\bigr)}
\simeq
\frac{
\Card \bigl( \widehat{\cX}^{(d)}_{\varepsilon_d,\mu_d}
\bigr)
}{
\Card \bigl( \widehat{\cZ}^{(d)}_{\varepsilon_d,\mu_d}
\bigr)
}
,
\end{align*}
which yields a self-consistency criterion and hence a natural value for $\varepsilon_{d,t+1}$ given $\mu_{d,t+1}$ and $\varepsilon_{d,t}$.
Furthermore, we observe that if $\varepsilon_{d,t}$ is too large compared to a reasonable value, the subsequent $\mu_{d,t}$ will be too large as well, and hence the next value $\varepsilon_{d,t+1}$ will be smaller. This strategy can hence be seen as a self-stabilizing retroactive loop.
The above heuristics and implementation should be considered as initial practical explorations. A thorough implementation of \codetection, particularly its application to real data, remains an important direction for future work.

\subsection{Further geometric generalizations}

The present work shows that with high probability, a $n$-sample of a mixture of immersed $\cC^2$-manifolds with bounded curvature allows to infer their $d$-dimensional structure at scale $O\bigl( (\log n / n)^{2/d} \bigr)$ in locations and $O\bigl( (\log n / n)^{1/d} \bigr)$ in directions.
We believe that the present formalism of general mixtures of immersed low-dimensional manifolds lies at a sweet spot in terms of statistical modeling.
It combines both the desirable properties induced by $\cC^2$ smoothness allowing the use of standard differential geometry, and the versatility of mixtures and self-intersections induced by general immersions.
Natural extensions of the model of Definition~\ref{def:model-mixture} shall include higher order $\cC^k$ smoothness for faster rates, and the ability to handle immersions of manifolds with boundary and corners, so as to include singular sets containing cusps or open books.
In the latter framework, stratification shall actually emerge via re-sampling and projecting in the ambient space onto the boundary/corners, along the same lines as~\cite{merigot2010voronoi}.

\subsection{About noise}

Even though our model allows for contamination data (see Remark~\ref{rem:clutter-as-ambient-manifold}), an apparent statistical limitation of this work is that it does not allow for any additive noise.
However, it is straightforward to show that all our analysis and results remain true with (possibly adversarial) additive noise of respective amplitudes $O\left( (\log n/n)^{2/d} \right)$ over the stratum of dimension $d$.
Said otherwise, the overall method is stable under small perturbations of data.
See~\cite{Aamari19,Divol21} for examples of such sample-dependent noise frameworks.
Handling additive noise of greater amplitude is a difficult open problem.
Even in the simpler case of a single embedded submanifold with positive reach, the question has not yet completely been understood. So far, fixed-size additive noise has been studied non-computationally in~\cite{Genovese12,Genovese12b,capitao2023support}.
An alternative uniform model was proposed by~\cite{Aizenbud21}, with a denoising process only outputting \emph{one} point close to the underlying submanifold.
Lastly,~\cite{puchkin2022structure} did solve the minimax manifold estimation problem in a computational way, but to the price of requiring a noise amplitude of order at most $O\left( (1/n)^{2/(3d+8)} \right)$.

\subsection{Related estimation tasks}
As pointed out in Section~\ref{sec:stat-setting}, many natural geometric and topological estimation questions are impossible to solve in stratified models.
However, an interesting feat would be to estimate the \emph{distributions} $P_1,\ldots,P_k$ rather than their supports only. For Wasserstein metrics, we believe that the kernel method of~\cite{divol2022measure} may be adapted in the presence of mixed dimensions and self-intersections. However, this would require substantial technical developments that we leave to future work.

\newpage
\appendix

\section{Geometry and volume of immersed manifolds}
\label{sec:appendix-geometry}

Below, immersions $\Psi : \mathbb{M} \to \bbR^D$ are always taken to be isometric in the Riemannian sense (see Section~\ref{sec-euclidean-immersions}).
Since $\bbM$ is compact and without boundary, the inherited Riemannian metric $(\bbM,g)$ yields a globally defined \emph{exponential map} at each point $x_0 \in \bbM$, which we denote by $\exp_{x_0}^{\bbM} : T_{x_0} \bbM \to \bbM$.

\subsection{Geodesic estimates}
\label{sec:geodesic-estimates}

\begin{lemma}[{\cite[Lemma 2.3 (iii)]{tinarrage2023recovering}}] 
\label{lem:geoimm}
Let $\Psi : \bbM \to \bbR^D$ be a $\cC^2$-immersion of a compact $d$-dimensional manifold $\bbM$ such that
$\sup_{x_0 \in \bbM} \|\II_{x_0}^{(\Psi)}\|_{\op} \leq \kappa_{\max}$.
Then for all $x_0,y_0 \in \bbM$ {\color{black}in the same connected component},
\[
\(1-\frac{\kappa_{\max}}{2} \dt_{\bbM}(x_0,y_0) \)\dt_{\bbM}(x_0,y_0) \leq \|\Psi(x_0)-\Psi(y_0)\|.
\]
\end{lemma}

As a result, for all $r < {\color{black}1/(4\kappa_{\max})}$, $\Psi$ is an embedding of the closed geodesic ball $\ball_{\bbM}(x_0,r)$ onto its image. 
{\color{black}The curvature bound gives sectional curvature in $[-2\kappa_{\max}^2,\kappa_{\max}^2]$. The intrinsic injectivity radius is at least $\pi/\kappa_{\max}$: the conjugate-radius bound follows from Rauch comparison, and Klingenberg's estimate together with Fenchel's theorem gives the same lower bound from the lengths of immersed closed geodesics.
}
Hence, if we lift the exponential map $\exp_{x_0}^{\bbM} : T_{x_0} \bbM \to \bbM$ onto $T_{x_0}^{(\Psi)} M := \diff \Psi(x_0) [T_{x_0} \bbM]$ via
\begin{align*}
\overline{\exp}_{x_0}^{\bbM} \colon\;
T_{x_0}^{(\Psi)} M
&\longrightarrow M 
\\
v &\longmapsto \Psi(\exp_{x_0}^{\bbM}(\diff\Psi(x_0)^{-1}[v])),
\end{align*}
we get that $\overline{\exp}_{x_0}^{\bbM}$ is an embedding from the tangent ball $\ball_{T_{x_0}^{(\Psi)} M}(0,r)$ to $\Psi(\ball_{\bbM}(x_0,r))$. 
This embedding fulfills the following estimates.

\begin{lemma}[Geometric approximation bounds]
\label{lem:geoineq}
Let $\Psi : \bbM \to \bbR^D$ be a $\cC^2$-immersion of a compact $d$-dimensional manifold $\bbM$ such that
$\sup_{x_0 \in \bbM} \|\II_{x_0}^{(\Psi)}\|_{\op} \leq \kappa_{\max}$. Fix $x_0 \in \bbM$ and write $x = \Psi(x_0)$.
\bitem
\item[1.] For all $y \in\Psi(\ball_{\bbM}(x_0,1/(4\kappa_{\max})))$, 
$$
\|x-y\| \leq \left\|\{\overline{\exp}_{x_0}^{\bbM}\}^{-1}(y)\right\| \leq 2\|x-y\|,
$$
\item[2.] For all $v \in \ball_{T_{x_0}^{(\Psi)} M}(0,1/(4\kappa_{\max}))$,
$$
\left\| \overline{\exp}_{x_0}^{\bbM}(v)-x-v\right\| \leq \kappa_{\max} \|v\|^2,
$$
\item[3.] For all $y \in \Psi(\ball_{\bbM}(x_0,1/(4\kappa_{\max})))$,
$$
\left\| y-x - \{\overline{\exp}_{x_0}^{\bbM}\}^{-1}(y) \right\| \leq 2\kappa_{\max} \|x-y\|^2,
$$
\item[4.] For all measurable set $A \subset \ball_{T_{x_0}^{(\Psi)} M}(0,1/(4\kappa_{\max}))$,
$$
2^{-d}\cH^d(A) \leq \cH^d(\overline{\exp}^{\bbM}_{x_0}(A)) \leq 2^d \cH^{d}(A).
$$
\item[5.] For all $y \in  \Psi(\ball_{\bbM}(x_0,1/(4\kappa_{\max})))$,
$$
\|y-x-\pr_{T_{x_0}^{(\Psi)} M}(y-x)\| \leq \frac{\kappa_{\max}}{2} \|y-x\|^2.
$$
\item[6.] For all $y_0 \in \ball_{\bbM}(x_0,1/(8\kappa_{\max}))$,
$$
\angle(T_{x_0}^{(\Psi)} M, T_{y_0}^{(\Psi)} M) \leq 5 \kappa_{\max} \|x -{\color{black}\Psi(y_0)}\|.
$$
\eitem
\end{lemma}

\begin{proof} {\color{black}Write $\kappa=\kappa_{\max}$. Integrating the acceleration bound along a geodesic gives
\[
\|\overline{\exp}_{x_0}^{\bbM}(v)-x-v\|\leq\frac\kappa2\|v\|^2.
\]
Together with Lemma~\ref{lem:geoimm} and the length bound, this proves points 1 to 3. For point 4, the sectional-curvature bounds and Jacobi-field comparison bound every singular value of the differential of the exponential map between $1/2$ and $2$ on the indicated ball, as in the proof of~\cite[Lemma 1]{Aamari19}. The area formula then applies, since both the exponential map and $\Psi$ are injective there.

For point 5, let $\gamma:[0,s]\to\bbR^D$ be the image of a minimizing unit-speed geodesic from $x_0$ to a preimage of $y$ in the indicated ball. Thus $s\leq1/(4\kappa)$, $\gamma(0)=x$, $\gamma(s)=y$, and $\|\ddot\gamma\|\leq\kappa$. For any unit normal $n$ to $T_{x_0}^{(\Psi)}M$, and $0\leq t,u\leq s$,
\[
|\langle n,\dot\gamma(t)\rangle|\leq\sin(\kappa t),\qquad
\langle\dot\gamma(t),\dot\gamma(u)\rangle\geq\cos(\kappa|t-u|).
\]
Integrating gives
\[
|\langle n,y-x\rangle|\leq\frac{1-\cos(\kappa s)}\kappa,
\qquad
\|y-x\|^2\geq\frac{2(1-\cos(\kappa s))}{\kappa^2},
\]
which proves point 5.
}
For point 6, we let $y = \Psi(y_0)$ and $v \in \ball_{T_{y_0}^{(\Psi)} M}(0,1/(8\kappa_{\max}))$ and let $z = \overline{\exp}_{y_0}^{\bbM}(v)$, which lies in ${\color{black}\Psi(\ball_{\bbM}(x_0,1/(4\kappa_{\max})))}$ by triangle inequality. 
Using point 5 twice, we first get that
\begin{align*}
\|z-y-\pr_{T_{x_0}^{(\Psi)} M}(z-y)\| \leq \frac{\kappa_{\max}}{2} \{\|z-x\|^2+\|y-x\|^2\}.
\end{align*} 
Furthermore, thanks to point 2, there holds
$$
\left\|z-y-v\right\| \leq \kappa_{\max} \|v\|^2
.
$$
Since $\|z-y\| \leq \|v\|$ from point 1{\color{black}, and since $\mathrm{Id}-\pr_{T_{x_0}^{(\Psi)} M}$ is an orthogonal projection}, we obtain
\begin{align*} 
\|v-\pr_{T_{x_0}^{(\Psi)} M}(v)\| &\leq \frac{\kappa_{\max}}{2} \{\|y-x\|^2+\|z-x\|^2\}+ {\color{black}\kappa_{\max} \|v\|^2} \\
&\leq  \frac{\kappa_{\max}}{2} \{2\|y-x\|^2+\|v\|^2+2\|x-y\|\|v\|\}+ {\color{black}\kappa_{\max} \|v\|^2.}
\end{align*} 
By taking $\|v\| = \|y-x\| \leq 1/(8\kappa_{\max})$, we find
$$
\|v-\pr_{T_{x_0}^{(\Psi)} M}(v)\| \leq 5 \kappa_{\max}\|y-x\|^2 = 5 \kappa_{\max}\|y-x\|\|v\|
.
$$
By symmetry between $x_0$ and $y_0$, we obtain $\angle(T_{x_0}^{(\Psi)} M, T_{y_0}^{(\Psi)} M) \leq 5\kappa_{\max} \|y-x\|$, ending the proof.
\end{proof}

\subsection{Fiber estimates}
\label{sec:sheaf-number}

Under our assumptions, intuition suggests that small enough Euclidean neighborhoods of the image $M = \Psi(\bbM)$ of a $\cC^2$ immersion of a compact manifold $\bbM$ write as a union of $L < \infty$ $d$-dimensional balls up to diffeomorphisms.
The following result provides a quantitative bound on this number $L$, depending on intrinsic dimension, curvature and volume of $(\bbM,\Psi)$.
We write $\omega_d := \cH^d(\ball_d(0,1))$ for the volume of the $d$-dimensional unit Euclidean ball.

\begin{lemma}[Upper bound on local fibers]
\label{lem:sheaf_number}
Let $M = \Psi(\bbM)$ be the image of a $\cC^2$-immersion $\Psi : \bbM \to \bbR^D$ of a compact $d$-dimensional manifold $\bbM$ such that
$\sup_{x_0 \in \bbM} \|\II_{x_0}^{(\Psi)}\|_{\op} \leq \kappa_{\max}$.
Then for all $x \in \bbR^D$ and $r < 1/(4\kappa_{\max})$, there exist $U_1,\ldots,U_L \subset \bbM$ such that
\begin{align*}
\Psi^{-1}(M \cap \ball(x,r)) = \bigcup_{\ell=1}^L U_\ell,
\end{align*}
with:
\begin{itemize}
\item
All the $U_\ell$'s are path-connected and have geodesic diameter at most $4r$.

\item
$L \leq N_0$, where $N_0 := \dfrac{\vol(\bbM)}{\omega_d}(4\kappa_{\max})^d$;

\end{itemize}

\end{lemma}

\begin{proof}[Proof of Lemma~\ref{lem:sheaf_number}]
We let
$$
\Psi^{-1}(M \cap \ball(x,r)) = \bigcup_{U \in \cU} U,
$$
be a path-connected decomposition of $\Psi^{-1}(M \cap \ball(x,r))$ that is maximal for the inclusion. That is, for all other such decomposition 
$$
\Psi^{-1}(M \cap \ball(x,r)) = \bigcup_{V \in \cV} V
,
$$
if for all $U \in \cU$ there exists a $V \in \cV$ such that $U \subset V$, then $\cU = \cV$. 
Such a decomposition exists from Zorn's lemma.

Note that the elements of such a decomposition are pairwise disjoint by minimality and connectedness. We let $U_1, U_2$ be two distinct elements of $\cU$.
We pick a point $x_i$ in each of the $U_i$. {\color{black}If the two points are in different connected components of $\bbM$, their intrinsic distance is infinite. Otherwise, let $\gamma_0$ be a unit-speed shortest path of $\bbM$ going from $x_1$ to $x_2$} and we let $\gamma := \Psi \circ \gamma_0$. The path $\gamma$ is also unit-speed, and because $\gamma_0$ is locally a geodesic, we have that $\|\ddot \gamma\| \leq \kappa_{\max}$ everywhere. 
We know that $\gamma$ has to leave $\ball(x,r)$, otherwise we could replace $U_1$ and $U_2$ by 
the path-connected set
$U_1 \cup {\color{black}\Im\gamma_0} \cup U_2$ 
,
which would contradict the maximality of $\cU$. 
We let $\sigma(t) := \inner{\gamma(t)-x}{\dot\gamma(t)}$. At the first time $t_0$ where $\gamma$ leaves $\ball(x,r)$, we have $\sigma(t_0) \geq 0$, and at the first time $t_1 \geq t_0$ when $\gamma$ returns to $\ball(x,r)$, we must have 
$\sigma(t_1) \leq 0$. 
However, notice that for $t \geq t_0$,
\begin{align*} 
\dot\sigma(t) 
&= 
\|\dot\gamma(t)\|^2
+
\inner{\gamma(t)-x}{\ddot\gamma(t)} 
\\
&\geq 
1 - \kappa_{\max} \|\gamma(t) - x\|
,
\end{align*}
so that for $\sigma$ to turn negative, one must have $\|\gamma(t) - x\| \geq 1/\kappa_{\max}$ for some $t \in [t_0,t_1]$. This implies that
\begin{align*} 
\dt_{\bbM}(x_1,x_2) 
&= {\color{black}\dt_{\bbM}(x_1,\gamma_0(t)) +\dt_{\bbM}(x_2,\gamma_0(t))}
\\
&\geq 
\|\Psi(x_1) - \gamma(t)\| +  \|\Psi(x_2) - \gamma(t)\| 
\\
&\geq 
(\|x - \gamma(t)\| - \|\Psi(x_1) - x\|) 
+
(\|x - \gamma(t)\| - \|\Psi(x_2) - x\|)
\\
&\geq 
2\(\frac{1}{\kappa_{\max}}- r\) 
\\
&\geq 
\frac1{\kappa_{\max}}
.
\end{align*} 
As a consequence, the open geodesic balls $\oball_{\bbM}(x_i,1/(2\kappa_{\max}))$ are mutually disjoint. {\color{black}To obtain the required volume bound at this radius, put $R=1/(2\kappa_{\max})$ and
\[
F(v):=\pr_{T_{x_i}^{(\Psi)}M}(\overline{\exp}_{x_i}^{\bbM}(v)-\Psi(x_i)).
\]
The geodesic acceleration bound gives $\|F(v)-v\|\leq\kappa_{\max}\|v\|^2/2$.
For every $w\in T_{x_i}^{(\Psi)}M$ with $\|w\|\leq1/(4\kappa_{\max})$, the map $v\mapsto w-(F(v)-v)$ sends the closed tangent ball of radius $R$ into its interior, since its norm is at most $3/(8\kappa_{\max})$. Brouwer's fixed point theorem gives $F(v)=w$ with $\|v\|<R$. Projection is $1$-Lipschitz and $\Psi$ is isometric, so
\[
\vol\bigl(\oball_{\bbM}(x_i,1/(2\kappa_{\max}))\bigr)\geq\omega_d(4\kappa_{\max})^{-d}.
\]
} By compactness of $\bbM$, this implies that $\cU$ is of finite cardinality $L$ and that
\begin{align*}
L \omega_d (4\kappa_{\max})^{-d}  \leq \vol(\bbM),
\end{align*}
which yields the bound on $L$.

To show that the $U$'s in $\cU$ have geodesic diameter less than $4r$, let us take $U \in \cU$ and $y,z \in U$ be fixed. 
Because $U$ is path-connected, there exists a continuous path $\gamma : [0,1] \to U$ (not necessarily geodesic) such that $\gamma(0) = y$ and $\gamma(1) =z$.
We let 
$$
\tau := \sup\cT~~~\text{with}~~~\cT := \{t \in [0,1]~|~\dt_{\bbM}(y,\gamma(t)) < 1/\kappa_{\max}\}.
$$
We clearly have $0 \in \cT$, and $\cT$ is open in $[0,1]$ since $\gamma$ is continuous. 
Furthermore, by continuity, $\dt_{\bbM}(y,\gamma(\tau)) \leq 1/\kappa_{\max}$. 
Using \lemref{geoimm}, we find that
\begin{align*}
\frac12 \dt_{\bbM}(y,\gamma(\tau)) \leq \(1-\frac{\kappa_{\max}}{2}  \dt_{\bbM}(y,\gamma(\tau))\)  \dt_{\bbM}(y,\gamma(\tau)) \leq \|\Psi(y)-{\color{black}\Psi(\gamma(\tau))}\| \leq 2r,
\end{align*}
so that $\dt_{\bbM}(y,\gamma(\tau)) \leq 4r < 1/\kappa_{\max}$, meaning that $\tau \in \cT$ which can only be the case if $\tau = 1$. As a result we have $\dt_{\bbM}(y,z) \leq 4r$, hence the result.
\end{proof}

\subsection{Geometric well-posedness of tangents and curvature}
\label{sec:appendix-curvature}

\begin{proof}[Proof of Lemma~\ref{lem:maximal-curvature-identification}]
Let $(\bbM,\Psi) \in \bbI(M)$  and $x \in M$. {\color{black}By compactness and the local embedded-manifold structure at the finitely many points of $\Psi^{-1}(x)$, we may choose $r<1/(16\kappa_{\max})$ sufficiently small that every path-connected component of $\Psi^{-1}(M\cap\ball(x,r))$ contains a preimage of $x$.}
From Lemma~\ref{lem:sheaf_number}, we can write $\Psi^{-1}(M\cap \ball(x,r))$ as a finite disjoint union $\bigcup_{\ell=1}^L U_\ell$ of path-connected $U_1,\dots,U_L \subset \bbM$, each of which is a closed neighborhood of $x_{0,\ell} \in \Psi^{-1}(x)$
with geodesic diameter less than $4r$. In particular, 
$
{\color{black}\Psi(U_\ell)}
\subset 
\overline{\exp}_{x_{0,\ell}}^{\bbM} 
\bigl( \ball_{ T_{x_{0,\ell}}^{(\Psi)} M}(0,4r) \bigr)
$
.

On the other hand if $x_0 \in \Psi^{-1}(x)$ is fixed, then Lemma~\ref{lem:geoimm} asserts that the map
$\Psi \circ \exp_{x_0}^{\bbM} =:\Phi_{x_0}^{\bbM} : \ball_{ T_{x_0} \bbM}(0,4 r) \longrightarrow M$ is a $\cC^2$-diffeomorphism onto its image. 
Furthermore, the first order Taylor expansion of $\Phi_{x_0}^{\bbM}$ nearby $0 \in T_{x_0} \bbM$ writes as
\begin{align*}
\Phi_{x_0}^{\bbM}(v)
&=
x
+
\diff
\Phi_{x_0}^{\bbM}(0)[v]
+
o(\|v\|)
\end{align*}
since $\Phi_{x_0}^{\bbM}(0) = \Psi(x_0) = x$.
Write $V_\ell := \Psi(U_\ell)$. Then for all $y \in V_\ell$,
\begin{align*}
y-x
&=
\diff
\Phi_{x_{0,\ell}}^{\bbM}(0)[\bigl(\Phi_{x_{0,\ell}}^{\bbM}\bigr)^{-1}(y)]
+
o(\|\bigl(\Phi_{x_{0,\ell}}^{\bbM}\bigr)^{-1}(y)\|)
\end{align*}
Furthermore, points 1. and 3. of Lemma~\ref{lem:geoineq}  assert that $\Vert y-x\Vert \leq \Vert \bigl(\Phi_{x_{0,\ell}}^{\bbM}\bigr)^{-1}(y) \Vert \leq (1+2\kappa_{\max} \Vert y-x\Vert) \Vert y-x\Vert$, so that
\begin{align*}
y-x
&=
\diff
\Phi_{x_{0,\ell}}^{\bbM}(0)[\bigl(\Phi_{x_{0,\ell}}^{\bbM}\bigr)^{-1}(y)]
+
o(\|y-x\|)
.
\end{align*}
To show that $T^{(\Psi)}_x M$ does not depend on $(\bbM,\Psi)$, note that by construction,
$
\diff \Psi(x_0) [T_{x_0} \bbM]
=
\diff\Phi_{x_0}^{\bbM}(0) [T_{x_0} \bbM]$.
At first order, the above expansion hence yields
{\color{black}\begin{multline*}
\diff \Psi(x_{0,\ell})[T_{x_{0,\ell}} \bbM]
= 
\{0\}\cup
\{
v \in \bbR^D\setminus\{0\}
\mid
\forall \varepsilon >0
,
\exists y_\varepsilon \in (V_\ell \cap \ball(x,\varepsilon))\setminus \{x\}
\text{~s.t.~}
\left|
\frac{y_\varepsilon-x}{\|y_\varepsilon-x\|}
-
\frac{v}{\|v\|}
\right|
<
\varepsilon
\}
.
\end{multline*}
}
As a consequence, since $M \cap \ball(x,\varepsilon) = \bigcup_{\ell = 1}^L (V_\ell \cap \ball(x,\varepsilon))$ for all $\varepsilon \leq r$,
\begin{align*}
T^{(\Psi)}_x M
:=&
\bigcup_{\ell = 1}^L
\diff \Psi(x_{0,\ell})[T_{x_{0,\ell}} \bbM]
\\
=&
\{0\}\cup
\{
v \in \bbR^D\setminus\{0\}
\mid
\forall \varepsilon >0
,
\exists y_\varepsilon \in (M \cap \ball(x,\varepsilon))\setminus \{x\}
\text{~s.t.~}
\left|
\frac{y_\varepsilon-x}{\|y_\varepsilon-x\|}
-
\frac{v}{\|v\|}
\right|
<
\varepsilon
\}
\end{align*}
does not depend on the chosen parametrisation $(\bbM,\Psi)$.
We note $T_x M := T^{(\Psi)}_x M$ from now on.

To show that $\kappa^{(\Psi)}(M)$ does not depend on $(\bbM,\Psi)$, proceed using the Taylor expansion of~\cite[Lemma B.5]{Aamari19b}, and get
\begin{align*}
\|\II_{x_{0,\ell}}^{(\Psi)}\|_{\op}
&=
\limsup_{\substack{y_{0,\ell} \in \bbM \\ y_{0,\ell} \to x_{0,\ell}}}
\frac{
2 
\dist\left(
\Psi(y_{0,\ell}) - x
\mid
 \diff_{x_{0,\ell}} \Psi [T_{x_{0,\ell}} \bbM]
\right)
}{\Vert \Psi(y_{0,\ell}) - x  \Vert^2}
\\
&=
\limsup_{\substack{y_{\ell} \in V_\ell \\ y_{\ell} \to x }}
\frac{
2 
\dist\left(
y_\ell - x
\mid
 T_\ell
\right)
}{\Vert y_\ell - x  \Vert^2}
,
\end{align*}
where $T_\ell := \diff \Psi(x_{0,\ell}) [T_{x_{0,\ell}} \bbM]$. Hence, writing 
$\kappa^{(\Psi)}(x) := \sup_{\substack{x_0 \in \bbM \\ \Psi(x_0) = x}} \|\II_{x_0}^{(\Psi)}\|_{\op}
$, we get

\begin{align*}
\kappa^{(\Psi)}(x) 
&=
\max_{1 \leq \ell \leq L}
\limsup_{\substack{y_{\ell} \in V_\ell \\ y_{\ell} \to x }}
\frac{
2 
\dist\left(
y_\ell - x
\mid
T_\ell
\right)
}{\Vert y_\ell - x  \Vert^2}
.
\end{align*}
This latter quantity might still depend on the immersion $\Psi$ since the decomposition $V_1,\ldots,V_L$ depend on it.
To show that $\kappa^{(\Psi)}(M) = \sup_{x \in M} \kappa^{(\Psi)}(x)$ does not depend on $\Psi$, we will show that $\kappa^{(\Psi)}(x)$ does not depend on $\Psi$.
For this, note that for all $\ell$ and $y_\ell \in V_\ell$, we have trivially
\begin{align*}
\dist(y_\ell-x \mid T_\ell)
&\geq
\dist\bigl(y_\ell-x \mid \bigcup_{1 \leq \ell' \leq L}
T_{\ell'}\bigr),
\end{align*}
so that we easily get that
$$
\kappa^{(\Psi)}(x) \geq \limsup_{ \substack{y \in M \\ y \to x}} \frac{2\dist(y-x \mid T_x M)}{\|y-x\|^2}. 
$$
We will show the converse inequality to obtain the result. {\color{black}The case $\kappa^{(\Psi)}(x)=0$ is immediate; assume otherwise.} We let $\ell \in \{1,\ldots,L\}$ be such that 
\begin{align}
\label{eq:achieving-limsup-sequence}
\kappa^{(\Psi)}(x) 
=
\limsup_{\substack{y_{\ell} \in V_\ell \\ y_{\ell} \to x }}
\frac{
2 
\dist\left(
y_\ell - x \mid  T_\ell
\right)
}{\Vert y_\ell - x  \Vert^2}.
\end{align}
By extraction of sub-sequences on the unit sphere of $T_\ell$, consider a unit direction $v \in T_{\ell}$ achieving the $\limsup$. That is, there exists $\ve_n \searrow 0$ such that if 
$y_n := \overline{\exp}_{x_{0,\ell}}^{\bbM}(\ve_n v)$, 
then 
\begin{align*}
\dist\left(
y_n - x \mid  T_\ell
\right) 
&=
\frac12 \|y_n-x\|^2 \kappa^{(\Psi)}(x) + o(\|y_n-x\|^2)
\\
&= 
\frac12\ve_n^2 \kappa^{(\Psi)}(x) + o(\ve_n^2)
,
\end{align*}
where we used that $y_n - x = \ve_n v+ O(\ve_n^2)$ from \lemref{geoineq} 2.
We can then find a direction $w \in T_\ell$ such that for all $\eta > 0$ sufficiently small, $v +\eta w \notin T_{\ell'}$ for all $T_{\ell'} \neq T_\ell$: that's because 
$$T_{\ell} \setminus \bigcup_{\ell'~|~T_{\ell'} \neq T_\ell} T_{\ell'}$$
is non empty so we can pick a $w$ in it. If there exists two distinct values $\eta_1,\eta_2$ such that $v +\eta_1 w, v+\eta_2 w \in T_{\ell'}$, then $(\eta_1-\eta_2)w \in T_{\ell'}$ which is a contradiction. Hence, for all $\ell'$, there exists at most one $\eta_{\ell'} > 0$ such that $v+ \eta_{\ell'} w$ is in $T_{\ell'}$ and thus $v + \eta w$ is in none of these subspaces for $\eta < \min_{\ell'} \eta_{\ell'}$. 
We then consider the small perturbation of $y_n$ in direction $w$ given by
$$
z_n := 
\overline{\exp}_{y_n}^{V_\ell} (\ve_n \eta \pr_{T_{y_n} V_{\ell}}(w)),
$$
Thanks to \lemref{geoineq} again, we know that
\begin{align*} 
 \pr_{T_{y_n} V_{\ell}}(w) &= w + O(\ve_n), \\
y_n - x &= \ve_n v+ O(\ve_n^2), \\
z_n - y_n &= \ve_n \eta \pr_{T_{y_n} V_{\ell}}(w) + O(\eta^2 \ve_n^2) \\
&= \ve_n \eta w + 
O(\eta \ve_n^2)
.
\end{align*} 
As a result, $\|z_n-x\| = \ve_n + O(\eta \ve_n)$ and hence
\begin{align*} 
\dist(z_n-x \mid T_\ell)^2 &= \|\pr_{T_\ell^\perp}(z_n-x)\|^2 \\
&= \|\pr_{T_\ell^\perp}(y_n-x)\|^2 + \|\pr_{T_\ell^\perp}(z_n-y_n)\|^2
\\
&\quad
+
2\inner{\pr_{T_\ell^\perp}(z_n-y_n)}{\pr_{T_\ell^\perp}(y_n-x)} 
\\
&= \frac14\kappa^{(\Psi)}(x)^2 \ve_n^4 + O(\eta \ve_n^4){\color{black}+o(\ve_n^4)}.
\end{align*} 
All in all, this gives
\begin{align*}
\frac{2\dist(z_n-x \mid T_\ell)}{\|z_n-x\|^2} 
&=
\frac{\kappa^{(\Psi)}(x) \ve^2_n + O(\eta \ve_n^2){\color{black}+o(\ve_n^2)}}{\ve_n^2+O(\eta \ve_n^2)} 
\\
&= 
\kappa^{(\Psi)}(x) + O(\eta){\color{black}+o(1)}
,
\end{align*}
so that the sequence $z_n \in V_\ell$ achieves the $\limsup$ \eqref{eq:achieving-limsup-sequence} up to $O(\eta)$.
Now for all $\ell'$ such that $T_{\ell'} \neq T_\ell$, there holds,
\begin{align*} 
\dist(z_n-x \mid T_{\ell'})
&=
\| \pr_{T_{\ell'}^\perp}(z_n-x)\| 
=  
\| \ve_n\pr_{T_{\ell'}^\perp}(v + \eta w) +  O(\ve_n^2)\| 
\\
&=
\ve_n \|\pr_{T_{\ell'}^\perp}(v + \eta w)\|
+ O(\ve_n^2)
,
\end{align*}
and because $v+\eta w \notin T_{\ell'}$ for all small $\eta > 0$, there holds that $ \|\pr_{T_{\ell'}^\perp}(v + \eta w)\|  > 0$. 
As a result, for $n$ large enough, for all $\ell'$ such that $T_{\ell'} \neq T_\ell$, we have 
$$
\dist(z_n-x \mid T_{\ell'}) > \dist(z_n-x \mid T_\ell)
,
$$
and in particular
$$
\dist(z_n-x \mid T_x M) = \dist(z_n-x \mid T_\ell),
$$
since $T_x M = \cup_{\ell'=1}^L T_{\ell'}$.
We thus conclude that
\begin{align*} 
\frac{2\dist(z_n-x \mid T_x M)}{\|z_n-x\|^2} 
&= 
\frac{2\dist(z_n-x \mid T_\ell)}{\|z_n-x\|^2} 
=
\kappa^{(\Psi)}(x) + O(\eta),
\end{align*} 
so that for all $\eta>0$,
$$
\limsup_{\substack{ y \in M \\ y \to x}} \frac{2\dist(y-x \mid T_x M)}{\|y-x\|^2} \geq \kappa^{(\Psi)}(x) - O(\eta).
$$
The proof is then complete after taking $\eta \to 0$.
\end{proof}

\subsection{Volume estimates}
\label{sec:volume-estimates}

This section is devoted to prove the heuristic developed in Section~\ref{sec:heuristic-slab-infinite-sample}. To obtain fine bounds, we will need the following notion of generalized angles.

\begin{definition}[Angles between subspaces] 
\label{def:angles}
Let $T,T' \subset \bbR^D$ be linear subspaces of respective dimensions $d \geq d'$. The \emph{principal angles} between $T$ and $T'$ is the sequence of numbers 
$\Theta(T,T') := (\theta_1,\dots,\theta_{d'}) \in [0,\pi/2]^{d'}$ 
defined by
\begin{align*}
\Theta(T,T') := \(\arccos\sigma_1(U^\top U'),\dots,\arccos\sigma_{d'}(U^\top U')\)
,
\end{align*}
where 
\begin{itemize}
\item
$U\in \bbR^{D \times d}$ and  $U'\in \bbR^{D \times d'}$ are orthonormal bases of $T$ and $T'$ respectively;
\item
$\sigma_k(U^\top U')$ is the $k$-th largest singular value of $U^\top U'$.
\end{itemize}
We extend $\Theta(T,T')$ as a vector of size $d$ by setting $\theta_k = \pi/2$ for $k \geq d'+1$. 
\end{definition}
When $\dim T = \dim T'$, the angles $\Theta(T,T')$ are linked to $\angle(T,T') = \Vert\pi_T-\pi_{T'}\Vert_{\mathrm{op}}$ through
\beq \label{eq:anglesin}
\angle(T,T') = \sqrt{1 - \cos^2 \Theta_d(T,T')} = \sin \Theta_d(T,T'),
\eeq
see for instance \cite[Proof of Thm 2.5.1]{golub2013matrix}. When $d > d'$, we have $\angle(T,T') = 1$ and $\Theta_d(T,T') = \pi/2$ so that \eqref{eq:anglesin} still holds for flats of different dimension.

\begin{lemma} \label{lem:angles} Let $T,T' \subset \bbR^D$ be linear subspaces of respective dimensions $d \geq d'$.
For all $h_{\para} \geq h_\perp > 0$, there holds
\begin{align*}
\cH^d\(T \cap \slab_{T'}(0,h_{\para},h_{\perp})\) \leq 4^d \prod_{k=1}^d \{h_{\para} \wedge \frac{h_{\perp}}{\sin\Theta_k(T,T')}\}
,
\end{align*}
where $\Theta(T,T') \in [0,\pi/2]^d$ is the sequence of principal angles between $T$ and $T'$ (Definition~\ref{def:angles}).
\end{lemma}
\begin{proof} Notice that since $h_\perp \leq h_{\para}$, $\slab_{T'}(0,h_{\para},h_{\perp}) \subset \ball(0,2h_{\para})$. Furthermore, there holds
\begin{align*} 
\{v \in T~|~\|\pr_{T'^\perp}(v)\| \leq h_\perp  \} &= \{v \in T~|~\inner{\pr_{T'^\perp}(v)}{\pr_{T'^\perp}(v)} \leq h^2_\perp  \} \\
&=  \{v \in T~|~\inner{v}{\cH v} \leq h^2_\perp  \} 
,
\end{align*} 
where $\cH := \pr_{T} \circ \pr_{T'^\perp} \circ \pr_T$ is a self-adjoint operator of $T$. Let $e_1,\dots,e_d$ be an orthonormal diagonalizing basis for $\cH$, with associated eigenvalues $\lambda_1,\dots,\lambda_d$. Writing $v_i := \inner{e_i}{v}$ for $v\in \bbR^D$, one gets
\begin{align*} 
\cH^d \{v \in \ball(0,2h_{\para}) \cap T~|~\inner{v}{\cH v} \leq h^2_\perp  \} &= \cH^d \{v \in \ball(0,2h_{\para})\cap T~|~\sum_{i=1}^d \lambda_i v_i^2 \leq h^2_\perp  \} \\
&\leq \cH^d \{v \in \ball(0,2h_{\para})\cap T~|~\forall i,~v_i^2 \leq h^2_\perp/\lambda_i  \} \\
&\leq 4^d \prod_{i=1}^d \(\frac{h_\perp}{\sqrt{\lambda_i}} \wedge h_{\para}\).
\end{align*} 
To obtain an explicit expression of the $\lambda_i$'s, let $U \in \bbR^{D\times d}$ and $V\in \bbR^{D\times d'}$ be orthonormal bases of $T$ and $T'$ respectively.
Then there holds
\begin{align*}
\cH 
&= 
UU^\top \(\Id_D - VV^\top\) UU^\top 
\\
&= U(\Id_d - (U^\top V)(U^\top V)^\top) U^\top
,
\end{align*}
so that as an endomorphism of $T$, the spectrum of $\cH$ is given by
$$
\spec\cH = 1 - \spec \{(U^\top V)(U^\top V)^\top\} = 1 - \cos^2 \Theta (T,T') = \sin^2\Theta(T,T'),
$$
which ends the proof.
\end{proof} 
From Lemma~\ref{lem:angles}, we deduce the following immediate corollary. 

\begin{corollary} \label{cor:slab_vol1}
Let $T,T' \subset \bbR^D$ be linear subspaces of respective dimensions $d$ and $d'$.
For all $h_{\para} \geq h_\perp > 0$, there holds
\bitem
\item[(i)] If $d = d'$, then
$$ 
\cH^d\(T \cap \slab_{T'}(0,h_{\para},h_{\perp})\) \leq 4^d h_{\para}^{d-1} \{\frac{h_\perp}{\angle(T,T')}\wedge h_{\para}\}. 
$$
\item[(ii)] If $d \geq d'$, then
$$
\cH^d\(T \cap \slab_{T'}(0,h_{\para},h_{\perp})\) \leq 4^d h_{\para}^{d'} h_{\perp}^{d-d'}. 
$$
\eitem
\end{corollary}

\begin{proof}
This is a direct application of Lemma~\ref{lem:angles} and Definition~\ref{def:angles}.
\end{proof}

Corollary~\ref{cor:slab_vol1} extends to the non-linear case in the following sense. 
\begin{lemma}
\label{lem:slab_vol2}
Let $M \in \modelM$ and $T \subset \bbR^{D}$ be a linear subspace of dimension $d' \leq d$. 
Take $h_{\para} < 1/(8\kappa_{\max})$ and
$h_{\perp} := \eta \kappa_{\max}h_{\para}^2$
for some constant $0 < \eta  < 1/(\kappa_{\max} h_{\para})$.
Then for all $x \in \bbR^D$, 
\bitem
\item[(i)] If $d' = d$, there holds
\begin{align*} \cH^{d}\bigl( M \cap \slab_T(x,h_{\para},h_\perp) \bigr) \leq C_d V_{\max} \kappa_{\max}^d h_{\para}^{d-1}
\times &\{\frac{(1+\eta)\kappa_{\max }h_{\para}^2}{\theta}\wedge h_{\para}\} 
,
\end{align*} 
where ${\color{black}\theta := \displaystyle \inf_{y \in M \cap \slab_T(x,h_{\para},h_\perp)} \angle{(T \mid T_y M)}}$.

\item[(ii)] If $d' < d$, there holds
$$
\cH^{d}\bigl( M \cap \slab_T(x,h_{\para},h_\perp) \bigr) \leq C_d (1+\eta)^{d-d'} V_{\max} \kappa_{\max}^{2d-d'} h_{\para}^{2d-d'}.
$$
\item[(iii)] If $d' = d$ and $\cH^{d}\bigl( M \cap \slab_T(x,h_{\para},h_\perp) \bigr) \geq \beta h_{\para}^d$ for some $\beta > 0$, then there exists $y \in M$ such that
\begin{align*} 
\angle( T \mid T_yM ) &\leq C_d \(1+\beta^{-1} V_{\max}\kappa_{\max}^{d}\)(1+\eta)\kappa_{\max} h_{\para}
 \\
\text{and}~~~~~~~~~~~~~\|x - y\| &\leq C_d \(1+\beta^{-1} V_{\max}\kappa_{\max}^{d}\)(1+\eta)\kappa_{\max} h_{\para}^2.
\end{align*}  
\eitem
\end{lemma}

\begin{proof} Let $(\bbM,\Psi) \in \bbI(M)$ be such that $\vol(\bbM) \leq 2 V_{\max}$. {\color{black}It suffices to prove the result for $h_{\para}<1/(64\kappa_{\max})$. For $1/64\leq\kappa_{\max}h_{\para}<1/8$, (i) and (ii) follow from $\cH^d(M)\leq2V_{\max}$ and $\theta\leq1$, after enlarging $C_d$; for (iii), any point of the nonempty intersection has angle at most $1$ and lies within $2h_{\para}$ of $x$, which gives both bounds after the same adjustment.} 
Since $h_{\perp} \leq h_{\para}$, \lemref{sheaf_number} yields
$$
M \cap \slab_T(x,h_{\para},h_\perp) \subset M \cap \ball(x,2h_{\para}) = \bigcup_{j=1}^N \Psi(U_j)
$$
with $\diam_{\bbM}(U_j) \leq 8h_{\para}$ and $N \leq N_0 := \vol(\bbM)(4\kappa_{\max})^d/\omega_d$.
A union bound hence yields, letting $V_j = \Psi(U_j)$,
\begin{align}
\cH^{d}\bigl( M \cap \slab_T(x,h_{\para},h_\perp) \bigr) 
\notag
&\leq 
\sum_{j=1}^N \cH^d\( V_j \cap \slab_T(x,h_{\para},h_\perp) \)
\\
&\leq
N_0 \max_{1 \leq j \leq N} 
\cH^d\( V_j \cap \slab_T(x,h_{\para},h_\perp) \)
. 
\label{eq:master-volume-bound}
\end{align}
If $V_j \cap \slab_T(x,h_{\para},h_\perp) = \emptyset$ then $ \cH^d\( V_j \cap \slab_T(x,h_{\para},h_\perp) \)=0$. Otherwise, take $x_j \in V_j \cap \slab_T(x,h_{\para},h_\perp)$, and use triangle inequality to get
\begin{align*}
V_j \cap \slab_T(x,h_{\para},h_\perp) \subset V_j \cap \slab_T(x_j,2h_{\para},2h_\perp).
\end{align*}
Let us write
$x_{0,j} 
:=
(\Psi_{\vert U_j})^{-1}(x_j)$.
From \lemref{geoineq} 3., we have that for all $y \in V_j \cap \slab_T(x_j,2h_{\para},2h_\perp)$,  the tangent vector $v = (\overline{\exp}_{x_{0,j}}^{\bbM})^{-1}(y) \in T_{x_{0,j}}^{(\Psi)} M$ satisfies 
$$
{\color{black}\|v-(y-x_j)\|} \leq 2\kappa_{\max}\|x_j-y\|^2  \leq 32\kappa_{\max}h_{\para}^2,
$$
where the last bound comes from
$\slab_T(x_j,2h_{\para},2h_\perp) \subset \ball(x_j,4h_{\para})$.
By triangle inequality, we hence get that 
$$
v \in 
\slab_T(0,2h_{\para} + 32\kappa_{\max}h_{\para}^2,2h_\perp + 32\kappa_{\max}h_{\para}^2) 
\subset
\slab_T(0,6h_{\para},(2\eta + 32)\kappa_{\max}h_{\para}^2)
.
$$
In summary, we have proven that
\begin{align*} 
V_j \cap \slab_T(x_j,2h_{\para},2h_\perp) 
&\subset \overline{\exp}_{x_{0,j}}^{\bbM}\({\color{black}\ball_{T_{x_{0,j}}^{(\Psi)}M}(0,8h_{\para})} \cap \slab_T(0,6h_{\para},(2\eta+32)\kappa_{\max}h_{\para}^2)\).
\end{align*} 

To prove point (i), assume that $d = d'$ and use \lemref{geoineq} 4. and {\color{black}Corollary~\ref{cor:slab_vol1}} (i) to get
\begin{align} 
\cH^d(V_j \cap \slab_T(x,h_{\para},h_\perp)) 
&\leq 
\cH^d\({\color{black}\overline{\exp}}_{x_{0,j}}^{\bbM}\({\color{black}\ball_{T_{x_{0,j}}^{(\Psi)}M}(0,8h_{\para})}
\cap 
\slab_T(0,6h_{\para}, (2\eta+32)\kappa_{\max}h_{\para}^2)\) 
\) \notag
\\
&\leq 
2^d \cH^d\(T_{x_{0,j}}^{(\Psi)} M \cap \slab_T(0,6h_{\para}, (2\eta+32)\kappa_{\max}h_{\para}^2)\) \label{eq:comp}\\
&\leq 48^d h_{\para}^{d-1}
\times \{\frac{(2\eta+32)\kappa_{\max}h_{\para}^2}{\angle{(T_{x_{0,j}}^{(\Psi)} M,T)}}\wedge 6h_{\para}\},
\notag
\end{align} 
The last inequality rewrites as
\beq
\cH^d(V_j \cap \slab_T(x,h_{\para},h_\perp)) \leq C_d h_{\para}^{d-1}\{\frac{(\eta+1)\kappa_{\max}h_{\para}^2}{\angle{(T_{x_{0,j}}^{(\Psi)} M,T)}}\wedge h_{\para}\}
,
\label{eq:volTx}
\eeq
which combined with \eqref{eq:master-volume-bound} concludes the proof of (i). 

For point (ii), use point (ii) of {\color{black}Corollary~\ref{cor:slab_vol1}} in the bound \eqref{eq:comp} to get
\begin{align*}
\cH^d(V_j \cap \slab_T(x,h_{\para},h_\perp)) &\leq 8^d (6h_{\para})^{d'}\((2\eta+32)\kappa_{\max}h_{\para}^2\)^{d-d'} 
\\
&\leq C_d (\eta+1)^{d-d'} \kappa_{\max}^{d-d'} h_{\para}^{2d-d'},
\end{align*} 
and combine it with \eqref{eq:master-volume-bound}. 

For point (iii), item (i) and its proof yield the existence of some $j \leq N$ such that $\cH^d(V_j \cap \slab_T(x,h_{\para},h_\perp)) \geq \beta h_{\para}^d/N_0$. If $(\eta+1)\kappa_{\max}h_{\para}^2/\angle{(T_{x_{0,j}}^{(\Psi)} M,T)} \geq h_{\para}$, then $\angle{(T_{x_{0,j}}^{(\Psi)} M,T)} \leq (\eta+1) \kappa_{\max} h_{\para}$. If not, then using \eqref{eq:volTx}, we get that 
$$\angle{(T_{x_{0,j}}^{(\Psi)} M, T)} \leq C_d V_{\max}\kappa_{\max}^{d+1} \frac{\eta+1}{\beta} h_{\para}.
$$
Furthermore, writing $y := \overline{\exp}_{x_{0,j}}^{\bbM}(v)$ with $v := \pr_{T_{x_{0,j}}^{(\Psi)} M}(x - x_j)$, we get that 
$\|y - x_j - v \| \leq 2\kappa_{\max}\|v\|^2$. 
In addition, 
\begin{align*} 
\|v -  (x - x_j) \| &\leq  \|\pr_{T_{x_{0,j}}^{(\Psi)} M}(\pr_{T}(x-x_j)) - \pr_{T}(x - x_j)\| + 2h_{\perp} \\
&\leq \angle{(T_{x_{0,j}}^{(\Psi)} M, T)}  \times \| \pr_{T}(x-x_j) \| + 2h_{\perp} \\
&\leq \(\beta^{-1}C_d V_{\max}\kappa_{\max}^{d}+1\)(\eta+1)\kappa_{\max} h^2_{\para} + 2 \eta \kappa_{\max} h_{\para}^2
\end{align*} 
which together with $\|v\|^2 \leq 8h_{\para}^2$ yields 
$$\|y-x\| \leq C_d \(\beta^{-1} V_{\max}\kappa_{\max}^{d}+1\)(\eta+1)\kappa_{\max} h^2_{\para}$$ 
Finally, using \lemref{geoineq} 6, we find that
\begin{align*}
\angle{(T|T_y M)} 
&\leq 
\angle{(T_{x_{0,j}}^{(\Psi)} M, T)} + 5 \kappa_{\max}\|y-x_j\| 
\\
&\leq 
C_d \(\beta^{-1} V_{\max}\kappa_{\max}^{d}+1\)(\eta+1)\kappa_{\max} h_{\para}
,
\end{align*}
ending the proof.
\end{proof}

Point (ii) of \lemref{slab_vol2} takes the simple, weaker form,
\beq \label{eq:simple}
\cH^d\(M \cap \slab_{T}(x,h_{\para},h_{\perp})\) \leq C_d (\eta+1)^{d} V_{\max} \kappa^{d+1}_{\max}h_{\para}^{d+1},
\eeq
which corresponds to the limiting case where $d = d'+1$ and which will be of use in the forthcoming proofs.

\begin{lemma}
\label{lem:slab_vol3}
Let $M \in \modelM$, $x \in M$, and $T\subset \bbR^D$ be a linear subspace of dimension $d$. 
Take $h_{\para} < 1/(8\kappa_{\max})$ and
$h_{\perp} := \eta \kappa_{\max}h_{\para}^2$
for some constant $0 < \eta  < 1/(\kappa_{\max} h_{\para})$.
If $\angle(T|T_x M) \leq b \kappa_{\max} h_{\para}$ with $b \geq \sqrt{\eta}$, then 
$$
\cH^{d}\bigl( M \cap \slab_T(x,h_{\para},h_\perp) \bigr) 
\geq
c_d((\eta/b) \wedge 1)^d h_{\para}^d. 
$$
\end{lemma}

\begin{proof} Writing $\gamma := (\eta/(2b)) \wedge 1$, let us show that there exists $x_0 \in \Psi^{-1}(x)$ such that the inclusion $\overline{\exp}_{x_0}^{\bbM}(T_{x_0}^{(\Psi)} M \cap \ball(0, \gamma h_{\para})) \subset M \cap \slab_T(x,h_{\para},h_\perp)$ holds. 
Combined with the lower bound of \lemref{geoineq} 4., this inclusion will be sufficient to conclude.

To prove the announced inclusion, take $x_0 \in \Psi^{-1}(x)$ such that $\angle(T,T_{x_0}^{(\Psi)} M) \leq b \kappa_{\max} h_{\para}$.
Given arbitrary $v \in T_{x_0}^{(\Psi)} M \cap \ball(0, \gamma h_{\para})$ write $y := \overline{\exp}^{\bbM}_{x_0}(v)$. 
Thanks to \lemref{geoineq} 1., we have $\|y-x\| \leq \|v\| \leq \gamma h_{\para} \leq h_{\para}$. Furthermore, \lemref{geoineq} 5. yields that
\begin{align*} 
\|\pr_{T_{x_0}^{(\Psi)} M^\perp}(y-x) \| &= \|y-x - \pr_{T_{x_0}^{(\Psi)} M}(y-x) \| \leq \frac{\kappa_{\max}}{2}\|y-x\|^2 \leq \frac{\kappa_{\max}}{2} \eta h_{\para}^2
,
\end{align*} 
where we used that $\gamma^2 \leq \eta^2/b^2 \leq \eta$ since $b \geq \sqrt{\eta}$.
Thus,
\begin{align*} 
\|\pr_{T^\perp}(y-x) \| &\leq \|\pr_{T_{x_0}^{(\Psi)} M^\perp}(y-x)  \| + \angle(T_{x_0}^{(\Psi)} M,T) \|y-x\| 
\\
&\leq   
\frac{\kappa_{\max}}{2} \eta h_{\para}^2 + b \kappa_{\max} \gamma h^2_{\para} 
\\
&\leq  
\eta\kappa_{\max} h_{\para}^2 
\\
&= 
h_{\perp}
,
\end{align*} 
which concludes the proof.
\end{proof}
\section{Concentration bounds}
\label{sec:appendix-concentration}

To link slab counting and (integrated) mass at the population level, we will use extensively the following concentration bound for the $K$ empirical distributions $P_{1,n},\ldots,P_{K,n}$ (see Section~\ref{sec:stat-setting} for notation).

\begin{lemma}
\label{lem:slab_concentration}
Let $\mathcal{S}$ denote the set of all slabs of $\bbR^D$:
\begin{align*}
\mathcal{S}
:=
\{
\slab_T(x,h_{\parallel},h_{\perp})
\mid
x \in \bbR^D,\ T \subset \bbR^D \text{~linear, and~}
h_{\parallel},h_{\perp} \in \bbR_+
\}.
\end{align*}
For $n$ so that
{\color{black}$n\alpha_{\min}\geq100D^2\log D$ and }
$(800D^2\log D)\frac{\log(n\alpha_{\min})}{n\alpha_{\min}}\leq1$,
and for all
$\sqrt{(800D^2\log D)\frac{\log(n\alpha_{\min})}{n\alpha_{\min}}}
\leq\ve\leq1$,
\begin{align*}
\bbP\(
\max_{k\in\{1,\ldots,K\}}
\sup_{S\in\mathcal{S}}
\frac{P_k(S)-P_{k,n}(S)}{\sqrt{P_k(S)}}
\vee
\frac{P_{k,n}(S)-P_k(S)}{\sqrt{P_{k,n}(S)}}
\geq\ve
\)
\leq
4K\exp\(-\ve^2n\alpha_{\min}/16\).
\end{align*}
As a result, if we let
\begin{align*}
\ve
:=
\(\gamma\frac{\log n}{n}\)^{1/2}
\qquad\text{with}\qquad
\gamma
\geq
\frac{16q}{\alpha_{\min}}
\vee
\frac{800D^2\log D}{\alpha_{\min}}
\end{align*}
for some $q>0$, then for $n$ large enough so that $\ve\leq1$,
the probability of the event
\begin{align*}
\cE_1
:=
\left\lbrace
\max_{k\in\{1,\ldots,K\}}
\sup_{S\in\mathcal{S}}
\frac{P_k(S)-P_{k,n}(S)}{\sqrt{P_k(S)}}
\vee
\frac{P_{k,n}(S)-P_k(S)}{\sqrt{P_{k,n}(S)}}
<\ve
\right\rbrace
\end{align*}
is at least $1-4Kn^{-q}$.
\end{lemma}

The proof is based on a Vapnik-Chernovenkis argument on the class of slabs, combined with concentration of multinomial distributions.

\begin{proof}
Write $\rho_0(t):=\ind_{t\geq0}$ for the Heaviside function on $\bbR$,
and $\rho_2(t):=t^2$ for the square function.
Let $T\subset\bbR^D$ be a linear subset of dimension
$d\in\{{\color{black}0,}1,\ldots,D\}$.
If $e_1,\ldots,e_d$ denotes an orthonormal basis of $T$
and $e_{d+1},\ldots,e_D$ an orthonormal basis of $T^\perp$,
then we have
\begin{align*}
&\ind_{\slab_T(x,h_{\parallel},h_{\perp})}(z)=1
\\
&\Leftrightarrow
\begin{cases}
\|\pi_T(z-x)\|^2\leq h_{\parallel}^2,\\
\|\pi_{T^\perp}(z-x)\|^2\leq h_{\perp}^2,
\end{cases}
\\
&\Leftrightarrow
\begin{cases}
\displaystyle
\sum_{j=1}^{d}\rho_2(\langle e_j,z-x\rangle)
\leq h_{\parallel}^2,\\
\displaystyle
\sum_{j=d+1}^{D}\rho_2(\langle e_j,z-x\rangle)
\leq h_{\perp}^2,
\end{cases}
\\
&\Leftrightarrow
\begin{cases}
\displaystyle
\rho_0\left(
h_{\parallel}^2-
\sum_{j=1}^{d}\rho_2(\langle e_j,z-x\rangle)
\right)=1,\\
\displaystyle
\rho_0\left(
h_{\perp}^2-
\sum_{j=d+1}^{D}\rho_2(\langle e_j,z-x\rangle)
\right)=1,
\end{cases}
\\
&\Leftrightarrow
\rho_0\left(
h_{\parallel}^2-
\sum_{j=1}^{d}\rho_2(\langle e_j,z-x\rangle)
\right)
+
\rho_0\left(
h_{\perp}^2-
\sum_{j=d+1}^{D}\rho_2(\langle e_j,z-x\rangle)
\right)
-2
\geq0.
\end{align*}

{\color{black}
This formulation shows that slab membership is determined by two
quadratic inequalities. Introducing the fixed feature map
\[
\Phi(z)
:=
\bigl(z_1,\ldots,z_D,(z_i z_j)_{1\leq i\leq j\leq D}\bigr)
\in\bbR^{m_D},
\qquad
m_D:=\frac{D(D+3)}2,
\]
these inequalities become affine inequalities in $\Phi(z)$.
Their conjunction can be represented by two threshold units with
activation $t\mapsto\ind_{t>0}$ recording violations of the respective
inequalities, followed by the linear output $1/2-u_1-u_2$.
Its sign is positive exactly when both inequalities hold.

Momentarily using notation from~\cite[Theorem~7]{bartlett2019nearly},
this network has $p=1$, $d=0$, $L=2$, $U=3$, and
$W=2(m_D+1)+3=D^2+3D+5$ parameters.
Since composition with a fixed feature map cannot increase
the VC dimension, we obtain
\begin{align*}
\VC(\mathcal{S})
&\leq
2+(D^2+3D+5)\log_2\bigl(12e\log_2(6e)\bigr)
\\
&\leq
50D^2\log D,
\end{align*}
where the last two inequalities hold for all $D\geq2$.
}

On the other hand, for all $k\in\{1,\ldots,K\}$,
conditionally on $N_k${\color{black}, with $N_k\geq v$},
\[
P_{k,n}
:=
\frac1{N_k}\sum_{i=1}^{n}\ind_{L_i=k}\delta_{X_i}
\]
is the empirical distribution of an i.i.d. sample
$X'_1,\ldots,X'_{N_k}$ with common distribution $P_k$.
Hence, from~\cite[Theorem~5.1]{boucheron2005theory}
and Sauer's lemma, we get that for all $\ve>0$,
\begin{align*}
\bbP\(
\sup_{S\in\mathcal{S}}
\frac{P_k(S)-P_{k,n}(S)}{\sqrt{P_k(S)}}
\geq\ve
\middle|N_k
\)
\leq
4\(\frac{2eN_k}{v}\)^v
\exp\(-\ve^2N_k/4\),
\end{align*}
where $v:={\color{black}50}D^2\log D$.
In particular, for all
$\ve\geq\sqrt{({\color{black}400}D^2\log D)\frac{\log N_k}{N_k}}$,
\begin{align*}
\bbP\(
\sup_{S\in\mathcal{S}}
\frac{P_k(S)-P_{k,n}(S)}{\sqrt{P_k(S)}}
\geq\ve
\middle|N_k
\)
\leq
\exp\(-\ve^2N_k/8\).
\end{align*}
{\color{black}
Indeed, $\log4+v\log(2e/v)\leq0$ and
$\ve^2N_k\geq8v\log N_k$.
Moreover, on $N_k>n\alpha_k/2$, we have
$N_k>n\alpha_{\min}/2\geq v$, and the lower bound on $\ve$
in the statement implies this conditional threshold, since
\[
\frac{\log N_k}{N_k}
\leq
\frac{2\log(n\alpha_{\min})}{n\alpha_{\min}}
\]
for the indicated range of $n$.
}

To obtain an unconditional bound, we apply Bernstein's inequality
\cite[Corollary~2.11 and~(2.10)]{Massart13}
to the Binomial variable $N_k$ with parameters $n$ and $\alpha_k>0$,
which yields
\begin{align*}
\bbP\(N_k\leq n\alpha_k/2\)
\leq
\exp(-3n\alpha_k/28).
\end{align*}
As a result,
\begin{align*}
\bbP
&\(
\sup_{S\in\mathcal{S}}
\frac{P_k(S)-P_{k,n}(S)}{\sqrt{P_k(S)}}
\geq\ve
\)
\\
&\leq
\bbP\(N_k\leq n\alpha_k/2\)
+
\bbP\(
\sup_{S\in\mathcal{S}}
\frac{P_k(S)-P_{k,n}(S)}{\sqrt{P_k(S)}}
\geq\ve
\middle|N_k>n\alpha_k/2
\)
\\
&\leq
\exp(-3n\alpha_k/28)
+
\exp\(-\ve^2n\alpha_k/16\)
\\
&\leq
2\exp\(-\ve^2n\alpha_k/16\),
\end{align*}
as soon as $\ve\leq1$ and
$\ve\geq
\sqrt{(800D^2\log D)\frac{\log(n\alpha_{\min})}{n\alpha_{\min}}}$.

A union bound over $k\in\{1,\ldots,K\}$ allows to conclude that
\begin{align*}
\bbP\(
\max_{k\in\{1,\ldots,K\}}
\sup_{S\in{\color{black}\mathcal{S}}}
\frac{P_k(S)-P_{k,n}(S)}{\sqrt{P_k(S)}}
\geq\ve
\)
\leq
2K\exp\(-\ve^2n\alpha_{\min}/16\).
\end{align*}
To conclude the proof, do the exact same reasoning for the relative
deviation $(P_{k,n}(S)-P_k(S))/\sqrt{P_{k,n}(S)}$
by using the other inequality of~\cite[Theorem~5.1]{boucheron2005theory}.
\end{proof}

We deduce the following thresholded version of Lemma~\ref{lem:slab_concentration}, which roughly asserts that slabs $S$ with $d_k$-dimensional mass $P_k(S)$ of order $\log n/n$ can be distinguished through $P_{k,n}(S)$.

\begin{corollary} \label{cor:conc}On the event $\cE_1$ of Lemma~\ref{lem:slab_concentration}, there holds, for all $c > 0$,
\benum
\item If $P_{k}(S) \geq c^2 \ve^2$, then $P_{k,n}(S) \geq \(1-1/c\) P_k(S)$.
\item If $P_k(S) \leq c^2 \ve^2$, then $P_{k,n}(S) \leq (1+c)^2 \ve^2$.
\eenum
\end{corollary}

\begin{proof} If $P_k(S) \geq c^2\ve^2$, then
$$
P_{k,n}(S) \geq P_k(S) - \ve \sqrt{P_k(S)} \geq \(1-1/c\) P_k(S).
$$
Reciprocally, if $P_k(S) \leq c^2 \ve^2$, then $P_{k,n}(S) - P_k(S) \leq  \ve \sqrt{P_{k,n}(S)}$ which implies
$$
\sqrt{P_{k,n}(S)} \leq \frac12\{\sqrt{\ve^2+4 P_k(S)}+\ve\} \leq \frac12\(\sqrt{1+4c^2}+1\)\ve \leq (c+1)\ve, 
$$
which ends the proof.
\end{proof}

Similarly, the following concentration bound asserts that enough geodesic balls of radii $(\log n/n)^{1/d_k}$ all contain at least approximately $\log n$ points drawn from the $d_k$-dimensional layer.

\begin{lemma} \label{lem:eventcE} 
{\color{black}Choose the immersions $(\bbM_k,\Psi_k)$ so that $\vol(\bbM_k)\leq2\nu_{\max}\kappa_{\max}^{-d_k}$.}
Given $k \in \{1,\ldots,K\}$ and $\rho_{d_k} >0$, let $\cP_{k}$ denote a maximal $(2\rho_{d_k})$-packing of $\bbM_k$, in the sense that ${\color{black}\dt_{\bbM_k}(p,p')} > 2\rho_{d_k}$ for all $p \neq p' \in \cP_k$.
Let $\cB_{k}$ denote the class of subsets of $\bbR^D$ defined by
\begin{align*}
\cB_{k} := \{\Psi_k(\ball_{\bbM_k}(p,\rho_{d_k}))~|~p\in \cP_{k}\}.
\end{align*}
Let $q >0$, $\upsilon \geq 0$ and
$$\ve := \(\gamma \frac{\log n}{n}\)^{1/2}~\text{and}~~\rho_d := \(\frac{\Upsilon_d\ve^2}{a_{\min}\kappa_{\max}^d}\)^{1/d}~~\text{with}~~\gamma \geq 
\frac{56{\color{black}(1\vee q)}}{3 \alpha_{\min}}
~~\text{and}~~\Upsilon_d \geq 2(1\vee \upsilon) \frac{2^d}{\omega_d}.$$ 
Then for $n$ large enough, the probability of the event
$$
\cE_2 := 
{\color{black}\{\min_{1 \leq k \leq K} N_k/(n\alpha_k)>1/2\}\cap}
\{
\min_{{\color{black}1 \leq k \leq  K}}
\min_{B \in \cB_{k}}
P_{k,n}(B)
\geq \upsilon \ve^2 \}
$$
is at least $1 - 2 K n^{-q}$.
\end{lemma}

\begin{proof} 
Let first $k\in \{1,\ldots,K\}$ and $B \in \cB_k$ be fixed{\color{black}, and write $G_k:=\{N_k>n\alpha_k/2\}$}.
As soon as $\ve^2 \leq a_{\min}/(4^{d_k}\Upsilon_{d_k})$,
\lemref{geoineq} 4. applies and yields
\begin{align}
\label{eq:ineq-mass-covering}
P_k(B)
&\geq 
a_{\min} \kappa_{\max}^{d_k} 2^{-d_k} \omega_{d_k}\rho_{d_k}^{d_k} 
\notag
\\
&= 
2^{-d_k} \omega_{d_k} \Upsilon_{d_k} \ve^2 
\notag
\\
&\geq 
2(1 \vee \upsilon) \ve^2.
\end{align}
Therefore, from Bernstein's bound for binomial distributions, we have
\begin{align*} 
\bbP\(P_{k,n}(B) \leq \frac12 P_k(B) \middle| N_k\) &\leq \exp\(-\frac3{28} N_k P_k(B)\) \leq \exp\(-\frac{3}{14} N_k \ve^2\).
\end{align*}
 Using the same bound as in the proof of \lemref{slab_concentration}, we proceed with
\begin{align*} 
\bbP\({\color{black}G_k^c\cup\lbrace}P_{k,n}(B) \leq \frac12 P_k(B){\color{black}\rbrace}\) &\leq 
\bbP(N_k \leq n\alpha_k/2)+\bbP\(P_{k,n}(B) \leq \frac12 P_k(B) \middle| N_k > \alpha_k n/2\) 
\\
&\leq 
\exp\(-3n\alpha_k/28\) + \exp\(-3 \alpha_k n\ve^2/28\)
\\
&\leq 
2\exp\(-\frac{3}{28}\alpha_k n\ve^2\),
\end{align*} 
for $n$ large enough such that $\ve \leq 1$. 
On the other hand, the cardinal of $\cB_k$ can be bounded from above thanks to \lemref{geoineq} 4. by
 \begin{align*} 
 \Card(\cB_k)
 &\leq
 \frac{\vol \bbM_k}{2^{-{d_k}} \omega_{d_k} \rho_{d_k}^{d_k}} 
 \leq
\frac{
 2 \nu_{\max}/\kappa_{\max}^{d_k}
 }{
 2^{-{d_k}} \omega_{d_k} 
 \(
 \frac{\Upsilon_{d_k} \ve^2}{a_{\min} \kappa_{\max}^{d_k}}
 \)
 }
 \leq
 \frac{{\color{black}a_{\min}}\nu_{\max} }{\ve^2}  
 \leq 
 \exp\(\frac{3}{56} \alpha_k n\ve^2\),
 \end{align*} 
where the last inequality holds for $n$ large enough.
As a result, for $n$ large enough, a union bound over the elements of $\cB_k$ yields
\begin{align*}
 \bbP\({\color{black}G_k^c\cup\lbrace}\min_{B \in \cB_k} \frac{P_{k,n}(B)}{P_k(B)} \leq \frac12{\color{black}\rbrace}\) 
 &\leq 
 2\exp\(-\frac{3}{28} \alpha_k n\ve^2\)
 \Card(\cB_k)
 \\
 &\leq
 2\exp\(-\frac{3}{56} \alpha_k n\ve^2\).
\end{align*}
By a union bound over $k \in \{1,\ldots,K\}$, we conclude that
\begin{align*}
\bbP\({\color{black}\bigcup_{k=1}^K G_k^c\cup\lbrace}\min_{1 \leq k \leq K}\min_{B \in \cB_k} \frac{P_{k,n}(B)}{P_k(B)} \leq \frac12{\color{black}\rbrace}\) \leq 2K\exp\(-\frac{3}{56} \alpha_{\min} n\ve^2\) \leq 2Kn^{-q}
,
\end{align*}
{\color{black}Retaining the polynomial packing bound before its exponential relaxation also gives, for later use,
\begin{align*}
\bbP(\cE_2^c)\leq\frac{2Ka_{\min}\nu_{\max}}{\ve^2}
\exp\left(-\frac{3}{28}\alpha_{\min}n\ve^2\right).
\end{align*}
}
This ends the proof since on the complement of this event, \eqref{eq:ineq-mass-covering} asserts that
\begin{align*} 
\min_{1 \leq k \leq K}
\min_{B \in \cB_{k}}
P_{k,n}(B) 
&\geq 
\frac12
\min_{1 \leq k \leq K}
\min_{B \in \cB_{k}}
P_k(B) \geq  \upsilon \ve^2
.
\end{align*} 
\end{proof}
\section{Algorithmic guarantees} \label{app:algo}

To prove the performances of \codetection, the following notation for thickenings of sets will be useful.
\begin{definition}[Offset]
\label{def:offset}
Given $A \subset \bbR^D$ and $r\geq 0$, the \emph{$r$-offset} of $A$ is
\begin{align*}
A^r
:=
\{x \in \bbR^D \mid \dist(x \mid A) \leq r
\}
.
\end{align*}
\end{definition}

We are now in position to prove the main technical result on \codetection.
Write:
\begin{itemize}
\item
$\wt\cX^{(d)}$ for the current value of set $\wt\cX$ at the \emph{start} of step $d \in \{1,\dots,d_{\max}\}$, 
\item
$\wh \cX^{(d)}$ for the set of labeled points at the end of step $d$.
\end{itemize}
In particular, $\wt\cX^{(1)} = \cX$ is the input point cloud and $\wt\cX^{(d+1)} = \wt\cX^{(d)}\setminus \wh\cX^{(d)}$.
\begin{lemma} 
\label{lem:recalgo} 
There exist constants $\Upsilon^*$, $\sigma^*$ depending on $D$, and $\zeta^*_d$ depending on $d$ such that, if we choose
\begin{align}  \label{eq:optiparam1}
\Upsilon \geq \nu_{\max} \Upsilon^*{\color{black}/\alpha_{\min}},\quad \gamma \geq \frac{800 D^2\log D}{\alpha_{\min}} \vee \frac{56 q}{3 \alpha_{\min}}, \quad 4D \leq \sigma \leq \sigma^*,\quad \zeta_d \geq \zeta_d^*,
\end{align} 
and
\begin{align}  \label{eq:optiparam2}
\begin{aligned}
\ve &:= \{\gamma \frac{\log n}{n}\}^{1/2},\quad\rho_d :=  \{\frac{\Upsilon \ve^2}{a_{\min} \kappa_{\max}^d }\}^{1/d},\quad h_{d} := 48\(1+{\color{black}\frac1d}\) d \rho_d, \quad
\kappa_d :=  \kappa_{\max}, \\
r_d &:= h_d, 
\quad
\delta_d := \zeta_d \kappa_{\max} h_{d}^2,
\quad
n_d := \sigma n \ve^2
,
\end{aligned}
\end{align} 
the following holds.

For $n$ large enough{\color{black}, with the tuning constants fixed independently of $n$}, when run with an i.i.d. point cloud $\cX_n = \{X_1,\ldots,X_n\}$ from some unknown distribution $P = \sum_{k=1}^K \alpha_k P_k \in \bar{\cP}$ and parameters~\eqref{eq:optiparam1} and~\eqref{eq:optiparam2},
the algorithm \emph{\codetection} operates in a way such that with probability at least $1-6Kn^{-q}$, at all steps $d \in \{1,\ldots,d_{\max}\}$, we have:
\begin{enumerate}
\item If a $(d+1)$-tuple of points $\bx = (x_1,\dots,x_{d+1})$ in $\wt \cX^{(d)}$ is co-detected, then $d = d_k$ for some $k \in \{1,\ldots,K\}$ and
\beq \label{def:xi_d}
\dH(\conv(x_{1},\dots,x_{d_k+1}) \mid M_k) \leq \xi_{d_k}\kappa_{\max} h^2_{d_k}\quad\text{where}\quad \xi_{d_k} = C_{d_k} \Upsilon \nu_{\max}\frac{a_{\max}}{a_{\min}}.
\eeq
Furthermore, for all $y \in \conv(x_{1},\dots,x_{d_k+1})$, there exists $x \in \ball(y, \xi_{d_k} \kappa_{\max} h^2_{d_k}) \cap M_k$ such that 
$$
\angle({\color{black}T(x_1,\dots,x_{d_k+1})} \mid T_x M_k) \leq \xi_{d_k} \kappa_{\max} h_{d_k}.
$$
\item If $d = d_k$, then for all $x \in M_k$ and all $d_k$-dimensional subspace $T \subset T_x M_k$, there exists a $(d_k+1)$-tuple $(x_1,\dots,x_{d_k+1})$ of $\wt\cX^{(d_k)}$ co-detected such that,
\begin{align*}
\dist(x \mid \conv(x_{1},\dots,x_{d_k+1})) \leq \zeta^*_{d_k} \kappa_{\max} h_{d_k}^2\\
\text{and}~~~~\angle(T, {\color{black}T(x_1,\dots,x_{d_k+1})}) \leq C_{d_k}  \kappa_{\max} h_{d_k}.
\end{align*} 
\end{enumerate}
\end{lemma}
For the proof of \lemref{recalgo}, we set ourselves on the intersection of the events $\cE_1$ and $\cE_2$ of \lemref{slab_concentration} and \lemref{eventcE}, which hold simultaneously with these parameters with probability higher than $1-6Kn^{-q}$, as soon as we take 
$$\Upsilon^* \geq \upsilon^* \max_{1 \leq d \leq D-1} \frac{2^{d+1}}{\omega_d}$$
for some $\upsilon^* \geq 2$. 
We also take $\Upsilon^* \geq 1$.
The value of $\upsilon^*$ will be specified later in the proof, see~\eqref{eq:defupsilon}.
The \lemref{recalgo} is proven by induction on $d \in \{1,\ldots,d_{\max}\}$. The initialization step is essentially the same as the induction step but easier, since we do not have to handle the points that have been removed from $\wt\cX^{(d)}$ in previous steps.
All steps are thus treated identically. 
We let $d \in \{1,\ldots,d_{\max}\}$ and we assume that points 1 and 2 of \lemref{recalgo} holds for all $d' \in \{1,\ldots,d-1\}$ (perhaps empty).
In particular, point $1$ and $2$ respectively yield the two inclusions
\beq 
\label{eq:excl}
\begin{cases}
\wt \cX^{(d)}
\supset 
\cX_n \setminus \bigcup_{d_k < d} M_k^{ \xi_{d_k} \kappa_{\max} h^2_{d_k}+\delta_{d_k}} 
\\
\wt \cX^{(d)} 
\subset 
\cX_n \setminus \bigcup_{d_k < d} M_k,
\end{cases}
\eeq
which will be the only thing we need for the induction to work --- whence the uselessness of treating the initialization step separately since~\eqref{eq:excl} holds automatically true for $d=1$. 

\begin{rem}[On the chosen immersions below]
In what follows, the $d_k$-dimensional support $M_k = \supp(P_k)$ of the unknown distribution $P = \sum_{k=1}^K \alpha_k P_k \in \bar{\cP}$ is parametrized by an arbitrary immersion $(\Psi_k,\bbM_k) \in \bbI(M_k)$ such that 
$$\vol(\bbM_k) \leq 2 \nu_{\max} (\kappa_{\max})^{-d_k},
$$
chosen once and for all.
\end{rem}
\begin{proof}[Proof of \lemref{recalgo} --- point 1, induction step.]
Let $\bx = (x_{1},\dots,x_{{d+1}})$ be a $(d+1)$-tuple of $\wt \cX^{(d)}$ co-detected by the algorithm, 
meaning that $\rad(\bx) \leq r_{d}$ and $\Card\(S \cap \wt \cX^{(d)}\) \geq \sigma n \ve^2$, where $S := \slab(\bx,h_{d},\kappa_{\max} h_{d}^2)$ is the slab associated to $\bx$.

\begin{itemize}
\item
\textbf{Proof that $d \in \{d_1,\ldots,d_K\}$.}
Let us write 
$
k := \inf\{k' \mid d_{k'} \geq d\}
$.
{\color{black}If this set is empty, then $\wt\cX^{(d)}=\emptyset$ by~\eqref{eq:excl}, so no tuple can be co-detected. Otherwise, all remaining points come from components of dimension at least $d$. Almost surely, every $(d+1)$-tuple of such points is affinely independent: conditionally on any preceding points, their affine span has dimension at most $d-1$ and therefore has zero mass under each such component.}
From the second inclusion of the inductive hypothesis~\eqref{eq:excl}, $\wt \cX^{(d)} \cap \cX_{k',n} \subset \wt \cX^{(d)} \cap M_{k'} = \emptyset$ for all $k'<k$. Hence, we have
$$
\sum_{k' \geq k} P_{k',n}(S) 
\geq 
\frac{1}{n}
\sum_{k' \geq k} 
\Card\(S \cap \wt \cX^{(d)} \cap \cX_{k',n}\)
=
\frac1n \Card\(S \cap \wt \cX^{(d)}\) 
\geq 
\frac{n_{d}}n = \sigma \ve^2.
$$
As a result, there exists $k_0 \geq k$ such that $P_{k_0,n}(S) \geq \sigma \ve^2/K$.
Using {\color{black}Corollary~\ref{cor:conc}}, and the fact that $\sigma \geq 4D \geq 4K$, we find that $P_{k_0}(S) \geq \ve^2$.
{\color{black}Suppose by contradiction that $d_{k_0}>d$. }On the other hand, letting $C_{{\color{black}d_{k_0}}}$ denote a generic constant depending on ${\color{black}d_{k_0}}$ only, point (ii) of \lemref{slab_vol2} yields that for $n$ large enough 
\begin{align*}
\ve^2
\leq
P_{k_0}(S) 
&\leq 
C_{d_{k_0}} (a_{\max} \kappa_{\max}^{d_{k_0}}) (\nu_{\max}/\kappa_{\max}^{d_{k_0}}) \kappa_{\max}^{2d_{k_0}-d} h_{d}^{2d_{k_0}-d}
\\
&=
C_{d_{k_0}} a_{\max} \nu_{\max} (\kappa_{\max} h_{d})^{2d_{k_0}-d}
,
\end{align*}
and since $h_d = {\color{black}A_d}(\Upsilon \ve^2/(a_{\min} \kappa_{\max}^d))^{1/d}${\color{black} with $A_d := 48(1+1/d)d$}, it yields
\begin{align*}
{\color{black}a_{\min}/(\Upsilon A_d^d)}
\leq
C_{d_{k_0}} a_{\max} \nu_{\max} 
\(\kappa_{\max} h_d \)^{{\color{black}2(d_{k_0}-d)}}
.
\end{align*}
{\color{black}The right-hand side tends to zero as $n\to\infty$, a contradiction.} Therefore $d_{k_0} = d$.
As $d \leq d_k \leq d_{k_0}$ by construction, we hence have proven that $d = d_k$.

\item
\textbf{Proof of the Hausdorff distance bound.} As $P_k(S) = P_{k_0}(S) \geq \ve^2$, we get
\begin{align*}
\cH^{(d_k)}\(M_k \cap S\) 
&\geq 
\frac{\ve^2}{a_{\max}\kappa_{\max}^{d_k}} 
= 
\frac{c_{d_k}}{\Upsilon} 
\frac{a_{\min}}{ a_{\max}} h_{d_k}^{d_k}
.
\end{align*}
Hence, point (iii) of \lemref{slab_vol2} yields the existence of some $x \in M_k$ such that 
$\|\bar{\bx} - x\| \leq \mu_k \kappa_{\max} h^2_{d_k}$
and $\angle(T \mid T_{x}{\color{black}M_k}) \leq \mu_k \kappa_{\max} h_{d_k}$, where $T := {\color{black}T(\bx)}$ {\color{black}is the linear direction of the affine span} of $\bx = (x_1,\ldots,x_{d+1})$, and
$$
\mu_k := C_{d_k}(1+\nu_{\max}\beta^{-1})\vee 1 ~~~\text{with}~~~\beta := 
\frac{c_{d_k}}{\Upsilon}
\frac{a_{\min}}{ a_{\max}} 
.
$$
In particular, we get
\beq \label{eq:boundmuk}
1 \leq \mu_k \leq C'_{d_k}\nu_{\max} \Upsilon \frac{a_{\max}}{a_{\min}},
\eeq
where we used that $\Upsilon$ and $\nu_{\max}$ were lower-bounded by numeric constants.
Because $\rad(\bx) \leq r_{d_k}$, we have in particular that $\|x_{j} - x\| \leq 2r_{d_k}+\mu_k \kappa_{\max} h^2_{d_k}$ for all $j \in \{1,\ldots,d_k+1\}$.
Recall that $T_x M_k = \bigcup_{x_0 \in \Psi_k^{-1}(x)} T_{x_0}^{(\Psi_k)} M_k$. 
Write $x_0 \in \Psi_k^{-1}(x)$ for any pre-image of $x$ such that $\angle(T,T_{x_0}^{(\Psi_k)} M_k) \leq \mu_k \kappa_{\max} h_{d_k}$.
{\color{black}Write $E:=T_{x_0}^{(\Psi_k)}M_k$. For every $z\in\conv(x_1,\ldots,x_{d_k+1})$, the enclosing-ball bound gives
\[
\|z-\bar{\bx}\|\leq2h_{d_k},\qquad
\|z-x\|\leq2h_{d_k}+\mu_k\kappa_{\max}h_{d_k}^2.
\]
Since $z-\bar{\bx}\in T$, we have
\begin{align*}
\|(I-\pr_E)(z-x)\|
&\leq\|\pr_T-\pr_E\|_{\op}\|z-\bar{\bx}\|+\|\bar{\bx}-x\|\\
&\leq3\mu_k\kappa_{\max}h_{d_k}^2.
\end{align*}
Let $v:=\pr_E(z-x)$ and $w:=\overline{\exp}_{x_0}^{\bbM_k}(v)\in M_k$.
For $n$ large enough, $\|v\|\leq1/(8\kappa_{\max})$, and Lemma~\ref{lem:geoineq} 2. yields
\begin{align*}
\|w-z\|
&\leq\kappa_{\max}\|v\|^2+\|(I-\pr_E)(z-x)\|\\
&\leq\kappa_{\max}(2h_{d_k}+\mu_k\kappa_{\max}h_{d_k}^2)^2
+3\mu_k\kappa_{\max}h_{d_k}^2\\
&\leq C\mu_k\kappa_{\max}h_{d_k}^2.
\end{align*}
Taking the supremum over $z$ proves
\[
\dH(\conv(x_1,\ldots,x_{d_k+1})\mid M_k)\leq C\mu_k\kappa_{\max}h_{d_k}^2,
\]
}
which is what we wanted to show. 
\item
\textbf{Proof of the angle bound.}
Finally, for any $z \in \conv(x_{1},\dots,x_{d_k+1})$, we can let $w = \overline{\exp}_{x_0}^{\bbM_k}({\color{black}\pr_E(z-x)}) \in M_k$. The same computations as above show that 
$\|z -w\| \leq C\mu_k\kappa_{\max} h_{d_k}^2$ 
and that 
\begin{align*} 
\angle(T \mid T_w M_k) &\leq \angle(T,T_{x_0}^{(\Psi_k)}M_k) + \angle(T_{x_0}^{(\Psi_k)}M_k \mid T_{w} M_k) \\
&\leq \mu_k\kappa_{\max}h_{d_k} + 5 \kappa_{\max} \|w-x\| \\
&\leq  \mu_k\kappa_{\max}h_{d_k} + 10 \kappa_{\max}r_{d_k} + 5\mu_k\kappa^2_{\max}h_{d_k}^2 \\
&\leq 
C \mu_k \kappa_{\max} h_{d_k},
\end{align*} 
which concludes the proof.
\end{itemize}
\end{proof}

\begin{proof}[Proof of \lemref{recalgo} --- point 2, induction step.] We let $x_0 \in (\Psi_k)^{-1}(x)$ be fixed. We identify $T_{x_0}^{(\Psi_k)} {\color{black}M_k}$ to $\bbR^{d_k}$ with its canonical basis $(e_1,\dots,e_{d_k})$. Define $e_{{d_k}+1}$ to be the unit vector
\begin{align*}
e_{{d_k}+1} := -\frac1{\sqrt{d_k}}\sum_{i=1}^{d_k} e_i
.
\end{align*}
Then, for all $i \in \{1,\dots,{d_k}+1\}$ and $s > 0$, let us write
$$\cQ_i(s) := \{v \in \bbR^{d_k}~\middle|~\|v-e_i\| \leq s\}$$
and
$$\wt\cQ_i(s) := \wt h_{d_k} \cQ_i(s)~~\text{with}~~ \wt h_{d_k} := \frac{1}{2(1+1/(8d_k))} h_{d_k}.
$$
Given any family of vectors $v_i \in \wt\cQ_i(1/8{d_k})$ for $1 \leq i \leq {d_k}+1$ with associated renormalized $\bar v_i := v_i/\|v_i\|$, we will show that
\begin{itemize}
\item[1.] For all $i \in\{1,\ldots,{d_k}+1\}$, $\bar v_i \in \cQ_i(1/4{d_k})$;
\item[2.] For all $i \neq j \in\{1,\ldots,{\color{black}d_k}\}$, there holds that $|\inner{\bar v_i}{\bar v_j}| \leq 1/2{d_k}$;
\item[3.] There holds that $0 \in \conv\(v_1,\dots,v_{{d_k}+1}\)$;
\item[4.] For all $i \in \{1,\dots,{d_k}+1\}$,
$$
\cX_{k,n} 
\cap 
\Bigl(
\overline{\exp}_{x_0}^{\bbM_k}(\wt\cQ_i(1/8{d_k})) 
\setminus 
\bigcup_{k' < k} M_{k'}^{\xi_{d_{k'}}\kappa_{\max} h_{d_{k'}}^2 + \delta_{d_{k'}}}\Bigr) 
\neq 
\emptyset
.
$$
\end{itemize}
For point 1, simply notice that for all $\lambda  > 0$, $\bar v_i$ is the projection of $\lambda v_i$ onto the unit sphere. Triangle inequality and this fact applied to $\lambda = 1/\wt h_{d_k}$ hence yield
\begin{align*}
\|\bar v_i - e_i\| 
&\leq 
\left\|\bar v_i - v_i/\wt h_{d_k}\right\|
+
\left\|v_i/\wt h_{d_k}-e_i\right\| 
\leq 
2\left\|v_i/\wt h_{d_k}-e_i\right\|
\\
&\leq 
\frac1{4d_k}
,
\end{align*} 
where the last inequality comes from the fact that $\|v_i - \wt h_{d_k} e_i \| \leq \wt h_{d_k}/(8d_k)$ by construction. 

For point 2, use point 1 to get
\begin{align*}
|\inner{\bar v_i}{\bar v_j}| 
&\leq 
|\inner{\bar v_i-e_i}{\bar v_j}|+
|\inner{e_i}{\bar v_j-e_j}|
\\
&\leq
1/2d_k
.
\end{align*}

For point 3, we let $V := \(\inner{\bar v_j}{e_i}\)_{1 \leq i,j \leq d_k} \in \bbR^{d_k \times d_k}$ and $A := I-V${\color{black}, where $\|\cdot\|_\infty$ denotes the maximum absolute entry}. 
We have $\|A\|_\infty \leq 1/4d_k$ so that $\|A^\ell \|_\infty \leq d_k^{\ell-1}/(4d_k)^\ell = 1/(d_k 4^\ell)$ for all $\ell \geq 1$. 
We deduce that $V$ is invertible, with the decomposition
$$
V^{-1} = I + R~~~~\text{and}~~~\|R\|_{\infty} \leq \sum_{\ell = 1}^\infty \frac1{d_k 4^\ell} \leq \frac{1}{2d_k}.
$$
From the invertibility of $V$, one can write
\begin{align*}
-\bar v_{d_k+1} = \sum_{i=1}^{d_k} \inner{-\bar v_{d_k+1}}{e_i} e_i = \sum_{j=1}^{d_k} a_j \bar v_j~~~\text{with}~~~a_j = \sum_{i=1}^{d_k} (V^{-1})_{{\color{black}j,i}} \inner{-\bar v_{d_k+1}}{e_i}. 
\end{align*}
Furthermore, point 1 yields that
$$
\left| {\color{black}\inner{-\bar v_{d_k+1}}{e_i}} - \frac1{\sqrt{d_k}}\right| \leq 1/(4d_k)
$$ 
for all $i \in \{1,\ldots,d_k\}$.
As a result,
\begin{align*} 
a_j 
&= 
(V^{-1})_{j,j} \inner{-\bar v_{d_k+1}}{e_j} + \sum_{i\neq j} (V^{-1})_{{\color{black}j,i}} \inner{-\bar v_{d_k+1}}{e_i} 
\\
&\geq
\(1 - \frac1{2d_k}\)\(\frac{1}{\sqrt{d_k}} - \frac{1}{4d_k}\) - (d_k-1)\times \frac{1}{2d_k} \times\({\color{black}\frac1{\sqrt{d_k}}} + \frac{1}{4d_k}\)
\\
&\geq 
\frac1{2\sqrt{d_k}} - \frac3{8d_k} > 0
.
\end{align*} 
From which we conclude that the convex coefficients $\bar{a}_j := {\color{black}(a_j/\|v_j\|)/(\sum_{m=1}^{d_k+1} a_m/\|v_m\|)}$ for $j \in \{1,\ldots,d_k+1\}$ with $a_{d_k+1} := 1$  satisfy
$$
0 = \sum_{j=1}^{{\color{black}d_k+1}} \bar{a}_j {\color{black}v_j} \in \conv\{ v_1,\dots, v_{d_k+1}\}
.
$$

For point 4, write
$$
\lambda_{k'} := \xi_{d_{k'}} \kappa_{\max} h_{d_{k'}}^2+\delta_{d_{k'}} =: \tau_{k'}\kappa_{\max}h_{d_{k'}}^2
$$
for all $k'<k$.
By construction, $\tau_{k'}$ is greater than some numeric constant, i.e. $\tau_{k'} \geq c$.
Let also $B_i := \overline{\exp}_{x_0}^{\bbM_k}(\wt\cQ_i(1/8{d_k}))$ for $i \in \{1,\ldots,{\color{black}d_k+1}\}$.
Then we have
\begin{align}
\label{eq:empirical-lower-bound-induction} 
P_{k,n}\bigl(B_i \setminus \bigcup_{k' < k} M_{k'}^{\lambda_{k'}}\bigr) \geq P_{k,n}\(B_i\) - \sum_{k' < k} P_{k,n}\bigl(B_i \cap M_{k'}^{\lambda_{k'}}\bigr).
\end{align} 
Let now $k' < k$ be fixed. 
Notice that $B_i \cap M_{k'}^{\lambda_{k'}} \subset \bigl(B_i^{\lambda_{k'}} \cap M_{k'}\bigr)^{\lambda_{k'}}$. 
For $n$ large enough, we have $\lambda_{k'} {\color{black}\leq \wt h_{d_k}/4}${\color{black} and $h_{d_{k'}} \leq \wt h_{d_k}$. Furthermore, \lemref{geoineq} 2. yields $\rad(B_i) \leq \wt h_{d_k}/(8d_k)+\kappa_{\max}h_{d_k}^2/2 \leq \wt h_{d_k}/(4d_k)$ for $n$ large enough}, and hence
\begin{align*}
\rad\bigl(B_i^{\lambda_{k'}}\bigr) 
&\leq 
\frac{\wt h_{d_k}}{4d_k} + \lambda_{k'} 
\\
&\leq 
\frac12 \wt h_{d_k} 
< 
\frac1{4\kappa_{\max}}
.
\end{align*} 
As a result, \lemref{sheaf_number} yields the existence of $z_1,\dots,z_{N'}$ in $\bbM_{k'}$ with $N' \leq N_0^{(d_{k'})}$ such that
$$
B_i^{\lambda_{k'}} \cap M_{k'} \subset \bigcup_{\ell=1}^{N'} {\color{black}\Psi_{k'}}(\ball_{\bbM_{k'}}(z_\ell,2\wt h_{d_k})).
$$
With points 1 and 5 of \lemref{geoineq}, we further find that for all $\ell \in \{1,\ldots,N'\}$,
$$
{\color{black}\Psi_{k'}}(\ball_{\bbM_{k'}}(z_\ell,{\color{black}2}\wt h_{d_k})) \subset \slab_{T_{z_\ell}^{(\Psi_{k'})} M_{k'}}(x_\ell, 2\wt h_{d_k}, 2\kappa_{\max}\wt h_{d_k}^2),
$$
where $x_\ell := \Psi_{k'}(z_\ell)$, so that in the end
\begin{align}
\label{eq:offset-in-slab-different-dims}
B_i \cap M_{k'}^{\lambda_{k'}}  
\subset 
\bigcup_{\ell=1}^{N'} \slab_{T_{z_\ell}^{(\Psi_{k'})} M_{k'}}(x_\ell, 3\wt h_{d_k}, (2+\tau_{k'})\kappa_{\max}\wt h_{d_k}^2)
. 
\end{align}
Using point (ii) of \lemref{slab_vol2}, we find that
\begin{align*} 
P_k&
\Bigl(
\slab_{T_{z_\ell}^{(\Psi_{k'})} M_{k'}}(x_\ell, 3\wt h_{d_k}, (2+\tau_{k'})\kappa_{\max}\wt h_{d_k}^2)\Bigr) 
\\
&\leq  (a_{\max} \kappa_{\max}^{d_k}) \times (C_{d_k} \tau_{k'}^{d_k-d_{k'}} \nu_{\max} \kappa_{\max}^{d_k-d_{k'}} \wt h_{d_k}^{2d_k-d_{k'}}).
\end{align*} 
Because $d_k > d_{k'}$, we can put the latter bound in the form $c_n \ve^2$ where $c_n = o(1)$ so that using~{\color{black}Corollary~\ref{cor:conc}}, we get that 
$$
P_{k,n}
\Bigl(\slab_{T_{z_\ell}^{(\Psi_{k'})} M_{k'}}({\color{black}x_\ell}, 3\wt h_{d_k}, (2+\tau_{k'})\kappa_{\max}\wt h_{d_k}^2)\Bigr) \leq {\color{black}(1+\sqrt{c_n})^2}\ve^2  \leq 2\ve^2,$$
for $n$ large enough. This latter bound for all $k'<k$ and all $\ell \leq N'$ together with~\eqref{eq:offset-in-slab-different-dims} yield the final bound
\begin{align*}
\sum_{k' < k} P_{k,n}\bigl(B_i \cap M_{k'}^{\lambda_{k'}}\bigr) \leq C_D \nu_{\max} \ve^2
,
\end{align*}
where we used that $K \leq D$ and $N_0^{(d_{k'})} \leq C_D \nu_{\max}$ (see Lemma~\ref{lem:sheaf_number}).
{\color{black}On the other hand, writing $u_i := \wt h_{d_k}e_i$ and $s_k := \wt h_{d_k}/(8d_k)$, our choice of $h_{d_k}$ gives
$$
s_k = \left(3+\frac{21}{8d_k+1}\right)\rho_{d_k}>3\rho_{d_k}.
$$
For $\|v-u_i\|=s_k$, \lemref{geoineq} 2. yields
$$
\left\|\overline{\exp}_{x_0}^{\bbM_k}(v)-\overline{\exp}_{x_0}^{\bbM_k}(u_i)\right\|
\geq s_k-\kappa_{\max}(\|v\|^2+\|u_i\|^2)
\geq s_k-\kappa_{\max}h_{d_k}^2/2>3\rho_{d_k}
$$
for $n$ large enough. As any geodesic path in $\bbM_k$ leaving the corresponding exponential-chart patch must cross its boundary, this patch contains the intrinsic ball of radius $3\rho_{d_k}$ centered at $\exp_{x_0}^{\bbM_k}(\diff\Psi_k(x_0)^{-1}[u_i])$. Since $\cP_k$ is a $2\rho_{d_k}$-net, $B_i$ therefore contains a ball of the form $\Psi_k\(\ball_{\bbM_k}(p,\rho_{d_k})\)$ with center $p \in \cP_k$ from \lemref{eventcE}.}
As $\Upsilon \geq \nu_{\max} \Upsilon^*{\color{black}/\alpha_{\min} \geq \nu_{\max} \Upsilon^*}$, {\color{black}we have $P_{k,n}(B_i)\geq\nu_{\max}\upsilon^*\ve^2$}, and thus~\eqref{eq:empirical-lower-bound-induction} yields
\begin{align} \label{eq:defupsilon}
P_{k,n}\Bigl(B_i \setminus \bigcup_{k' < k} M_{k'}^{\lambda_{k'}}\Bigr) \geq (\upsilon^*-C_D) \nu_{\max} \ve^2,
\end{align}
and it is sufficient to take $\upsilon^* > C_D$ to make sure that $P_{k,n}\bigl(B_i \setminus \bigcup_{k' < k} M_{k'}^{\lambda_{k'}}\bigr) > 1/n$ for $n$ large enough, which ends the proof of point 4.

\begin{itemize}
\item
\textbf{Proof of the distance bound.}
For all $1 \leq i \leq {\color{black}d_k+1}$ we let $x_i \in \cX_{k,n}$ be a point that belongs to ${\color{black}\overline{\exp}}_{x_0}^{\bbM_k}(\wt\cQ_i(1/8d_k)) \setminus \bigcup_{k' < k} M_{k'}^{\lambda_{k'}}$. We denote $v_i := \{\overline{\exp}_{x_0}^{\bbM_k}\}^{-1}(x_i)$ and let  $\lambda_i \in [0,1]$ summing to $1$ such that $0 = \sum_{i=1}^{{\color{black}d_k+1}} \lambda_i v_i$ (point 3). 
From \lemref{geoineq}{\color{black}, as $\|v_i\|\leq (1+1/(8d_k))\wt h_{d_k}=h_{d_k}/2$}, we have $\|x_i - x - v_i\| \leq {\color{black}\kappa_{\max}h_{d_k}^2/4}$, from which we deduce
\begin{align*}
\left\|x - \sum_{i=1}^{{\color{black}d_k+1}} \lambda_i x_i\right\| 
&= 
\left\|\sum_{i=1}^{{\color{black}d_k+1}} \lambda_i (x + v_i - x_i)\right\| 
\\
&\leq 
\sum_{i=1}^{{\color{black}d_k+1}} \lambda_i {\color{black}\frac{\kappa_{\max}}4 h_{d_k}^2} 
\\
&= 
{\color{black}\frac{\kappa_{\max}}4 h_{d_k}^2},
\end{align*}
so that 
\begin{align*}
\dist(x \mid \conv(x_1,\dots,x_{d_k+1})) 
&\leq 
{\color{black}\kappa_{\max} h_{d_k}^2/4} 
\\
&=: 
\zeta_{d_k}^* \kappa_{\max} h_{d_k}^2
,
\end{align*}
where $\zeta_{d_k}^* {\color{black}= 1/4}$.

\item
\textbf{Proof of the angle bound.}
We now let {\color{black}$A := \operatorname{aff}\{x_1,\dots,x_{d_k+1}\}$, $T := A-A$ and $p := \pr_A(x)$}.
Write {\color{black}$w_i := x_i-p$ for $1\leq i\leq d_k$}. We have 
\begin{align*}
\|w_i - v_i\| 
&= 
{\color{black}\|x_i-p-v_i\|} 
\\
&\leq 
{\color{black}\|x-p\|}+\|x_i-x-v_i\| 
\\
&\leq 
{\color{black}\frac{\kappa_{\max}}2 h_{d_k}^2}
,
\end{align*}
so that{\color{black}, writing $\bar w_i:=w_i/\|w_i\|$ and using $\|v_i\|\geq (1-1/(8d_k))\wt h_{d_k}$, for $n$ large enough,}
\begin{align*}
\| \bar w_i - \bar v_i \| 
&\leq 
\frac{\|w_i - v_i\|}{\|v_i\|} + \|w_i\|\left | \frac{1}{\|w_i\|}-\frac{1}{\|v_i\|}\right| \\
&\leq 
2\frac{\|w_i-v_i\|}{\|v_i\|} 
\\
&\leq 
{\color{black}6}\kappa_{\max}\wt h_{d_k}
.
\end{align*}
We complete $(e_i)_{1\leq i \leq d_k}$ to a full orthonormal basis of $\bbR^D$ and extend $V \in \bbR^{D \times d_k}$ (see point 3) with $0$. Let also $W := (\inner{\bar w_j}{e_i}) \in \bbR^{D \times d_k}$.
From the last bound, we have $\|W - V\|_{\op} \leq {\color{black}6}\sqrt{d_k}\kappa_{\max} \wt h_{d_k}$, and thanks to point 2 and the Gershgorin circle theorem,
$$
\|(V^\top V)^{-1}\|_{\op} \leq \(1- \frac{d_k-1}{2d_k}\)^{-1} \leq 2
.
$$ 
{\color{black}For $n$ large enough, $W$ has rank $d_k$, so that $T = \Vect\{w_1,\dots,w_{d_k}\}$. }We thus have that 
\beq\label{eq:cstar}
\| W(W^\top W)^{-1} W^{\top}- V(V^{\top} V)^{-1}V^{\top} \|_{\op} \leq C \sqrt{d_k} \kappa_{\max} \wt h_{d_k}
\eeq
for some explicit numerical constant $C$, whence $\angle(T,T_{x_0}^{(\Psi_k)} M_k) \leq  C \sqrt{d_k} \kappa_{\max} \wt h_{d_k}$. 

\textbf{
\item
Proof that $x_1,\dots,x_{{\color{black}d_k+1}}$ are codetected.
}
First, we know that
\begin{align*}
\|x - x_i \|\leq 
\|v_i\| 
\leq 
\wt h_{d_k} + \frac{1}{8 d_k} \wt h_{d_k}  
\leq 
h_{d_k}/2
\leq
r_{d_k}/2,
\end{align*}
so that $\rad\(x_1,\dots,x_{d_k+1}\) \leq r_{d_k}$.
{\color{black}Second, writing $\bar{\bx}:=(d_k+1)^{-1}\sum_{i=1}^{d_k+1}x_i$, the preceding bounds yield
$$
\|\bar{\bx}-x\|\leq h_{d_k}/2,\qquad
\|\pr_{T^\perp}(\bar{\bx}-x)\|=\dist(x\mid A)\leq\kappa_{\max}h_{d_k}^2/4.
$$
}{\color{black}Hence} the slab defined by $(x_1,\dots,x_{{\color{black}d_k+1}})$ contains the slab
$$S := \slab_T\(x,h_{d_k}/2,\kappa_{\max} h_{d_k}^2/2\).$$
Using \lemref{slab_vol3} and the angle bound we just showed, we get that 
$$
P_k(S) \geq C_{d_k} a_{\min} \kappa_{\max}^{d_k} h_{d_k}^{d_k} 
\geq 
C_{D} \nu_{\max} \upsilon^* \ve^2{\color{black}/\alpha_{\min}}.
$$
Assuming that we took $\upsilon^*$ large enough so that 
$
4 \leq C_{D}\nu_{\max}\upsilon^*{\color{black}/\alpha_{\min}},
$
which we can (recall that $\nu_{\max}$ is lowered-bounded by a numeric constant without loss of generality, see Remark~\ref{rem:model-parameters}), Corollary~\ref{cor:conc} yields{\color{black}, after adjusting $C_D$,}
\begin{align*} 
P_{k,n}(S) &\geq \frac12P_k(S) \geq C_{D} \nu_{\max} \upsilon^* \ve^2{\color{black}/\alpha_{\min}}.
\end{align*} 
On the other hand, in the same spirit as what as been done above, one can show that
\begin{align*} 
P_{k,n}\bigl(S \cap \bigcup_{k'<k} M_{k'}^{\lambda_{k'}}\bigr) \leq C'_D \nu_{\max} \ve^2. 
\end{align*} 
{\color{black}On $\cE_2$, we also have $N_k>\alpha_kn/2\geq\alpha_{\min}n/2$. }In the end, we thus have
\begin{align*} 
\Card\(S \cap \wt\cX^{(d_k)}\) 
&\geq 
\Card\Bigl( {\color{black}\cX_{k,n}} \cap \bigl(S \setminus \bigcup_{k'<k} M_{k'}^{\lambda_{k'}}\bigr)\Bigr) 
\\
&{\color{black}= N_k P_{k,n}\bigl(S \setminus \bigcup_{k'<k} M_{k'}^{\lambda_{k'}}\bigr)}
\\
&{\color{black}\geq N_k(C_{D}\upsilon^*/\alpha_{\min} - C'_D) \nu_{\max} \ve^2}
\\
&{\color{black}\geq \frac12(C_{D}\upsilon^* - \alpha_{\min}C'_D) \nu_{\max} n\ve^2},
\end{align*} 
Using {\color{black}$\alpha_{\min}\leq1$ and }again the fact that $\nu_{\max}$ is lower-bounded by a numeric constant (see Remark~\ref{rem:model-parameters}), we can write{\color{black}, after adjusting the constants,} 
$$
\Card\(S \cap \wt\cX^{(d_k)}\)  \geq {\color{black}\frac12}(C_{D}\upsilon^* - C'_D)  n\ve^2 \geq {\color{black}\frac14} C_D \upsilon^* n \ve^2,
$$
as soon as $\upsilon^* \geq 2C_D'/C_D$. 
We conclude by setting $\sigma^* := {\color{black}\frac14} C_D \upsilon^*$ and by taking $\upsilon^*$ large enough so that $\sigma^* \geq 4D$. 
The tuple $(x_1,\dots,x_{d_k+1})$ is then co-detected by the algorithm at step $d_k$, which ends the proof.

\end{itemize}
\end{proof}

\section{Proof of the main results} 
\label{sec:proof-main-results}

The proofs of \thmref{dimension-and-component-estimation}, \thmref{main-dimension-labeling}, \thmref{main-manifold-estimation} and \thmref{main-tangent-estimation} are straightforward applications of \lemref{recalgo}. 
{\color{black}For Theorem~\ref{thm:main-tangent-estimation}, apply each part of Lemma~\ref{lem:recalgo} separately to each tangent branch, and take $\Delta_k$ large enough to contain the corresponding distance bounds. This is precisely the order of the supremum and infimum in the localized losses.}
It only remains to prove \thmref{main-clustering-dimensions}. The proof relies on the following lemma.

\begin{lemma} \label{lem:offset} Let $k' < k$. Then for all ${\color{black}0\leq\lambda < 1/(256\kappa_{\max})}$
\begin{align*}
P_k(M_{k'}^{\lambda}) \leq C_{d_k}a_{\max}\nu_{\max}^2 \kappa_{\max}^{d_k-d_{k'}} \lambda^{d_k-d_{k'}},
\end{align*}
\end{lemma}

\begin{proof} 
{\color{black}The case $\lambda=0$ follows from $\cH^{d_k}(M_{k'})=0$. For $\lambda>0$, }{\color{black}we} first write that 
\begin{align*}
P_k(M_{k'}^{\lambda}) \leq (a_{\max}\kappa_{\max}^{d_k})
\cH^{(d_k)}(M_k \cap M_{k'}^{\lambda})
,
\end{align*}
hence reducing the result to a volume bound.
From the volume bound $\vol(\bbM_{k'}) \leq 2\nu_{\max} \kappa_{\max}^{-d_{k'}}$, we can cover $\bbM_{k'}$ with at most 
$$N' \leq C_{d_{k'}}\nu_{\max} \kappa_{\max}^{-d_{k'}} / \(\sqrt{\lambda/\kappa_{\max}}\)^{d_{k'}} = C_{d_{k'}}\nu_{\max} (\kappa_{\max} \lambda)^{-d_{k'}/2}$$
balls $\ball_{\bbM_{k'}}\(x_{0,j},\sqrt{\lambda/\kappa_{\max}}\)$ with $x_{0,j} \in {\color{black}\bbM_{k\prime}}$.
Let also write $T_j := T_{x_{0,j}}^{(\Psi_{k'})} M_{k'}$ and $y_j := \Psi_{k'}(x_{0,j})$. Then from \lemref{geoineq} 5.,
\begin{align*}
\Psi_{k'}\(\ball_{\bbM_{k'}}({\color{black}x_{0,j}},\sqrt{\lambda/\kappa_{\max}})\) \subset {\color{black}\slab_{T_j}}\(y_j, \sqrt{\lambda/\kappa_{\max}}, \lambda\)
\end{align*}
for all $j \in \{1,\ldots,N'\}$.
Furthermore, using the fact that $\lambda \leq \sqrt{\lambda/\kappa_{\max}}$, we get
\begin{align*} 
M_k \cap M_{k'}^{\lambda} 
&\subset 
\bigcup_{j=1}^{N'} M_k \cap {\color{black}\slab_{T_j}}\(y_j, \sqrt{\lambda/\kappa_{\max}}, \lambda\)^\lambda  
\\
&\subset 
\bigcup_{j=1}^{N'} M_k \cap {\color{black}\slab_{T_j}}\(y_j, 2\sqrt{\lambda/\kappa_{\max}}, 2\lambda\).
\end{align*}
Applying point (ii) of \lemref{slab_vol2} to each element of the union (which we can, since $2\sqrt{\lambda/\kappa_{\max}} < {\color{black}1/(8\kappa_{\max})}$ by assumption), we get
\begin{align*}
\cH^{(d_k)}\(M_k \cap {\color{black}\slab_{T_j}}\(y_j, 2\sqrt{\lambda/\kappa_{\max}}, 2\lambda\)\) 
\leq 
C_{d_k} \nu_{\max} \kappa_{\max}^{d_k-d_{k'}} \(\frac{\lambda}{\kappa_{\max}}\)^{\frac{2d_k-d_{k'}}{2}}
,
\end{align*}
so that using an union bound over $j \in \{1,\ldots,N'\}$ we find
\begin{align*}
\cH^{(d_k)}(M_k \cap M_{k'}^{\lambda}) 
&\leq 
 (C_{d_{k'}}\nu_{\max} (\kappa_{\max} \lambda)^{-d_{k'}/2}) \times C_{d_k} \nu_{\max} \kappa_{\max}^{d_k-d_{k'}} \(\frac{\lambda}{\kappa_{\max}}\)^{\frac{2d_k-d_{k'}}{2}}
 \\
&=
C_{d_{k}} \nu_{\max}^2 \kappa_{\max}^{-d_{k'}}\lambda^{d_k -d_{k'}}
,
\end{align*}
which ends the proof.
\end{proof}

\begin{proof}[Proof of \thmref{main-clustering-dimensions}] 
As in the proof of Lemma~\ref{lem:recalgo}, let us write
\begin{align*}
\lambda_k := \xi_{d_k} \kappa_{\max} h_{d_k}^2 + \delta_{d_k}
,
\end{align*}where $\xi_{d_k}$ is given by~\eqref{def:xi_d}. 
We write $\cE := \cE_1 \cap \cE_2$, where $\cE_1$ and $\cE_2$ are the events of \lemref{slab_concentration} and \lemref{eventcE} respectively. 
In particular, we have $\bbP(\cE^c) \leq 6Kn^{-q}$, and \lemref{recalgo} yields that on the event on $\cE$,
$$
\cX_{k,n} \setminus \bigcup_{k' < k} M_{k'}^{\lambda_{k'}} \subset \wh \cX_{k,n} \subset \cX_{k,n} \cup \bigcup_{k' > k} \cX_{k',n} \cap M_k^{\lambda_k}.
$$
Hence,
\begin{align*}
\# \{\cX_{k,n} \Delta \wh \cX_{k,n}\}  \ind_{\cE}
&=
\# \{\cX_{k,n} \setminus \wh \cX_{k,n}\} \ind_{\cE}
+
\# \{\wh \cX_{k,n} \setminus \cX_{k,n}\} \ind_{\cE}
\\
&\leq \sum_{k' < k} \#\{\cX_{k,n} \cap M_{k'}^{\lambda_{k'}}\} \ind_{\cE} + \sum_{k' > k} \#\{\cX_{k',n} \cap M_{k}^{\lambda_{k}}\} \ind_{\cE}.
\end{align*}
Taking the expectation, we find that
\begin{align}
\label{eq:sum-cardinals}
\notag
\bbE&\left[\frac{\# \{\cX_{k,n} \Delta \wh \cX_{k,n}\}}{1 \vee N_k}  \ind_{\cE}\right] 
\\
&\leq \sum_{k' < k} \bbE\left[\frac{\#\{\cX_{k,n} \cap M_{k'}^{\lambda_{k'}}\}}{1\vee N_k} \ind_{\cE}\right] + \sum_{k' > k} \bbE\left[\frac{\#\{\cX_{k',n} \cap M_{k}^{\lambda_{k}}\}}{1\vee N_k} \ind_{\cE}\right].
\end{align}
For $k' < k$, it is easy to see that
\begin{align*} 
\bbE\left[\frac{\#\{\cX_{k,n} \cap M_{k'}^{\lambda_{k'}}\}}{1\vee N_k} \ind_{\cE}\right] 
&= 
\bbE\left[
\frac{1}{1 \vee N_k}
\bbE\left[\#\{\cX_{k,n} \cap M_{k'}^{\lambda_{k'}}\} \ind_{\cE} \middle| N_k\right]\right] 
\\
&\leq
\bbE\left[\frac{N_k}{1 \vee N_k}\right]
P_k(M_{k'}^{\lambda_{k'}})
\\
&\leq
C_{d_k}a_{\max}\nu_{\max}^2 \kappa_{\max}^{d_k-d_{k'}} \lambda_{k'}^{d_k-d_{k'}},
\end{align*} 
where we used \lemref{offset}. 
{\color{black}For $k'>k$, we use $\cE\subset G:=\{\min_jN_j/(n\alpha_j)>1/2\}$ and \lemref{offset} to get
\begin{align*}
\bbE\left[\frac{\#\{\cX_{k',n}\cap M_k^{\lambda_k}\}}{1\vee N_k}\ind_{\cE}\right]
&\leq\frac{2}{\alpha_kn}\bbE[\#\{\cX_{k',n}\cap M_k^{\lambda_k}\}]\\
&\leq\frac{\alpha_{k'}}{\alpha_k}C_{d_{k'}}a_{\max}\nu_{\max}^2
(\kappa_{\max}\lambda_k)^{d_{k'}-d_k}.
\end{align*}
With $\zeta_d=1/4$ as in~\eqref{eq:parameter-choice}, the constants in $\tau_j$ can be chosen so that
$\lambda_j\leq\tau_j(\log n/n)^{2/d_j}$.
The leading terms of the two finite sums in~\eqref{eq:sum-cardinals} correspond to $k'=k-1$ and $k'=k+1$ when those indices exist. For fixed parameters and $n$ large enough, the other terms are absorbed into these bounds, with the value of $\eta_k$ stated in the theorem.

It remains to control the loss outside $\cE$. Write $L_k:=\#(\cX_{k,n}\triangle\wh\cX_{k,n})/(1\vee N_k)$. Since the symmetric difference is contained in $\cX_n$, we have $L_k\leq2/\alpha_{\min}$ on $G$, and $L_k\leq n$ everywhere. Lemma~\ref{lem:slab_concentration} and the sharper bound in the proof of Lemma~\ref{lem:eventcE} therefore give
\begin{align*}
\bbE[L_k\ind_{\cE^c}]
&\leq\frac{2}{\alpha_{\min}}\{\bbP(\cE_1^c)+\bbP(\cE_2^c)\}+n\bbP(G^c)\\
&\leq\frac{8K}{\alpha_{\min}}n^{-\alpha_{\min}\gamma/16}
+\frac{4Ka_{\min}\nu_{\max}}{\alpha_{\min}\gamma\log n}
n^{1-3\alpha_{\min}\gamma/28}
+nK e^{-3\alpha_{\min}n/28}.
\end{align*}
Our choice of $\gamma$ ensures $\alpha_{\min}\gamma/16>q$ and $3\alpha_{\min}\gamma/28-1>q$: use $\alpha_{\min}\gamma\geq56q/3$ for $q\geq2$ and $\alpha_{\min}\gamma\geq800D^2\log D$ otherwise. Hence this bound is $o(n^{-q})$, and is at most $6Kn^{-q}$ for $n$ large enough. Combining the two parts proves the result.
}
\end{proof}

\section*{Acknowledgments}
We thank Alexandra Carpentier for making this work bloom.
EA would like to warmly thank Emmanuel Giroux for insightful discussions.
CB is grateful to Yann Chaubet for helpful discussions. 

\section*{Funding}
The work of EA is supported by the PEPR IA project PERSNET.
The work of CB was supported by the Deutsche Forschungsgemeinschaft (German Research Foundation) on the French-German PRCI ANR ASCAI CA 1488/4-1 ``Aktive und Batch-Segmentierung,
Clustering und Seriation: Grundlagen der KI''.

\bibliographystyle{imsart-number}
\bibliography{biblio}
\end{document}